\documentclass[11pt,twoside]{amsart}
\usepackage{amsmath, amsthm, amscd, amsfonts, amssymb, graphicx, color}
\usepackage[all,knot,arc,import,poly]{xy}

\usepackage[bookmarksnumbered, plainpages]{hyperref}
\newtheorem{theorem}{Theorem}[section]

\numberwithin{equation}{section}
\newtheorem{lemma}[theorem]{Lemma}
\newtheorem{corollary}[theorem]{Corollary}
\newtheorem{proposition}[theorem]{Proposition}

\def\pn{\par\noindent}

\begin{document}
%------------------------------------------------------------------------------------%

%------------------------------------------------------------------------------------%

\title{The $n$-ary adding machine and solvable groups }
\author{J. da Silva Rocha$^*$ and S. Najati Sidki }

\thanks{{\scriptsize
\hskip -0.4 true cm MSC(2010): Primary: 20F05; Secondary: 05C05; 20E08; 20F16.
$*$Corresponding author}}
\maketitle

\begin{abstract} We describe under various conditions abelian subgroups of the automorphism
group $\mathrm{Aut}(T_{n})$ of the regular $n$-ary tree $T_{n}$, which are
normalized by the $n$-ary adding machine $\tau =(e,...,e,\tau )\sigma _{\tau
}$ where $\sigma _{\tau }$ is the $n$-cycle $\left( 0,1,...,n-1\right) $. As
an application, for $n=p$ a prime number, and for $n=4$ , we prove that
every   soluble subgroup of $\mathrm{Aut}(T_{n})$,
containing $\tau $ is an extension of a torsion-free metabelian group by a
finite group.
\end{abstract}

\vskip 0.2 true cm

%------------------------------------------------------------------------------------%

\pagestyle{myheadings}
\markboth{\rightline {\sl  \hskip 10 cm  J. S. Rocha and S. N. Sidki }}
         {\leftline{\sl  \hskip 0 cm J. S. Rocha and S. N. Sidki}}

\bigskip
\bigskip

%------------------------------------------------------------------------------------%
%------------------------------------------------------------------------------------%

\section{\bf   Introduction}

\vskip 0.4 true cm

Adding machines have played an important role in dynamical systems, and in
the theory of groups acting on trees : see \cite{Bass95, Sid98, Sid00,
GNS00, Nek05}.

An element $\alpha $ in the automorphism group $\mathcal{A}_{n}=\mathrm{Aut}%
(T_{n})$ of the $n$-ary tree $T_{n}$, is represented as $\alpha =\alpha
|_{\phi }=\left( \alpha |_{0},...,\alpha |_{n-1}\right) \sigma _{\alpha }$
where $\phi $ is the empty sequence from the free monoid $\mathcal{M}$
generated by $Y=\left\{ 0,1,..,n-1\right\} $, where $\alpha |_{i}\in 
\mathcal{A}_{n},$ for $i\in Y,$ are called $1$st level states of $\alpha $
and where $\sigma _{\alpha }$ (the activity of $\alpha $) is a permutation
in the symmetric group $\Sigma _{n}$ on $Y$ extended `rigidly' to act on the
tree; if $\sigma _{\alpha }=e$, we say that $\alpha $ is inactive.

In applying the same representation to $\alpha |_{0}$ we produce $\alpha
|_{0i}$ for all $i\in Y$ and we produce in general $\left\{ \alpha |_{u}\mid
u\in \mathcal{M}\text{ }\right\} $ the set of \textit{states }of $\alpha $.
Following this notation, the $n$-ary adding machine is represented as $\tau
=(e,...,e,\tau )\sigma _{\tau }$ where $e$ is the identity automorphism and $%
\sigma _{\tau }$ is the regular permutation $\sigma =\left(
0,1,...,n-1\right) $. In this sense, the adding machine is an infinite
variant of the regular permutation which appears often in geometric and
combinatorial contexts.

A characteristic feature of $\tau $ is that its $n$-th power $\tau ^{n}$ is
the diagonal automorphism of the tree $\left( \tau ,...,\tau \right) $. This
fact implies that the centralizer of the cyclic group $\left\langle \tau
\right\rangle $ in $\mathcal{A}_{n}$ is equal to its topological closure $%
\overline{\langle\tau \rangle}$ in the group $\mathcal{A}_{n}$ when considered as a
topological group with respect to the the natural topology induced by the
tree. The pro-cyclic group $\overline{\langle\tau \rangle}$ is isomorphic to $\mathbb{Z}%
_{n}$, the ring of $n$-adic integers $\xi $ $=\sum_{i\geq 0}a_{i}n^{i}$ ($%
0\leq a_{i}\leq n-1$ for all $i$).

A large variety of subgroups of $\mathcal{A}_{n}$ which contain $\tau $ have
been constructed, including   groups which are torsion-free
and just non-solvable without free subgroups of rank $2$ (see, \cite{BSV99,
Sid-Sil01} and generalizations thereof \cite{Sid05}). Furthermore, the free
group of rank $2$ has been represented on the binary tree as a group
generated by two conjugates of the adding machine $\tau $ each having a
finite number of states \cite{Voro-Voro07}. On the other hand, the
restricted structure of its centralizer indicate that solvable groups which
contain $\tau $ have restricted structure. For nilpotent groups we show

\textbf{Proposition. }\textit{Let }$G$\textit{\ be a nilpotent subgroup of }$%
\mathcal{A}_{n}$\textit{\ which contains the }$n$\textit{-adic adding
machine }$\tau $\textit{. Then }$G$\textit{\ is a subgroup of }$\overline{%
<\tau >}$\textit{\ .}

The most visible examples of   solvable groups containing $%
\tau $ are conjugate to subgroups of those belonging to the infinite
sequence of groups 
\begin{eqnarray*}
\Gamma _{0} &=&N_{\mathcal{A}_{n}}(\overline{\langle \tau \rangle}),\text{ } \\
\Gamma _{i+1} &=&\left( \times _{n}\Gamma _{i}\right) \rtimes G_{i+1}\text{ }%
(i\geq 0)
\end{eqnarray*}%
where $\times _{n}\Gamma _{i}$ is a direct product of $n$ copies of $\Gamma
_{i}$ (seen as a subgroup of the $1$st level stabilizer of the tree) and
where $G_{i}$ is a solvable subgroup of the symmetric group $\Sigma _{n}$ in
its canonical action on the tree and containing the cycle $\sigma _{\tau }$.
We observe that for all $i$, the groups $\Gamma _{i}$ are metabelian by
'finite solvable subgroups of $\Sigma _{n}$'. It was shown by the second
author that for $n=2$,   solvable groups which contain the
binary adding machine are conjugate to some subgroups of $\Gamma _{i}$
acting on the binary tree \cite{Sid03}. \ \ This appears to be the general
pattern. However, the description for degrees $n>2$ requires a
classification of solvable subgroups of $\Sigma _{n}$ which contain the
cycle $\sigma =\left( 0,1,...,n-1\right) $\cite{Jones02}. This in itself is
an open problem, even for metabelian groups. On the other hand, the answer
for primitive solvable subgroups of $\Sigma _{n}$ is simple and classical.
For then, $n$ is a prime number $p$ or $n=4$. In case $n=p$, the solvable
subgroups $G_{i}$ can all be taken to be the normalizer $F=N_{\Sigma
_{n}}\left( \left\langle \sigma \right\rangle \right) $ of order $p\left(
p-1\right) $ and in case $n=4$, the $G_{i}$'s can all be taken to be the
symmetric group $\Sigma _{4}$.

Given this background, the main theorem of this paper is

\textbf{Theorem A. }\textit{Let }$n=p$\textit{, a prime number, or }$n=4$%
\textit{. Then any   solvable subgroup of }$\mathcal{A}_{n}$%
\textit{which contains the }$n$\textit{-ary machine }$\tau $\textit{\ is
conjugate to a subgroup of }$\Gamma _{i}$\textit{\ for some }$i$.

The result follows first from general analysis, first of the conditions $%
[\beta ,\beta ^{\tau ^{x}}]=e$ (for some $\beta \in \mathcal{A}_{n}$ and all 
$x\in \mathbb{Z)}$, then their impact on the $1$st level states of the
subgroup $\left\langle \beta ,\tau \right\rangle $ and on how these in turn
translate successively to conditions on states at lower levels. It is
somewhat surprising that the process converges to a clear global description
for trees of degrees $p$ and $4$.

The first step of this analysis lead to the following description of the
normal closure of $\left\langle \beta \right\rangle $ under the action of $%
\tau $.

\textbf{Theorem B. } {\em Let $B$ be an abelian subgroup of $%
\mathcal{A}_{n}$ normalized by $\tau,$  let $\beta
=(\beta |_{0},\beta |_{1},\ldots ,\beta |_{n-1})\sigma _{\beta }\in B $%
and define the subgroup $H=\left\langle \beta |_{i}\left(
i\in Y\right) ,\tau \right\rangle $ generated by the first level
states of $\beta $ and $\tau. $ \newline
(I) Suppose $\sigma _{\beta }=\left( \sigma _{\tau }\right) ^{s}$
for some integer $s.$  Then $H$ is metabelian-by-finite.
More precisely, let $\;m=\frac{n}{\gcd (n,s)}$, define the product 
$\pi _{i}=\beta |_{i}\beta |_{ i+s}\beta |_{ i+2s}\cdots \beta |_{ i+\left( m-1\right) s}$%
 (the notation $\beta|_{j}$ means $\beta|_{\overline{j}}, $ where  $\overline{j}$ is the representative of $j$
in $Y$ modulo $n$) and define the subgroup  
\begin{equation*}
K=\left\langle 
\begin{array}{c}
\lbrack \beta |_{i},\tau ^{k}],\;\pi _{i}\mid k\in \mathbb{Z},\;i\in Y%
\end{array}%
\right\rangle
\end{equation*}%
Then $K$ is an abelian group and $H$ affords the normal series 
\begin{equation*}
H\trianglerighteq K\left\langle \tau \right\rangle \left( =O\right)
\trianglerighteq K
\end{equation*}%
where the quotient group $\frac{H}{O}$  is a homomorphic image of a
subgroup of the wreath product $C_{m}\wr C_{n}$ of the cyclic
groups $C_{m},C_{n}$.\newline
(II) Let $n$  be an even number. Then $H$ is a metabelian
group if $s=\frac{n}{2}$  or if $\sigma _{\beta }$ is a
transposition. } 

Part (I) of Theorem B will be proven in Sections 4 and 5 and part (II) in
Section 7.

Let $P$ be a subgroup of $\Sigma _{n}$. The \textit{layer closure} of $P$ in 
$\mathcal{A}_{n}$ is the group $L\left( P\right) $ formed by elements of $%
\mathcal{A}_{n}$ whose states have activities in $P$. The following result
is yet another characterization of the adding machine.

\textbf{Theorem C. }\textit{Let }$n$\textit{\ be an odd number} \textit{\
and let }$L=L\left( \left\langle \sigma \right\rangle \right) $,\textit{\
the layer closure of }$\left\langle \sigma \right\rangle $\textit{\ in }$%
A_{n}$\textit{.\ Let }$s$\textit{\ be an integer which is relatively prime
to }$n$\textit{\ and let} $\beta =(\beta |_{0},\beta |_{1},\ldots ,\beta
|_{n-1})\sigma ^{s}\in L$\textit{\ be such }that\textit{\ }$[\beta ,\beta
^{\tau ^{x}}]=e$\textit{\ for all }$x\in \mathbb{Z}.$\textit{\ Then }$\beta $%
\textit{\ is a conjugate of }$\tau $\textit{\ in }$L$\textit{.}

\section{\bf Preliminaries}

\vskip 0.4 true cm

We start by introducing definitions and notation. The $n$-ary tree $T_{n}$
can be identified with the free monoid $\mathcal{M}=<$ $0,1,..,n-1>^{\ast }$
of finite sequences from $Y=$ $\left\{ 0,1,...,n-1\right\} $, ordered by $%
v\leq u$ provided $u$ is an initial subword of $v$.

The identity element of $\mathcal{M}$ is the empty sequence $\phi .$ The
level function for $T_{n},$ denoted by $|m|$ is the length of $m\in \mathcal{%
M}$; the root vertex $\phi $ has level $0$.

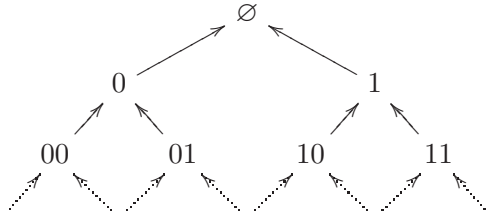
\begin{figure}[htbp]
\begin{center}
\xymatrix @R=.4cm @C=.2cm{&&&&&&&&&&&\varnothing&&& \\
&&&&&&& &&0\ar[urr] & & & & 1\ar[ull]& \\
&&&&&&&&00\ar[ur] & & 01\ar[ul] & & 10\ar[ur] & & 11\ar[ul]& \\
&&&&&&&\ar@{.>}[ur] &&\ar@{.>}[ul] \;\;\; \ar@{.>}[ur]&&\ar@{.>}[ul] \;\;\; %
\ar@{.>}[ur]&&\ar@{.>}[ul] \;\;\; \ar@{.>}[ur]&& \ar@{.>}[ul] }
\end{center}
\caption{The Binary Tree}
\end{figure}

The action $\rho :i\rightarrow j$ of a permutation $\rho \in \Sigma _{n}$
will be from the right and written as $\left( i\right) \rho =j$ or as $%
i^{\rho }=j$. If $i,j$ are integers then the action of $\rho $ on $i$ is to
be identified with its action on its representatives $\overline{i}$ in $Y$,
modulo $n$ . Permutations $\sigma $ in $\Sigma _{n}$ are extended `rigidly'
to automorphisms of $\mathcal{A}_{n}$ by 
\begin{equation*}
(y.u)\rho =(y)\sigma .u,\;\forall \;y\in Y,\;\forall \;u\in \mathcal{M}\text{%
.}
\end{equation*}

An automorphism $\alpha \in \mathcal{A}_{n}$ induces a permutation $\sigma
_{\alpha }$ on the set $Y$. Consequently, $\alpha $ affords the
representation $\alpha =\alpha ^{\prime }\sigma _{\alpha }$ where $\alpha
^{\prime }$ fixes $Y$ point-wise and for each $i\in Y$, $\alpha ^{\prime }$
induces $\alpha |_{i}$ on the subtree whose vertices form the set $i\cdot 
\mathcal{M}$. If $j$ is an integer the $\alpha |_{j}$ will be understood as $%
\alpha |_{\overline{j}}$ where $\overline{j}$ is the representative of $j$
in $Y$ modulo $n$.

Given $i$ in $Y$, we use the canonical isomorphism $i\cdot {}u\mapsto u$
between $i\cdot \mathcal{M}$ and the tree $T_{n}$, and thus identify $\alpha
|_{i}$ with an automorphism of $T_{n}$; therefore, $\alpha ^{\prime }\in 
\mathcal{F}(Y,\mathcal{A}_{n})$, the set of functions from $Y$ into $%
\mathcal{A}_{n}$, or what is the same, the $1$st level stabilizer $Stab(1)$
of the tree. This provides us with the factorization $\mathcal{A}_{n}=%
\mathcal{F}(Y,\mathcal{A}_{n})\cdot \Sigma _{n}$.

Let $\alpha ,\beta ,\gamma \in \mathcal{A}_{n}$. Then the following
formulas hold

\begin{equation}
\sigma _{\alpha ^{-1}}=\left( \sigma _{\alpha }\right) ^{-1},\text{ }\sigma
_{\alpha }\sigma _{\beta }=\sigma _{\alpha \beta }\text{,}  \label{eq1}
\end{equation}%
\begin{equation}
(\alpha ^{-1})|_{u}=\left(\alpha |_{\left( u\right) ^{{\alpha }%
^{-1}}}\right)^{-1}\text{,}  \label{eq2}
\end{equation}%
\begin{equation}
\text{ }(\alpha \beta )|_{u}=\left( \alpha |_{u}\right) \left( \gamma
|_{u}\right) \text{ where }\gamma |_{u}=\beta |_{(u)^{\alpha }}  \label{eq3}
\end{equation}%
\begin{equation}\label{eq4}
\gamma =\alpha ^{-1}\beta \alpha \Leftrightarrow \left( \sigma _{\gamma }=\sigma
_{\alpha }^{-1}\sigma _{\beta }\sigma _{\alpha } \text{ and }
\gamma |_{\left( i\right) \sigma _{\alpha }}=\alpha |_{i}^{-1}\beta
|_{i}\alpha |_{\left( i\right) \sigma _{\beta }},\forall i\in Y\right) \text{.}
\end{equation}%
\begin{equation} \label{eq6}
\theta =[\beta ,\alpha ]=\beta ^{-1}\beta ^{\alpha }\Leftrightarrow \left\{ \begin{array}{l} \sigma
_{\theta }=[\sigma _{\beta },\sigma _{\alpha }]\text{, }  
\\ 
\theta |_{\left( i\right) \sigma _{\alpha \beta }}=\left( \beta |_{\left(
i\right) \sigma _{\alpha }}\right) ^{-1}\left( \alpha |_{i}\right)
^{-1}\left( \beta |_{i}\right) \left( \alpha |_{\left( i\right) \sigma
_{\beta }}\right) ,\forall i\in Y\text{.}  
\end{array} \right.
\end{equation}

\begin{equation}
\left( \alpha ^{m}\right) |_{i}=\left( \alpha |_{i}\right) \left( \alpha
|_{(i){\sigma _{\alpha }}}\right) \left( \alpha |_{(i){\sigma _{\alpha }^{2}}%
}\right) \cdots \left( \alpha |_{(i){\sigma _{\alpha }^{m-1}}}\right)
\label{eq8}
\end{equation}

\begin{equation}
\left( \beta ^{\alpha }\right) |_{u}=\left( \beta |_{\left( u\right) \alpha
^{-1}}\right) ^{\alpha |_{\left( u\right) \alpha ^{-1}}},\text{where }\beta
\in Stab(k)\text{ and }|u|\leq k.  \label{eq9}
\end{equation}

An automorphism $\alpha \in \mathcal{A}_{n}$ corresponds to an input-output
automaton with alphabet $Y$ and with set of states $\mathrm{Q}(\alpha
)=\{\alpha |_{u}\mid u\in \mathcal{M}\}$. The automaton $\alpha $ transforms
the letters as follows: if the automaton is in state $\alpha |_{u}$ and
reads a letter $i\in Y$ then it outputs the letter $j=\left( i\right) \alpha
|_{u}$ and its state changes to $\alpha |_{ui}$; these operations can be
best described by the labeled edge $\alpha |_{u}\overset{i|j}{%
\longrightarrow }\alpha |_{ui}$. Following terminology of automata theory,
every automorphism $\alpha |_{u}$ is called the \textit{state} of $\alpha $
at $u$.

The tree $T_{n}\ $is a topological space which is the direct limit of its
truncations at the $n$-th levels. Thus the group $\mathcal{A}_{n}$ is the
inverse limit of the permutation groups it induces on the $n$-th level
vertices. This transforms $\mathcal{A}_{n}$ into a topological group. An
infinite product of elements $\mathcal{A}_{n}$ is a well-defined element of $%
\mathcal{A}_{n}$ provided that for any given level $l$, only finitely many
of the elements in the product have non-trivial action on vertices at level $%
l$. Thus, if $\alpha \in \mathcal{A}_{n}$ and $\xi $ $=\sum_{i\geq
0}a_{i}n^{i}\in \mathbb{Z}_{n\text{ }}$ then $\alpha ^{\xi }=\alpha
^{a_{0}}\cdot \alpha ^{na_{1}}\cdots \alpha ^{n^{i}a_{i}}\cdots $ is a well
defined element of $\mathcal{A}_{n}$. The notation $\alpha |_{\xi }$ is to
be understood as $\alpha |_{i}$ where $i=a_{0}$.

The topological closure of a subgroup $H$ in $\mathcal{A}_{n}$ \ will be
indicated by $\overline{H}$. We note that if $H$ is abelian then 
\begin{equation*}
\overline{H}=\{h^{\xi }\mid h\in H,\xi \text{ }\in \mathbb{Z}_{n\text{ }}\}%
\text{.}
\end{equation*}%
One of the characterizing aspects of the $n$-ary adding machine is that the
centralizer of $\tau $ is a pro-cyclic group; namely, 
\begin{equation*}
C_{\mathcal{A}_{n}}(\tau )=\overline{\left\langle \tau \right\rangle }%
=\{\tau ^{\xi }\mid \xi \in \mathbb{Z}_{n}\}\text{.}
\end{equation*}

Let $v=yu$ where $y\in Y,u\in \mathcal{M}$. The image of $v$ under the
action of $\alpha $ is 
\begin{equation*}
(v)\alpha =(yu)\alpha =\left( y\right) \sigma _{\alpha }.(u)\alpha |_{y}%
\text{.}
\end{equation*}%
The action extends to infinite sequences (or boundary points of the tree) in
the same manner. A boundary point of the tree $c=c_{0}c_{1}c_{2}\ldots, $%
where $c_{i}\in Y$ for all $i$, corresponds also to the $n$-adic integer $%
\xi =\sum \{c_{i}n^{i}|i\geq 0\}\in \mathbb{Z}_{n}$. Thus the action of the
tree automorphism $\alpha $ can thus be translated to an action on the ring
of $n$-adic integers. We will indicate $c_{0}$ by $\overline{\xi }$ which is 
$\xi $ modulo $n$. In the case of the automorphism $\tau =(e,e,...,e,\tau
)\sigma $, the action of $\tau $ on $c$ is 
\begin{equation*}
\left( c\right) \tau =\left\{ 
\begin{array}{ll}
\left( c_{0}+1\right) c_{1}c_{2}\ldots & \text{ if }0\leq c_{0}\leq n-2, \\ 
0(c_{1}c_{2}\ldots )\tau , & \text{ if }c_{0}=n-1,%
\end{array}%
\right.
\end{equation*}%
which translates to the $n$-ary addition 
\begin{equation*}
\xi ^{\tau } = \xi + 1 \text{.}
\end{equation*}

\begin{figure}[htbp]
\begin{center}
\xymatrix @C=1.5cm @R=2cm{&& *++[o][F-]{\tau }\ar@(rd,ld)^{1/0} \ar[rr]%
^{0/1} & & *++[o][F-]{e}\ar@(ru,lu)_{0/0,\; 1/1} }
\end{center}
\caption{The binary adding machine}
\end{figure}
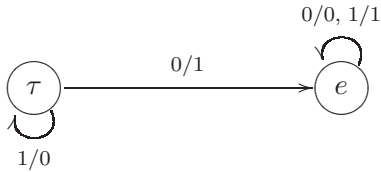

\section{\bf  Normalizer of the topological closure $\overline{\left\langle 
\protect\tau \right\rangle }$}

\vskip 0.4 true cm

An element $\xi =\sum_{i\geq 0}a_{i}n^{i}\in \mathbb{Z}_{n}$ is a unit in $%
\mathbb{Z}_{n}$ if and only if $\overline{\xi }\left( =a_{0}\right) $ is a
unit in $\mathbb{Z}$ modulo $n$. The group of automorphisms of $\mathbb{Z}%
_{n}$ is isomorphic to the multiplicative group of units $U(\mathbb{Z}_{n})$%
. The subgroup of $U(\mathbb{Z}_{n})$ consisting of elements $\xi $ with $%
\overline{\xi }=1$ is denoted by by $\mathbb{Z}_{n}^{1}$. This subgroup has
the transversal $\left\{ j\mid 1\leq j\leq n-1,\gcd \left( j,n\right)
=1\right\} $ in $U(\mathbb{Z}_{n})$ and therefore has index$\ \left[ U(\mathbb{Z%
}_{n}):\mathbb{Z}_{n}^{1}\right] =\varphi \left( n\right) $ where $\varphi $
is the Euler function.

Given $\alpha \in \mathcal{A}_{n}$ we denote the diagonal automorphism $%
\left( \alpha ,...,\alpha \right) $ by $\alpha ^{\left( 1\right) }$ and
define inductively $\alpha ^{\left( i+1\right) }=\left( \alpha ^{\left(
i\right) }\right) ^{\left( 1\right) }$ for all $i\geq 1$.

\subsection{Powers of $\protect\tau $}

Let $\xi =\sum_{i\geq 0}a_{i}n^{i}\in \mathbb{Z}_{n}$. Then $\sum_{i\geq
1}a_{i}n^{i-1}=\frac{\xi -\overline{\xi }}{n}$.

\begin{lemma}
Let $\xi \in \mathbb{Z}_{n}$. Then 
\begin{equation*}
\tau ^{\xi }=(\tau ^{\frac{\xi -a_{0}}{n}},\ldots ,\tau ^{\frac{\xi -a_{0}}{n%
}},\underbrace{\tau ^{\frac{\xi -a_{0}}{n}+1},\ldots ,\tau ^{\frac{\xi -a_{0}%
}{n}+1}}_{a_{0}\;\mathrm{terms}})\sigma _{\tau }^{a_{0}}.
\end{equation*}
\end{lemma}

\begin{proof}
For $j$ an integer with $1\leq j\leq n-1$, we have 
\begin{equation*}
\tau ^{j}=\left( e,...,e,\underbrace{\tau ,\ldots ,\tau }_{j\;\mathrm{terms}%
}\right) \sigma _{\tau }^{j}
\end{equation*}%
and $\tau ^{n}=\left( \tau ,...,\tau \right) =\tau ^{\left( 1\right) }$.

Given $\xi =\sum_{i\geq 0}a_{i}n^{i}$, then

\begin{eqnarray}
\tau ^{a_{0}} &=&(e,\ldots ,e,\underbrace{\tau ,\ldots ,\tau }_{a_{0}\;%
\mathrm{terms}})\sigma _{\tau }^{a_{0}},  \label{i1} \\
\tau ^{a_{j}n^{j}} &=&\tau ^{(a_{j}n^{j-1})n}=\left( \tau
^{a_{j}n^{j-1}}\right) ^{\left( 1\right) }, \\
\tau ^{\xi } &=&(\tau ^{\frac{\xi -a_{0}}{n}},\ldots ,\tau ^{\frac{\xi -a_{0}%
}{n}},\underbrace{\tau ^{\frac{\xi -a_{0}}{n}+1},\ldots ,\tau ^{\frac{\xi
-a_{0}}{n}+1}}_{a_{0}\;\mathrm{terms}})\sigma _{\tau }^{a_{0}} \\
&=&(\tau ^{\frac{\xi -\overline{\xi }}{n}},\ldots ,\tau ^{\frac{\xi -%
\overline{\xi }}{n}},\underbrace{\tau ^{\frac{\xi -\overline{\xi }}{n}%
+1},\ldots ,\tau ^{\frac{\xi -\overline{\xi }}{n}+1}}_{\overline{\xi }\;%
\mathrm{terms}})\sigma _{\tau }^{\overline{\xi }}\text{.}
\end{eqnarray}
\end{proof}

As we have seen, the description of $\tau ^{\xi }$ involves the partition of
the interval $[0,...,n-1]$ into two subintervals. It is convenient to use
here the carry 2-cocycle

\bigskip\ $\delta :\mathbb{Z}_{n}\times \mathbb{Z}_{n}\rightarrow \{0,1\}$%
\textit{\ }defined by%
\begin{equation*}
\delta (\eta ,\kappa )=\frac{\overline{\eta }+\overline{\kappa }-\overline{%
\eta +\kappa }}{n}=\left\{ 
\begin{array}{ll}
0, & \text{ if }\overline{\eta }+\overline{\kappa }<n \\ 
1, & \text{otherwise}%
\end{array}%
\right. \text{.}
\end{equation*}%
We call this 2-valued function by \textit{Delta-2 ( }later on we will
introduce a 3-valued function \textit{Delta-3)}.

Using \textit{Delta-2}, the notation for the power of $\tau $ becomes 
\begin{equation}
\tau ^{\xi }=\left( \tau ^{\frac{\xi -\overline{\xi }}{n}+\delta (i,\xi
)}\right) _{0\leq i\leq n-1}\sigma _{\tau }^{\overline{\xi }}\text{.}
\label{funtal}
\end{equation}

\begin{figure}[h]
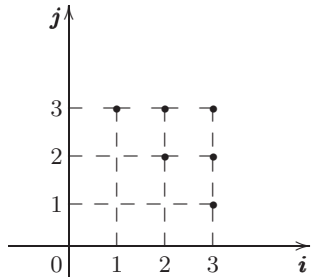

\label{grafx} \hspace{2cm} {\footnotesize \xy 
/r.8mm/:{(0,0)*{}; (50,0)*{} **\dir{-} ?>*\dir{>}; (10, -5)*{}; (10,40)*{} **%
\dir{-} ?>* \dir{>}; (8,-3)*{0};(18,-3)*{1};(26,-3)*{2};(34,-3)*{3};(49,-3)*{%
\pmb{i}}; (8,7)*{1};(8,15)*{2};(8,23)*{3};(8,38)*{\pmb{j}};
(18,23)*{};(18,0)*{};**\dir{--} ?>*\dir{*}; (18,23)*{};(10,23)*{};**\dir{--}%
; (26,23)*{};(26,0)*{};**\dir{--} ?>*\dir{*}; (26,23)*{};(18,23)*{};**%
\dir{--}; (34,23)*{};(34,0)*{};**\dir{--} ?>*\dir{*};
(34,23)*{};(26,23)*{};**\dir{--}; (26,15)*{};(10,15)*{};**\dir{--} ?>*\dir{*}%
; (34,15)*{};(26,15)*{};**\dir{--} ?>*\dir{*}; (34,7)*{};(10,7)*{};**\dir{--}
?>*\dir{*};} \endxy
}
\caption{Delta 2 function for $n = 4.$}
\end{figure}

\subsection{Centralizer of $\protect\tau $}

\begin{lemma}\label{centralizadortau}
$C_{\mathcal{A}_{n}}\left( \tau \right) =\overline{\left\langle \tau
\right\rangle }$.
\end{lemma}

\begin{proof}
Let $\alpha \in \mathcal{A}_{n}$ commute with $\tau $. Then, $[\sigma
_{\alpha },\sigma _{\tau }]=e$ and therefore $\sigma _{\alpha }=\left(
\sigma _{\tau }\right) ^{s_{0}}$ for some integer $0\leq s_{0}\leq n-1$.
Therefore, $\beta =\alpha \tau ^{-s_{0}}=\left( \beta |_{0},...,\beta
|_{n-1}\right) $ commutes with $\tau $ and $\sigma _{\beta }=e$. Now, 
\begin{equation*}
\beta ^{\tau }=\left( \left( \beta |_{n-1}\right) ^{\tau },\beta
|_{0},...,\beta |_{n-2}\right) =\beta
\end{equation*}%
implies $\beta |_{i}=\beta |_{0}\ $ for all $0\leq i \leq n-1$ and $\beta
|_{0}$ commutes with $\tau $. Therefore $\beta =\left( \beta |_{0}\right)
^{\left( 1\right) }$ and $\beta |_{0}$ replaces $\alpha $ in the previous
argument. Hence, there exists an integer $s_{1}$ such that $0\leq s_{1}\leq
n-1$ and $\gamma =\beta |_{0}\tau ^{-s_{1}}=\left( \gamma |_{0}\right)
^{\left( 1\right) }$. From this we conclude 
\begin{eqnarray*}
\alpha &=&\beta \tau ^{s_{0}}=\left( \beta |_{0}\right) ^{\left( 1\right)
}\tau ^{s_{0}} \\
&=&\left( \left( \gamma |_{0}\right) ^{\left( 1\right) }\tau
^{s_{1}},\ldots,\left( \gamma |_{0}\right) ^{\left( 1\right) }\tau
^{s_{1}}\right) \tau ^{s_{0}} \\
&=&\left( \gamma |_{0}\right) ^{\left( 2\right) }\tau ^{ns_{1}}\tau
^{s_{0}}=\left( \gamma |_{0}\right) ^{\left( 2\right) }\tau
^{ns_{1}+}{}^{s_{0}}\text{.}
\end{eqnarray*}%
We obtain the desired form inductively, $\alpha =\tau ^{\xi }$ where 
\begin{equation*}
\xi =s_{0}+s_{1}n+ s_{2}n^{2} + \ldots
\end{equation*}
\end{proof}

The characterization of nilpotent groups which contain $\tau $, announced in
the introduction, follows.

\begin{proposition}
\textit{Let }$G$\textit{\ be a nilpotent subgroup of }$\mathcal{A}_{n}$%
\textit{\ which contains the }$n$\textit{-adic adding machine }$\tau $%
\textit{. Then }$G$\textit{\ is a subgroup of }$\overline{<\tau >}$\textit{\
.}
\end{proposition}

\begin{proof}
Suppose $G$ is a nilpotent group of class $k>1$ which contains $\tau $.
Then, the center $Z\left( G\right) $ is contained in $\overline{\left\langle
\tau \right\rangle }$. Let $j$ be the maximum index such that $Z_{j}\left(
G\right) \leq \overline{\left\langle \tau \right\rangle }$ and let $\alpha
\in Z_{j+1}\left( G\right) \backslash\overline{\left\langle \tau \right\rangle }$. Then $\left[
\tau ,\alpha \right] $ $\in Z_{j}\left( G\right) $ and therefore $\left[
\tau ,\alpha \right] =\tau ^{\xi }$ for some $\xi \in \mathbb{Z}%
_{n}\backslash \left\{ 0\right\} $. Therefore 
\begin{eqnarray*}
\left[ \tau ,2\alpha \right] &=&\left[ \tau ,\alpha ,\alpha \right] =\left[
\tau ^{\xi },\alpha \right] \\
&=&\left[ \tau ,\alpha \right] ^{\xi }=\tau ^{\xi ^{2}}\in Z_{j-1}\left(
G\right)
\end{eqnarray*}%
and more generally, for $l\geq 1$, we have $\left[ \tau ,l\alpha \right]
=\tau ^{\xi ^{l}}\in Z_{j-l+1}\left( G\right) $. It follows that $\tau ^{\xi
^{j-1}}\in Z_{0}\left( G\right) =\left\{ e\right\} $. Thus, $\xi ^{j-1}=0$
and $\xi =0$; a contradiction.
\end{proof}

\subsection{ Normalizer of $\overline{\left\langle \protect\tau %
\right\rangle }$}

\begin{lemma}
The group $\Gamma _{0}=N_{\mathcal{A}_{n}}\left( \overline{\left\langle \tau
\right\rangle }\right) $ is metabelian. Indeed, the derived subgroup $\Gamma
_{0}^{\prime }$ is contained in $\overline{\left\langle \tau \right\rangle }$%
.
\end{lemma}

\begin{proof}
Let $\alpha ,\beta \in \Gamma _{0},$ then $\tau ^{\alpha }=\tau ^{\xi }$ and 
$\tau ^{\beta }=\tau ^{\eta }$ for some $\eta ,\xi \in U(\mathbb{Z}_{n})$.
Therefore,

\begin{equation*}
\tau ^{\alpha }=\tau ^{\xi },\tau =(\tau ^{\xi })^{\alpha ^{-1}}=(\tau
^{\alpha ^{-1}})^{\xi },
\end{equation*}%
\begin{equation*}
\tau ^{\alpha ^{-1}}=\tau ^{\xi ^{-1}}\text{.}
\end{equation*}

Likewise, $\tau ^{\beta ^{-1}}=\tau ^{\eta ^{-1}}$. Thus, $\tau ^{\lbrack
\alpha ,\beta ]}=\tau $ and $\Gamma _{0}^{\prime }\leq C_{\mathcal{A}%
_{n}}(\tau )=\overline{\left\langle \tau \right\rangle }$ follows.
\end{proof}

We present properties of the Delta-2 function which we will use in the
sequel.

\begin{lemma}
For all $0\leq i,j<n$ and $\xi \in \mathbb{Z}_{n}$ we have 
\begin{eqnarray*}
\sum_{i=0}^{n-1}\delta (i,j) &=&\overline{j}, \\
\delta (i,j\xi ) &=&j\left( \frac{\xi -\overline{\xi }}{n}\right) -\frac{%
j\xi -\overline{j\xi }}{n}+\sum_{k=0}^{j-1}\delta (i+k\xi ,\xi )\text{.}
\end{eqnarray*}
\end{lemma}

\begin{proof}
The first assertion is easy to verify.

The second is obtained from 
\begin{equation*}
(\tau ^{\xi })^{j}|_{i}=(\tau ^{\xi })|_{i}(\tau ^{\xi })|_{i+\xi }\cdots
(\tau ^{\xi })|_{i+(j-1)\xi },
\end{equation*}%
by substituting 
\begin{equation*}
(\tau ^{\xi })|_{i}=\tau ^{\frac{\xi -\overline{\xi }}{n}+\delta (i,\xi )}
\end{equation*}%
in its right hand side and 
\begin{equation*}
\tau ^{\xi j}|_{i}=\tau ^{\frac{j\xi -\overline{j\xi }}{n}+\delta (i,j\xi )}
\end{equation*}%
in its left.
\end{proof}

\begin{proposition}
\label{conjtau} Given $\alpha \in \mathcal{A}_{n}$ and $\xi \in U(\mathbb{Z}%
_{n})$. Then the condition $\tau ^{\alpha }=\tau ^{\xi }$ is equivalent to
conditions (i), (ii) and (iii) below.\newline
(i)%
\begin{equation*}
\alpha |_{i}=\left( \alpha |_{0}\right) \tau ^{\mu _{i}}\text{ }\left( 1\leq
i\leq n-1\right)
\end{equation*}%
where 
\begin{equation*}
\mu _{i}=i\frac{(\xi -\overline{\xi })}{n}+\sum_{k=0}^{i-1}\delta ((v(\alpha
)+k)\xi ,\xi )
\end{equation*}%
and where $v(\alpha )$ is defined by%
\begin{eqnarray*}
0 &\leq &v(\alpha )\leq n-1, \\
\left( 0\right) \sigma _{\alpha } &=&\overline{v(\alpha )\xi }\text{;}
\end{eqnarray*}%
(ii) (recursion)%
\begin{equation*}
\tau ^{\alpha |_{0}}=\tau ^{\xi }\text{;}
\end{equation*}%
(iii)%
\begin{equation*}
(j)\sigma _{\alpha }=\overline{(v(\alpha )+j)\xi }\text{ }\left( 0\leq j\leq
n-1\right) \text{.}
\end{equation*}%
Furthermore, if $\xi \in \mathbb{Z}_{n}^{1}$, then $v(\alpha )=0,\;(j)\sigma
_{\alpha }=\overline{j\xi }=j$ and $\mu _{i}=i\frac{\xi -1}{n}.$
\end{proposition}

\begin{proof}
Since $\sigma _{\tau }^{\sigma _{\alpha }}=\sigma _{\tau }^{\xi }$, we have
an equality between the permutations 
\begin{equation*}
(\left( 0\right) \sigma _{\alpha },\left( 1\right) \sigma _{\alpha },\ldots
,(n-1)\sigma _{\alpha })=(0,\overline{\xi },\overline{2\xi },\ldots ,%
\overline{(n-1)\xi })\text{.}
\end{equation*}%
Therefore, there exists $v(\alpha )\in Y$ such that $\left( 0\right) \sigma
_{\alpha }=\overline{v(\alpha )\xi }\ $and so, 
\begin{equation*}
(j)\sigma _{\alpha }=\overline{(v(\alpha )+j)\xi },\;\forall j\in Y\text{.}
\end{equation*}

Now, $\tau ^{\alpha }=\tau ^{\xi }$ is equivalent to $\alpha =\tau
^{-s}\alpha \tau ^{s\xi }$ for every $s\in \mathbb{Z},$ which in turn is
equivalent to 
\begin{equation*}
\alpha |_{\left( i\right) \sigma _{\tau }^{s}}=\left( (\tau
^{s})|_{i}\right) ^{-1}\left( \alpha |_{i}\right) (\tau ^{\xi {s}})|_{\left(
i\right) \sigma _{\alpha }},\forall i\in Y,\forall s\in \mathbb{Z}\text{.}
\end{equation*}

The latter conditions are equivalent to%
\begin{eqnarray*}
\alpha |_{0} &=&\alpha |_{\left( 0\right) \sigma _{\tau }^{n}}=\left( (\tau
^{n})|_{0}\right) ^{-1}\left( \alpha |_{0}\right) (\tau ^{\xi {n}})|_{\left(
0\right) \sigma _{\alpha }}, \\
\alpha |_{i} &=&\alpha |_{\left( 0\right) \sigma _{\tau }^{i}}=\left( (\tau
^{i})|_{0}\right) ^{-1}\alpha |_{0}(\tau ^{\xi {i}})|_{\left( 0\right)
\sigma _{\alpha }}\text{ }\forall i\in Y\backslash \{0\}
\end{eqnarray*}%
and these in turn are equivalent to%
\begin{eqnarray*}
\text{ }\alpha |_{i} &=&\alpha |_{0}\tau ^{\frac{\xi {i}-\overline{\xi {i}}}{%
n}+\delta (v(\alpha )\xi ,\xi {i})}=\alpha |_{0}\tau ^{\mu _{i}}, \\
\mu _{i} &=&i\left( \frac{\xi -\overline{\xi }}{n}\right)
+\sum_{k=0}^{i-1}\delta ((v(\alpha )+k)\xi ,\xi )\text{ }\forall i\in
Y\backslash \{0\}\text{.}
\end{eqnarray*}

If $\xi \in \mathbb{Z}_{n}^{1},$ then $\sum_{k=0}^{i-1}\delta (k\xi ,\xi
)=\sum_{k=0}^{i-1}\delta (k,1)=0$. The rest of the assertion follows
directly.
\end{proof}

\begin{corollary}
Let $\xi \in U\left( \mathbb{Z}_{n}\right), \sigma_{\alpha}$  and $\mu _{i}$ be as above.
Then $\alpha =\left( \alpha \right) ^{\left( 1\right) }\left( e,\tau ^{\mu
_{1}},...,\tau ^{\mu _{n-1}}\right)\sigma_{\alpha} $  conjugates $\tau $ to $\tau ^{\xi }.$ 
In particular,  if 
$\xi \in \mathbb{Z}_{n}^{1}, $ then $
\alpha =\left( \alpha \right) ^{\left( 1\right) }(e,\tau ^{\frac{\xi -1}{n}%
},\tau ^{2\frac{\xi -1}{n}},\ldots ,\tau ^{(n-1)\frac{\xi -1}{n}})
$ 
 (denoted by $\lambda _{\xi }$) conjugates $\tau $ to $\tau
^{\xi }$.
\end{corollary}

Although we have computed above an automorphism which inverts $\tau $, we
give another with a simpler description. Define the permutation 
\begin{equation*}
\varepsilon =\left( 0,n-1\right) \left( 1,n-2\right) ...\left( \left[ \frac{%
n-2}{2}\right] ,\left[ \frac{n+1}{2}\right] \right) \text{.}
\end{equation*}%
Then $\varepsilon $ inverts $\sigma _{\tau }=\left( 0,1,...,n-1\right) $ and 
\begin{equation*}
\iota =\iota ^{\left( 1\right) }\varepsilon
\end{equation*}%
inverts $\tau $.

Define 
\begin{eqnarray*}
\Lambda &=&\{\lambda _{\xi }\mid \xi \in \mathbb{Z}_{n}^{1}\}, \\
\Psi &=&\{\lambda _{\xi }\tau ^{t}\mid \xi \in \mathbb{Z}_{n}^{1},\text{ }%
t\in \mathbb{Z}_{n}\}
\end{eqnarray*}%
and call $\Psi $ the \textit{monic normalizer} of $\overline{\left\langle
\tau \right\rangle }$.

\begin{proposition}
\label{prop18} (i) $\Lambda $ is an abelian group isomorphic to $\mathbb{Z}%
_{n}^{1}$;\newline
(ii) $\Psi =\Lambda \ltimes \overline{\left\langle \tau \right\rangle }\cong 
\mathbb{Z}_{n}^{1}\ltimes \mathbb{Z}_{n}$;\newline
(iii) on letting $\Psi ^{\prime }$ denote the derived subgroup of $\Psi $,
we have $\Psi ^{\prime }=\overline{\left\langle \tau ^{n}\right\rangle }$.
\end{proposition}

\begin{proof}
(i) Let $\xi ,\theta \in \mathbb{Z}_{n}^{1}.$ Then, as $\lambda _{\xi
},\lambda _{\theta }$ and $\lambda _{\xi \theta }$ are inactive, it follows
that 
\begin{equation*}
(\lambda _{\xi }\lambda _{\theta }\lambda _{\xi \theta }^{-1})|_{i}=(\lambda
_{\xi })|_{i}(\lambda _{\theta })|_{i}\left( (\lambda _{\xi \theta
})|_{i}\right) ^{-1}
\end{equation*}%
\begin{equation*}
=\lambda _{\xi }\tau ^{i\frac{\xi -1}{n}}\lambda _{\theta }\tau ^{i\frac{%
\theta -1}{n}}\left( \lambda _{\xi \theta }\tau ^{i\frac{\xi \theta -1}{n}%
}\right) ^{-1}=\lambda _{\xi }\lambda _{\theta }\lambda _{\theta }^{-1}\tau
^{i\frac{\xi -1}{n}}\lambda _{\theta }\tau ^{i\frac{\theta -1}{n}}\tau ^{-i%
\frac{\xi \theta -1}{n}}\lambda _{\xi \theta }^{-1}
\end{equation*}%
\begin{equation*}
=\lambda _{\xi }\lambda _{\theta }\left( \tau ^{i\theta \frac{\xi -1}{n}%
}\tau ^{i\frac{\theta -1}{n}}\tau ^{-i\frac{\xi \theta -1}{n}}\right)
\lambda _{\xi \theta }^{-1}=\lambda _{\xi }\lambda _{\theta }\lambda _{\xi
\theta }^{-1},\forall i\in \{0,\ldots ,n-1\}\text{.}
\end{equation*}

Therefore, $\lambda _{\xi }\lambda _{\theta }=\lambda _{\xi \theta }$. In
addition, $\lambda _{\xi }=e$ if and only if $\xi =1$.

(ii) This factorization is clear.

(iii) Let $\theta =1+\theta ^{\prime }n,\eta \in \mathbb{Z}_{n}$. Then%
\begin{eqnarray*}
\lbrack \tau ^{\eta },\lambda _{\theta }] &=&\tau ^{-\eta }\lambda _{\theta
^{-1}}\tau ^{\eta }\lambda _{\theta }= \\
\tau ^{-\eta }\tau ^{\eta \theta } &=&\tau ^{\eta \left( \theta -1\right)
}=\left( \tau ^{n}\right) ^{\eta \theta ^{\prime }}\text{.}
\end{eqnarray*}
\end{proof}

We prove below the existence of conjugates $\tau ^{\alpha }$ of $\tau $ in $%
N_{\mathcal{A}_{n}}\left( \overline{\left\langle \tau \right\rangle }\right) 
$, which lie outside $\overline{\left\langle \tau \right\rangle }$. This
fact allows us to construct the first important type of metabelian groups $%
\overline{\left\langle \tau \right\rangle }\left\langle \tau ^{\alpha
}\right\rangle $ containing $\tau $.

\begin{proposition}
Given $\xi ,\rho \in \mathbb{Z}_{n}^{1}$ with $\xi \neq 1$. Then for all $n$
odd and for all $n$ even such that $2n\mid (\xi -1)$, an element $\alpha
=(\alpha |_{0},\ldots ,\alpha |_{n-1})$ in $A_{n}$ satisfies $\tau ^{\alpha
}=\lambda _{\xi }\tau ^{\rho }$ if and only if%
\begin{equation*}
\left\{ 
\begin{array}{ll}
\alpha |_{i+1}=\left( \alpha |_{0}\right) \lambda _{\xi ^{i+1}}\tau ^{\frac{1%
}{n}\left[ \rho \frac{\xi ^{i+1}-1}{\xi -1}-(i+1)\right] }\text{ }\left(
0\leq i\leq n-2\right) , &  \\ 
\tau ^{\alpha |_{0}}=\lambda _{\xi ^{n}}\tau ^{\frac{1}{n}\left[ \rho \frac{%
\xi ^{n}-1}{\xi -1}\right] }\text{.} & 
\end{array}%
\right.
\end{equation*}
\end{proposition}

\begin{proof}
From $\displaystyle\tau ^{\alpha }=\lambda _{\xi }\tau ^{1+\kappa n}$, we
obtain using (\ref{eq4}),

\begin{equation*}
\left\{ 
\begin{array}{l}
\lambda _{\xi }\tau ^{i\frac{\xi -1}{n}+\kappa }=\left( \alpha
|_{i}^{-1}\right) \alpha |_{i+1},\;\;\text{if }i\in Y-\{n-1\} \\ 
\lambda _{\xi }\tau ^{(n-1)\frac{\xi -1}{n}+\kappa +1}=\left( \alpha
|_{n-1}^{-1}\right) \tau \left( \alpha |_{0}\right) \text{.}%
\end{array}%
\right.
\end{equation*}%
Therefore,

\begin{equation*}
\begin{array}{l}
\alpha |_{i+1}=\left( \alpha |_{0}\right) \lambda _{\xi }\tau ^{\kappa
}\lambda _{\xi }\tau ^{\frac{\xi -1}{n}+\kappa }\cdots \lambda _{\xi }\tau
^{i\frac{\xi -1}{n}+\kappa },\;\;\text{for }i=0,1,\ldots ,n-2\text{,} \\ 
\alpha |_{0}=\tau ^{-1}\left( \alpha |_{n-1}\right) \lambda _{\xi }\tau
^{(n-1)\frac{\xi -1}{n}+\kappa +1}\text{.}%
\end{array}%
\end{equation*}%
The first equations can be expressed as 
\begin{eqnarray*}
\alpha |_{i+1} &=&\left( \alpha |_{0}\right) \lambda _{\xi ^{i+1}}\tau
^{\kappa \left( \sum_{j=0}^{i}\xi ^{j}\right) +\frac{\xi -1}{n}\xi
^{i}\left( \sum_{j=1}^{i}j(\xi ^{-1})^{j}\right) } \\
&=&\left( \alpha |_{0}\right) \lambda _{\xi ^{i+1}}\tau ^{\frac{1}{n}\left[
(1+\kappa n)\frac{\xi ^{i+1}-1}{\xi -1}-(i+1)\right] }
\end{eqnarray*}%
and the last as

\begin{eqnarray*}
\alpha |_{0} &=&\tau ^{-1}\left( \alpha |_{0}\right) \lambda _{\xi ^{n}}\tau
^{\frac{\xi }{n}\left[ (1+\kappa n)\frac{\xi ^{n-1}-1}{\xi -1}-(n-1)\right]
}\tau ^{(n-1)\frac{\xi -1}{n}+\kappa +1} \\
&=&\tau ^{-1}\left( \alpha |_{0}\right) \lambda _{\xi ^{n}}\tau ^{\frac{1}{n}%
\left[ (1+\kappa n)\frac{\xi ^{n}-1}{\xi -1}\right] }\text{.}
\end{eqnarray*}

Now, we need to show that $\tau ^{\alpha |_{0}}=\lambda _{\xi ^{n}}\tau ^{%
\frac{1}{n}\left[ (1+\kappa {}n)\frac{\xi ^{n}-1}{\xi -1}\right] }$
satisfies the same conditions as those for $\alpha ;$ that is, both $\xi
^{n},\;\rho ^{\prime }=\frac{1}{n}\left[ (1+\kappa {}n)\frac{\xi ^{n}-1}{\xi
-1}\right] \in \mathbb{Z}_{n}^{1}.$

Of course, $\xi ^{n}\in \mathbb{Z}_{n}^{1},$ so let us consider $\rho (\xi
^{n}-1)/n(\xi -1).$ Since $\xi \in \mathbb{Z}_{n}^{1},$ we can write $\xi
=1+\ell n,$ and then 
\begin{equation*}
\frac{\xi ^{n}-1}{\xi -1}\equiv n+\binom{n}{2}\ell n\;(\bmod{\;n^{2}}),
\end{equation*}%
by using the \emph{Binomial Theorem.} Since $\rho \equiv 1\;(\bmod{\;n}),$
it follows that 
\begin{equation*}
\frac{\rho (\xi ^{n}-1)}{n(\xi -1)}\equiv 1+\binom{n}{2}\ell \;(\bmod{\;n}),
\end{equation*}%
and consequently, $\rho (\xi ^{n}-1)/n(\xi -1)\in \mathbb{Z}_{n}^{1}$ if and
only if $n\mid \binom{n}{2}\ell $ (that is, if and only if $(n-1)\ell $ is
even). So $\rho (\xi ^{n}-1)/n(\xi -1)\in \mathbb{Z}_{n}^{1}$ holds for odd $%
n,$ and for even $n$ provided that $2n\mid (\xi -1)$.
\end{proof}

\section{\bf Abelian groups $B$ normalized by $\protect\tau$}

\vskip 0.4 true cm

Let $B$ be an abelian subgroup of $\mathcal{A}_{n}$ normalized by $\tau $.
For a fixed $\beta \in B$, we define the `1st level state closure' of $%
\left\langle \beta ,\tau \right\rangle $ as the group

\begin{equation*}
H=\left\langle \beta |_{i}\text{ }\left( i\in Y\right) ,\tau \right\rangle 
\text{.}
\end{equation*}%
We will be dealing frequently with the following subgroups of $H$ ,%
\begin{eqnarray*}
N &=&\left\langle [\beta |_{i},\tau ^{k_{i}}]\mid k_{i}\in \mathbb{Z},i\in
Y\right\rangle \\
M &=&N\left\langle \tau \right\rangle \text{.}
\end{eqnarray*}%
When $\sigma _{\beta }=\left( \sigma _{\tau }\right) ^{s}$\textit{\ for some
integer }$s$, $m=\frac{n}{\gcd (n,s)}$ and \[ \pi _{i}=\beta |_{i}\beta |_{%
 i+s}\beta |_{ i+2s}\cdots \beta |_{%
i+\left( m-1\right) s}\] \textit{we will also be dealing
with the subgroups }%
\begin{eqnarray*}
K &=&\left\langle N,\;\pi _{i}\mid i\in Y\right\rangle \text{,} \\
O &=&K\left\langle \tau \right\rangle.
\end{eqnarray*}

We show below that when $n$ is a power of a prime number $p^{k}$, the
activity range of $\beta $ narrows down to a Sylow $p$-subgroup of $\Sigma
_{n}$. This is used to restrict the location of an abelian group $B$
normalized by $\tau $, within $\mathcal{A}_{n}.$

\begin{proposition}
\label{casop^k_1} Let $n=p^{k}$, $\sigma =\left( 0,1,...,n-1\right) $ and $P$
be a Sylow $p$-subgroup $P$ of $\Sigma _{n}$ which contains $\sigma $. Then%
\newline
(i) $P$ is isomorphic to $((...\left( ...C_{p})\wr)C_{p}\right) \wr C_{p}$, a
wreath product of the cyclic group $C_{p}$ of order $p$ iterated $k-1$
times; the normalizer of $P$ in $\Sigma _{n}$ is $N_{\Sigma
_{n}}(P)=P\left\langle c\right\rangle $ where $c$ is cyclic of order $p-1$; 
\newline
(ii) $P$ is the unique Sylow $p$-subgroup $P$ of $\Sigma _{n}$ which
contains $\sigma $;\newline
(iii) if $\ W$ is an abelian subgroup of $\Sigma _{n}$ normalized by $\sigma 
$ then $W$ is contained in $P$.\newline
\end{proposition}

\begin{proof}
(i) The structure of $P$ as an iterated wreath product is well-known. The
center of $P$ is $Z=\left\langle z\left( =\sigma ^{p^{k-1}}\right)
\right\rangle $ and $C_{\Sigma _{n}}(z)=P$. Therefore, $N_{\Sigma
_{n}}(P)=N_{\Sigma _{n}}(Z)=P\left\langle c\right\rangle $ where $c$ is
cyclic of order $p-1$.

(ii) If $\sigma \in P^{g}$ for some $g\in \Sigma _{n}$ then $z^{g}\in
C_{\Sigma _{n}}(\sigma )=\left\langle \sigma \right\rangle $ and therefore $%
\left\langle z^{g}\right\rangle =\left\langle z\right\rangle ,$ $P^{g}=P$.
Thus, $P$ is the unique Sylow $p$-subgroup of $\Sigma _{n}$ to contain $%
\sigma $.

(iii) Let $W$ be an abelian subgroup of $\Sigma _{n}$ normalized by $\sigma $%
. Let $V=W\langle \sigma \rangle$ and $V_{0}$ be the stabilizer of $0$ in $V.$ Then,
since $\sigma $ is a regular cycle, it follows that $V=V_{0}\left\langle
\sigma \right\rangle ,$ $V_{0}\cap \left\langle \sigma \right\rangle =\{e\}$%
. Suppose that there exists a prime $q$ different from $p$ which divides the
order of $W$ and let $Q$ be the unique Sylow $q$-subgroup of $W$. Then $Q$
is the unique Sylow $q$-subgroup of $V$ and $Q\leq V_{0}.$ Therefore, $%
Q=\{e\}$ and $W$ is a $p$-group. As $\sigma \in V$, we conclude $W\leq P.$
\end{proof}

\begin{lemma}
\label{lemma5} (a) Let $\gamma \in \mathcal{A}_{n}$. Conditions (i), (ii)
below are equivalent:\newline
(i) $[\gamma ,\gamma ^{\tau ^{k}}]=e$ for all $k\in \mathbb{Z}$;\newline
(ii) $[\tau ^{k},\gamma ,\gamma ]=e$ for all $k\in \mathbb{Z}$.\newline
Condition (i) implies \newline
(iii) $\left\langle [\gamma ,\tau ^{k}]\mid k\in \mathbb{Z}\right\rangle $
is a commutative group.\newline
Condition (iii) implies \newline
$\left\langle [\gamma |_{u},\tau ^{k}]\mid k\in \mathbb{Z}\right\rangle $ is
a commutative group for all indices $u$.\newline
(b) Let $n=p^{k}$. Then any abelian subgroup $B$ normalized by $\tau $ is
contained in the layer closure $L=L\left( N_{\Sigma _{n}}(P)\right) $.
\end{lemma}

\begin{proof}
(a) First,

\begin{equation*}
\begin{array}{ll}
\lbrack \gamma ,\gamma ^{\tau ^{k}}] & =\gamma ^{-1}\left( \tau ^{-k}\gamma
^{-1}\tau ^{k}\right) \gamma \left( \tau ^{-k}\gamma \tau ^{k}\right) \\ 
& =\gamma ^{-1}\left( \tau ^{-k}\gamma ^{-1}\tau ^{k}\gamma \right) \gamma
\left( \gamma ^{-1}\tau ^{-k}\gamma \tau ^{k}\right) \\ 
& =[\tau ^{k},\gamma ]^{\gamma }[\gamma ,\tau ^{k}]%
\end{array}%
\end{equation*}%
and so,%
\begin{equation*}
\lbrack \gamma ,\gamma ^{\tau ^{k}}]=e\Leftrightarrow \lbrack \gamma ,\tau
^{k}]^{\gamma }=[\gamma ,\tau ^{k}]\text{.}
\end{equation*}

Furthermore, since 
\begin{equation}
\lbrack \gamma ,\tau ^{k_{1}}]^{\tau ^{k_{2}}}=[\gamma ,\tau
^{k_{2}}]^{-1}[\gamma ,\tau ^{k_{1}+k_{2}}]  \label{basico}
\end{equation}%
for all integers $k_{1},k_{2}$, condition (ii) implies 
\begin{eqnarray*}
\lbrack \gamma ,\tau ^{k_{1}}]^{[\gamma ,\tau ^{k_{2}}]} &=&[\gamma ,\tau
^{k_{1}}]^{\gamma ^{-1}\tau ^{-k_{2}}\gamma \tau ^{k_{2}}}=[\gamma ,\tau
^{k_{1}}]^{\tau ^{-k_{2}}\gamma \tau ^{k_{2}}} \\
&=&\left( [\gamma ,\tau ^{-k_{2}}]^{-1}[\gamma ,\tau ^{k_{1}-k_{2}}]\right)
^{\gamma \tau ^{k_{2}}}=\left( [\gamma ,\tau ^{-k_{2}}]^{-1}[\gamma ,\tau
^{k_{1}-k_{2}}]\right) ^{\tau ^{k_{2}}}
\end{eqnarray*}%
\begin{equation*}
=[\gamma ,\tau ^{k_{1}}]\text{.}
\end{equation*}

Finally, we note that by (\ref{eq6}), 
\begin{equation*}
\begin{array}{ll}
([\gamma ,\tau ^{nk}])|_{\left( i\right) \sigma _{\gamma }} & =(\gamma
^{-1})|_{_{\left( i\right) \sigma _{\gamma }}}(\tau ^{-nk})|_{i}\left(
\gamma |_{i}\right) (\tau ^{nk})|_{_{\left( i\right) \sigma _{\gamma }}} \\ 
& =\left( \gamma |_{i}^{-1}\right) \tau ^{-k}\left( \gamma |_{i}\right) \tau
^{k} \\ 
& =[\gamma |_{i},\tau ^{k}]\text{.}%
\end{array}%
\end{equation*}

Since $[\gamma ,\tau ^{kn}]$ is inactive for all $k\in \mathbb{Z},$ we
obtain $\{[\gamma |_{i},\tau ^{k}]\mid k\in \mathbb{Z}\}$ is a commutative
set for all $i$. The rest of the assertion follows by induction on the tree
level.

(b) Let $\beta \in B$. Since the normal closure of $\left\langle \sigma
_{\beta }\right\rangle $ under the action of $\left\langle \sigma _{\tau
}\right\rangle $ is an abelian subgroup, it follows that $\sigma _{\beta
}\in P$. Furthermore, as $\left\langle [\beta |_{u},\tau ^{k}]\mid k\in 
\mathbb{Z}\right\rangle $ is an abelian group normalized by $\tau $, it
follows that $[\sigma _{\beta |_{u}},\sigma ]\in P$ and therefore $\sigma
^{\sigma _{\beta |_{u}}}\in P$. Thus, we conclude $\sigma _{\beta |_{u}}\in
N_{\Sigma _{n}}(P)$ and $\beta \in L$.
\end{proof}

\begin{proposition}
\label{inativos} Let $l\geq 1$ and suppose $\alpha ,\gamma \in \mathrm{Stab}%
(l)$ satisfy $[\alpha ,\gamma ^{\tau ^{x}}]=e$ for all $x\in \mathbb{Z}.$
Then 
\[ \displaystyle 
\lbrack \alpha |_{u},\left( {\gamma |_{v}}\right) ^{\tau ^{x}}] =e, %
\forall u,v\in \mathcal{M}  
\text{ having }  |u| =|v|\leq l \text{ and  } \forall x\in \mathbb{Z}.  \] %

\end{proposition}

\begin{proof}
We start with the case $l=1$. Write $x=r+kn$ where $r=\overline{x}$.

By (\ref{eq4}),%
\begin{eqnarray*}
\left( \gamma ^{\tau ^{x}}\right) |_{\left( i\right) \tau ^{x}} &=&\left(
(\tau ^{x})|_{i}^{-1}\right) \left( \gamma |_{i}\right) (\tau ^{x})|_{i}, \\
\left( \gamma ^{\tau ^{x}}\right) |_{i} &=&\tau ^{-k-\delta (i-r,r)}\left(
\gamma |_{\overline{i-r}}\right) \tau ^{k+\delta (i-r,r)}\text{.}
\end{eqnarray*}

As $[\alpha ,\gamma ^{\tau ^{x}}]=e$ and $\alpha ,\gamma ^{\tau ^{x}}\in 
\mathrm{Stab}(1),$we have, for all $i,j,r\in Y$ and all $k,x\in \mathbb{Z}$, 
\begin{eqnarray*}
\lbrack \alpha |_{i},(\gamma ^{\tau ^{x}})|_{i}] &=&e,\text{ }[\alpha
|_{i},\left( \gamma |_{\overline{i-r}}\right) ^{\tau ^{k+\delta
(i-r,r)}}]=e,\; \\
\lbrack \alpha |_{i},\left( \gamma |_{j}\right) ^{\tau ^{x}}] &=&e\text{.}
\end{eqnarray*}%
The general case $l\geq 1$ follows by induction.
\end{proof}

We apply the above proposition to $\beta \in B$.

\begin{corollary}
\label{inativo} Let $\sigma _{\beta }=e.$ Then \ for all $i,j\in Y$ and for
all $x\in \mathbb{Z}$%
\begin{equation*}
\lbrack \left( \beta |_{i}\right) ,\left( \beta |_{j}\right) ^{\tau ^{x}}]=e%
\text{.}
\end{equation*}
\end{corollary}

We derive further relations in $H=\left\langle \beta |_{i}\text{ }\left(
i\in Y\right) ,\tau \right\rangle $.

\begin{proposition}
\label{abp1} Let $\beta \in B$. Then the following relations hold in $H$ for
all $v\in \mathbb{Z}$ and for all $i\in Y$:

\begin{itemize}
\item[(I)] 
\begin{eqnarray*}
&&\left( \tau ^{v}|_{\left( i\right) {\sigma _{\tau }^{-v}}}\right)
^{-1}\left( \beta |_{\left( i\right) {\sigma _{\tau }^{-v}}}\right) \left(
\tau ^{v}|_{\left( i\right) {\sigma _{\tau }^{-v}\sigma _{\beta }}}\right)
\left( \beta |_{\left( i\right) {\sigma _{\tau }^{-v}\sigma _{\beta }\sigma
_{\tau }^{v}}}\right) \\
&=&\left( \beta |_{i}\right) \left( \tau ^{v}|_{\left( i\right) {\sigma
_{\beta }\sigma _{\tau }^{-v}}}\right) ^{-1}\left( \beta |_{\left( i\right) {%
\sigma _{\beta }\sigma _{\tau }^{-v}}}\right) \left( \tau ^{v}|_{\left(
i\right) {\sigma _{\beta }\sigma _{\tau }^{-v}\sigma _{\beta }}}\right) ,
\end{eqnarray*}%
\begin{equation*}
\lbrack \sigma _{\beta },\sigma _{\beta }^{\sigma _{\tau }^{v}}]=e;
\end{equation*}

\item[(II)] 
\begin{equation*}
\lbrack \beta |_{i},\tau ^{v}]^{\beta |_{\left( i\right) {\sigma _{\beta }}%
}}=[\beta |_{\left( i\right) {\sigma _{\beta }}},\tau ^{v}];
\end{equation*}

\item[(III)] 
\begin{equation*}
\left( \beta |_{\left( i\right) \sigma _{\beta }}\right) \left( \beta
|_{\left( i\right) {\sigma _{\beta }^{2}}}\right) \cdots \left( \beta
|_{\left( i\right) {\sigma _{\beta }^{s_{i}}}}\right) \text{ commutes with }%
[\beta |_{i},\tau ^{v}]
\end{equation*}%
where $s_{i}$ is the size of the orbit of $i$ under the action of $%
\left\langle \sigma _{\beta }\right\rangle $.
\end{itemize}
\end{proposition}

\begin{proof}
(I) Clearly $[\beta ,\beta ^{\tau ^{v}}]=e$ implies $[\sigma _{\beta
},\sigma _{\beta }^{\sigma _{\tau }^{v}}]=e$. It also implies 
\begin{equation*}
\left( \beta |_{\left( i\right) {\sigma _{\beta ^{\tau ^{v}}}}}\right)
^{-1}\left( \beta ^{\tau ^{v}}|_{i}\right) ^{-1}\left( \beta |_{i}\right)
\left( \beta ^{\tau ^{v}}|_{\left( i\right) {\sigma _{\beta }}}\right) =e,
\end{equation*}%
\begin{equation*}
\left( \beta ^{\tau ^{v}}|_{i}\right) \left( \beta |_{\left( i\right) {%
\sigma _{\beta ^{\tau ^{v}}}}}\right) =\left( \beta |_{i}\right) \left(
\beta ^{\tau ^{v}}|_{\left( i\right) {\sigma _{\beta }}}\right) ,
\end{equation*}

\begin{eqnarray*}
&&\left( \tau ^{v}|_{\left( i\right) {\sigma _{\tau ^{v}}^{-1}}}\right)
^{-1}\left( \beta |_{\left( i\right) {\sigma _{\tau ^{v}}^{-1}}}\right)
\left( \tau ^{v}|_{\left( i\right) \left( {\sigma _{\tau ^{v}}^{-1}}\right)
\left( {\sigma _{\beta }}\right) }\right) \left( \beta |_{\left( i\right) {%
\sigma _{\beta ^{\tau ^{v}}}}}\right) \\
&=&\left( \beta |_{i}\right) \left( \tau ^{v}|_{\left( i\right) \left( {%
\sigma _{\beta }}\right) \left( {\sigma _{\tau ^{v}}^{-1}}\right) }\right)
^{-1}\left( \beta |_{\left( i\right) \left( {\sigma _{\beta }}\right) \left( 
{\sigma _{\tau ^{v}}^{-1}}\right) }\right) \left( (\tau ^{v})|_{\left(
i\right) {\sigma _{\beta }}\left( {\sigma _{\tau ^{v}}^{-1}}\right) {\sigma
_{\beta }}}\right) \text{.}
\end{eqnarray*}

(II) On changing $v$ to $nv$ in (I), we obtain:%
\begin{equation*}
\tau ^{-v}\left( \beta |_{i}\right) \tau ^{v}\left( \beta |_{\left( i\right) 
{\sigma _{\beta }}}\right) =\left( \beta |_{i}\right) \tau ^{-v}\left( \beta
|_{\left( i\right) {\sigma _{\beta }}}\right) \tau ^{v},
\end{equation*}%
\begin{eqnarray*}
&&\left( \beta |_{\left( i\right) {\sigma _{\beta }}}\right) ^{-1}\left(
\left( \beta |_{i}^{-1}\right) \tau ^{-v}\left( \beta |_{i}\right) \tau
^{v}\right) \left( \beta |_{\left( i\right) {\sigma _{\beta }}}\right) \\
&=&(\left( \beta |_{\left( i\right) {\sigma _{\beta }}}\right) ^{-1}\left(
\beta |_{i}^{-1}\right) )\left( \beta |_{i}\right) \tau ^{-v}\left( \beta
|_{\left( i\right) {\sigma _{\beta }}}\right) \tau ^{v}\text{.}
\end{eqnarray*}

(III) From (II), we derive%
\begin{equation*}
\lbrack \beta |_{i},\tau ^{v}]^{\left( \left( \beta |_{\left( i\right) {%
\sigma _{\beta }}}\right) \left( \beta |_{\left( i\right) {\sigma _{\beta
}^{2}}}\right) \cdots \left( \beta |_{\left( i\right) {\sigma _{\beta
}^{s_{i}}}}\right) \right) }=[\beta |_{\left( i\right) {\sigma _{\beta }}%
},\tau ^{v}]^{\left( \left( \beta |_{\left( i\right) {\sigma _{\beta }^{2}}%
}\right) \cdots \left( \beta |_{\left( i\right) {\sigma _{\beta }^{s_{i}}}%
}\right) \right) }=...=[\beta |_{i},\tau ^{v}].
\end{equation*}
\end{proof}

\section{\bf The case $\protect\beta \in B$ with $\protect\sigma _{\protect\beta %
}\in \left\langle \protect\sigma _{\protect\tau }\right\rangle $}

\vskip 0.4 true cm

This section is devoted to the proof of the second part of Theorem B. For
this purpose, we introduce the following $3$-variable combination of Delta-2%
\textit{\ }functions 
\begin{equation*}
\Delta _{s}(i,t)=\delta (i,t-i)-\delta (i-s,t-i)
\end{equation*}%
which we call the \textit{Delta-3 }function\textit{.}

\begin{lemma}
\label{lema6} Let $\beta \in \mathcal{A}_{n}$ such that $[\beta ,\beta
^{\tau ^{x}}]=e$ for any $x\in \mathbb{Z}$ and let $\sigma _{\beta }=\sigma
_{\tau }^{s}$ for some $s\in Y.$ Then,%
\begin{eqnarray*}
&&\tau ^{\Delta _{s}(i,t)}\left( \beta |_{i-s}\right) [\beta |_{i-s},\tau
^{z}]\left( \beta |_{t}\right) \\
&=&\left( \beta |_{t-s}\right) \left( \beta |_{i}\right) [\beta |_{i},\tau
^{z}]\tau ^{\Delta _{s}(i+s,t+s)},
\end{eqnarray*}
\end{lemma}

for all $i,t\in \{0,1,\ldots ,n-1\},z\in \mathbb{Z}.$

\begin{proof}
Since $\sigma _{\beta }=\sigma _{\tau }^{s}$, we have $\sigma _{\beta ^{\tau
^{x}}}=\sigma _{\beta }=\sigma _{\tau }^{s}.$

From (\ref{eq4}) and  (\ref{eq6}),  we obtain

\begin{equation}
\begin{array}{ll}
& \tau ^{-\frac{x-\overline{x}}{n}-\delta (j-x,x)}\left( \beta
|_{j-x}\right) \tau ^{\frac{x-\overline{x}}{n}+\delta (j-x+s,x)}\left( \beta
|_{j+s}\right) \\ 
= & \left( \beta |_{j}\right) \tau ^{-\frac{x-\overline{x}}{n}-\delta
(j+s-x,x)}\left( \beta |_{j+s-x}\right) \tau ^{\frac{x-\overline{x}}{n}%
+\delta (j+2s-x,x)}%
\end{array}
\label{ajuste1}
\end{equation}

Setting $\displaystyle k=\frac{x-\overline{x}}{n}$ and $r=\overline{x}$ and
using (\ref{ajuste1}), we have

\begin{equation}
\begin{array}{ll}
& \tau ^{-k-\delta (j-r,r)}\left( \beta |_{j-r}\right) \tau ^{k+\delta
(j+s-r,r)}\left( \beta |_{j+s}\right) \\ 
& =\left( \beta |_{j}\right) \tau ^{-k-\delta (j+s-r,r)}\left( \beta
|_{j+s-r}\right) \tau ^{k+\delta (j+2s-r,r)},%
\end{array}
\label{ajuste2}
\end{equation}%
for all $r,j\in Y$ and all $k\in \mathbb{Z}$.

Also on setting $t=\overline{j+s},i=\overline{j+s-r}$ and $z=k+\delta
(j+s-r,r)\left( =k+\delta (i,t-i)\right) $ and using (\ref{ajuste2}), we
obtain 
\begin{equation*}
\begin{array}{ll}
& \tau ^{-z+\delta (i,t-i)-\delta (i-s,t-i)}\beta |_{i-s}\tau ^{z}\beta |_{t}
\\ 
= & \beta |_{t-s}\tau ^{-z}\beta |_{i}\tau ^{z-\delta (i,t-i)+\delta
(i+s,t-i)},%
\end{array}%
\end{equation*}%
for all $t,i\in \{0,1,\ldots ,n-1\}$ and all $z\in \mathbb{Z}$.

It follows that 
\begin{eqnarray*}
&&\tau ^{\delta (i,t-i)-\delta (i-s,t-i)}\left( \beta |_{i-s}\right) [\beta
|_{i-s},\tau ^{z}]\left( \beta |_{t}\right) \\
&=&\left( \beta |_{t-s}\right) \left( \beta |_{i}\right) [\beta |_{i},\tau
^{z}]\tau ^{-\delta (i,t-i)+\delta (i+s,t-i)}
\end{eqnarray*}%
for all $t,i\in \{0,1,\ldots ,n-1\}$ and all $z\in \mathbb{Z}$.
\end{proof}

We develop below some properties of the $\Delta _{s}$ function to be used in
the sequel.

\begin{proposition}
\label{PropDelta} The Delta-3 function satisfies

\begin{itemize}
\item[(i)] $\displaystyle\Delta _{s}(i,t)=\delta (i,-s)-\delta
(t,-s)=\left\{ 
\begin{array}{rl}
0, & \mathrm{if}\;\overline{t},\overline{i}\geq \overline{s}\;\;\mathrm{or}%
\;\;\overline{t},\overline{i}<\overline{s} \\ 
1, & \text{if }\overline{t}<\overline{s}\leq \overline{i} \\ 
-1, & \text{if }\overline{i}<\overline{s}\leq \overline{t}%
\end{array}%
\right. ,$

\item[(ii)] $\Delta _{s}(i,t)=-\Delta _{s}(t,i),$

\item[(iii)] $\Delta _{s}(i+s,t+s)=-\Delta _{-s}(i,t),$

\item[(iv)] $\Delta _{s}(i,t)=\Delta _{s}(i,z)+\Delta _{s}(z,t),$

\item[(v)] $\displaystyle\sum_{k=0}^{\frac{n}{(s,n)}-1}\Delta
_{s}(i+ks,t+ks)=0,$

\item[(vi)] $\displaystyle\sum_{k=0}^{n-1}\Delta _{s}(k,t)=\left\{ 
\begin{array}{ll}
n-\overline{s}, & \text{if }\overline{t}<\overline{s} \\ 
-\overline{s} & \text{if }\overline{t}\geq \overline{s}%
\end{array}%
\right. $
\end{itemize}

for all $i,t,z\in \mathbb{Z}.$
\end{proposition}

\begin{proof}
$\:\:\:$

\begin{itemize}
\item[(i)] Using the definition $\delta (i,j)=\frac{\overline{i}+\overline{j}%
-\overline{i+j}}{n}$ we have 
\begin{equation*}
\begin{array}{ll}
\Delta _{s}(i,t) & \displaystyle=\frac{\overline{i}+\overline{t-i}-\overline{%
t}}{n}-\frac{\overline{i-s}+\overline{t-i}-\overline{t-s}}{n} \\ 
&  \\ 
& \displaystyle=\frac{\overline{i}+\overline{-s}-\overline{i-s}}{n}-\frac{%
\overline{t}+\overline{-s}-\overline{t-s}}{n} \\ 
&  \\ 
& =\delta (i,-s)-\delta (t,-s) \\ 
& =\left\{ 
\begin{array}{rl}
0, & \mathrm{if}\;\overline{t},\overline{i}\geq \overline{s}\;\;\mathrm{or}%
\;\;\overline{t},\overline{i}<\overline{s} \\ 
1, & \mathrm{if}\;\overline{t}<\overline{s}\leq \overline{i} \\ 
-1, & \mathrm{if}\;\overline{i}<\overline{s}\leq \overline{t}%
\end{array}%
\right. \text{.}%
\end{array}%
\end{equation*}

\item[(ii)] Follows from (i).

\item[(iii)] Calculate%
\begin{equation*}
\begin{array}{ll}
\Delta _{s}(i+s,t+s) & =\delta (i+s,t-i)-\delta (i,t-i) \\ 
& =-\left( \delta (i,t-i)-\delta (i+s,t-i)\right) \\ 
& =-\Delta _{-s}(i,t)\text{.}%
\end{array}%
\end{equation*}

\item[(iv)] This part follows from (i).

\item[(v)] From the definition of the \textit{Delta-2} function 
\begin{equation*}
\sum_{k=0}^{\frac{n}{(n,s)}-1}\delta (i+ks,t-i)=\sum_{k=0}^{\frac{n}{(n,s)}%
-1}\delta (i+(k-1)s,t-i)\text{.}
\end{equation*}

\item[(vi)] Finally, we have%
\begin{equation*}
\begin{array}{ll}
\displaystyle\sum_{k=0}^{n-1}\Delta _{s}(k,t) & =\sum_{k=0}^{\overline{s}%
-1}\Delta _{s}(k,t)+\sum_{k=\overline{s}}^{n-1}\Delta _{s}(k,t) \\ 
& \overset{(i)}{=}\left\{ 
\begin{array}{ll}
n-\overline{s}, & \text{if }\overline{t}<\overline{s} \\ 
-\overline{s}, & \text{if }\overline{t}\geq \overline{s}%
\end{array}%
\right. \text{.}%
\end{array}%
\end{equation*}
\end{itemize}
\end{proof}

With the use of the \textit{Delta-3} function we obtain

\begin{proposition}
\label{ptc} The following relations are verified in $H=\left\langle \beta
|_{i}\text{ }\left( i\in Y\right) ,\tau \right\rangle $, for all $x,z\in 
\mathbb{Z}$ and all $i,t\in Y:$

\begin{itemize}
\item[(I)] $\displaystyle\tau ^{\Delta _{s}(i,t)}\left( \beta |_{i-s}\right)
\left( \beta |_{t}\right) =\left( \beta |_{t-s}\right) \left( \beta
|_{i}\right) \tau ^{\Delta _{s}(i+s,t+s)};$

\item[(II)] $\displaystyle\lbrack \left( \beta |_{i-s}\right) ,\tau
^{z}]^{\left( \beta |_{t}\right) \tau ^{-\Delta _{s}(i+s,t+s)}}=[\beta
|_{i},\tau ^{z}];$

\item[(III)] $[[\beta |_{i},\tau ^{z}],[\beta |_{t},\tau ^{x}]]=e$.
\end{itemize}
\end{proposition}

\begin{proof}
Returning to Lemma \ref{lema6}, we have%
\begin{eqnarray*}
&&\tau ^{\Delta _{s}(i,t)}\left( \beta |_{i-s}\right) [\beta |_{i-s},\tau
^{z}]\left( \beta |_{t}\right) \\
&=&\left( \beta |_{t-s}\right) \left( \beta |_{i}\right) [\beta |_{i},\tau
^{z}]\tau ^{\Delta _{s}(i+s,t+s)}\text{.}
\end{eqnarray*}

Consequently, 
\begin{equation}
\tau ^{\Delta _{s}(i,t)}\left( \beta |_{i-s}\right) \left( \beta
|_{t}\right) =\left( \beta |_{t-s}\right) \left( \beta |_{i}\right) \tau
^{\Delta _{s}(i+s,t+s)}  \label{outra0}
\end{equation}%
and 
\begin{equation}
\lbrack \beta |_{i-s},\tau ^{z}]^{\left( \beta |_{t}\right) \tau ^{-\Delta
_{s}(i+s,t+s)}}=[\beta |_{i},\tau ^{z}],  \label{outra}
\end{equation}%
for all $t,i\in Y$ and all $z\in \mathbb{Z}$.

From (\ref{outra}) and (\ref{basico}), $N=\left\langle [\beta |_{i},\tau
^{k_{i}}]\mid k_{i}\in \mathbb{Z},i\in Y\right\rangle $ is a normal subgroup
of $H$. Moreover, by applying alternately the above equations, we obtain 
\begin{equation*}
\lbrack \beta |_{i},\tau ^{z}]^{[\beta |_{t},\tau ^{k}]}=[\beta |_{i},\tau
^{z}]^{\left( \beta |_{t}^{-1}\right) \tau ^{-k}\left( \beta |_{t}\right)
\tau ^{k}}
\end{equation*}%
\begin{equation*}
=[\beta |_{i},\tau ^{z}]^{\left( \tau ^{-\Delta _{s}(i+s,t+s)}\tau ^{\Delta
_{s}(i+s,t+s)}\left( \beta |_{t}^{-1}\right) \tau ^{-k}\left( \beta
|_{t}\right) \tau ^{k}\right) }
\end{equation*}%
\begin{equation*}
\overset{(\ref{basico})}{=}\left( [\beta |_{i},\tau ^{-\Delta
_{s}(i+s,t+s)}]^{-1}{\large .}[\beta |_{i},\tau ^{z-\Delta
_{s}(i+s,t+s)}]\right) ^{\displaystyle\left( \tau ^{\Delta
_{s}(i+s,t+s)}\left( \beta |_{t}^{-1}\right) \tau ^{-k}\left( \beta
|_{t}\right) \tau ^{k}\right) }
\end{equation*}%
\begin{equation*}
\overset{(\ref{outra})}{=}\left( [\beta |_{i-s},\tau ^{-\Delta
_{s}(i+s,t+s)}]^{-1}.[\beta |_{i-s},\tau ^{z-\Delta _{s}(i+s,t+s)}]\right) ^{%
\displaystyle\tau ^{-k}\left( \beta |_{t}\right) \tau ^{k}}
\end{equation*}%
\begin{equation*}
\overset{(\ref{basico})}{=}\left( 
\begin{array}{c}
\left( \lbrack \beta |_{i-s},\tau ^{-k}]^{-1}.[\beta |_{i-s},\tau
^{-k-\Delta _{s}(i+s,t+s)}]\right) ^{-1} \\ 
\left( \lbrack \beta |_{i-s},\tau ^{-k}]^{-1}.[\beta |_{i-s},\tau
^{-k+z-\Delta _{s}(i+s,t+s)}]\right)%
\end{array}%
\right) ^{\displaystyle\beta |_{t}\tau ^{k}}
\end{equation*}%
\begin{equation*}
=\left( [\beta |_{i-s},\tau ^{-k-\Delta _{s}(i+s,t+s)}]^{-1}.[\beta
|_{i-s},\tau ^{-k+z-\Delta _{s}(i+s,t+s)}]\right) ^{\displaystyle\left(
\beta |_{t}\right) \tau ^{k}}
\end{equation*}%
\begin{eqnarray*}
&&\overset{(\ref{outra})}{=}\left( [\beta |_{i},\tau ^{-k-\Delta
_{s}(i+s,t+s)}]^{-1}.[\beta |_{i},\tau ^{-k+z-\Delta _{s}(i+s,t+s)}]\right)
^{\displaystyle\tau ^{k+\Delta _{s}(i+s,t+s)}} \\
&&\overset{(\ref{basico})}{=}[\beta |_{i},\tau ^{z}]\text{.} \\
&&
\end{eqnarray*}
\end{proof}

\begin{corollary}
\label{me(C1)} Let $\beta \in A_{n}$ such that $[\beta ,\beta ^{\tau
^{x}}]=e $ for every $x\in \mathbb{Z}$ with $\sigma _{\beta }=\sigma _{\tau
}^{s}$ for some $s\in \{0,1,\ldots ,n-1\}.$ Then 
\begin{equation*}
M=\left\langle [\beta |_{i},\tau ^{k_{i}}],\tau \mid k_{i}\in \mathbb{Z}%
,0\leq i\leq n-1\right\rangle
\end{equation*}%
is a normal metabelian subgroup of $H.$
\end{corollary}

\begin{proof}
By Proposition \ref{ptc}, $N=\left\langle [\beta |_{i},\tau ^{k_{i}}]\mid
k_{i}\in \mathbb{Z},0\leq i\leq n-1\right\rangle $ is abelian and normal in $%
H$. Since $N\tau \in Z(H/N)$, it follows that $M=N\left\langle \tau
\right\rangle $ is a normal subgroup of $H$ and is clearly metabelian.
\end{proof}

We are ready to prove part \textit{(I) of }Theorem B.

\begin{theorem}
\label{principal2} Let $\beta \in \mathcal{A}_{n}$ be such that $[\beta
,\beta ^{\tau ^{x}}]=e,\forall x\in \mathbb{Z}$ and $\sigma _{\beta }=\sigma
_{\tau }^{s}$ for some $s\in Y$ and $H=\left\langle \beta |_{0},\ldots
,\beta |_{n-1},\tau \right\rangle $. Recall $\pi _{j}=\beta |_{j}\beta |_{ j+s}\beta |_{ j+2s}\cdots \beta |_{%
 j+\left(m-1\right) s}$. Then,$\newline
$(i) the group $K = \left\langle [\beta |_{i},\tau ^{x}],\pi _{j} \mid
i,j\in Y,x\in \mathbb{Z}_{n}\right\rangle $ is an abelian normal subgroup of $H$ and the group $O = K\langle \tau\rangle$ is a  metabelian normal subgroup of 
$H$;\newline
(ii) the quotient group $H/O$ is a homomorphic image of a subgroup of $%
C_{m}\wr C_{n}$.\newline

 In particular, $H$ is metabelian-by-finite.\newline

\end{theorem}

\begin{proof}
(i) Recall 
\begin{eqnarray*}
N &=&\left\langle [\beta |_{i},\tau ^{k_{i}}]\mid k_{i}\in \mathbb{Z},i\in
Y\right\rangle , \\
K &=&N\left\langle \pi _{j}\mid j\in Y\right\rangle
\end{eqnarray*}%
where $m=\frac{n}{\gcd (n,s)}$. Then, by Proposition \ref{ptc}, $N$ is an
abelian normal subgroup of $H$.

By (\ref{outra}), we have 
\begin{equation*}
\begin{array}{cl}
& [\beta |_{i},\tau ^{z}]^{\pi _{j}} \\ 
= & [\beta |_{i+s},\tau ^{z}]^{\tau ^{\Delta _{s}(i+2s,j+s)}\beta |_{
j+s}\cdots \beta |_{ j+\left( m-1\right) s}}
\\ 
= & [\beta |_{i+2s},\tau ^{z}]^{\tau ^{\Delta _{s}(i+2s,j+s)+\Delta
_{s}(i+3s,j+2s)}\beta |_{ j+2s}\cdots \beta |_{
j+\left( m-1\right) s}} \\ 
= & [\beta |_{i},\tau ^{z}]^{\tau ^{\left( \sum_{k=0}^{m-1}\Delta
_{s}(i+(k+1)s,j+ks)\right) }} \\ 
\overset{\mathrm{Prop.}\ref{PropDelta}(v)}{=} & [\beta |_{i},\tau ^{z}]\text{%
.}%
\end{array}%
\end{equation*}

Thus, 
\begin{equation}
\lbrack \lbrack \beta |_{i},\tau ^{z}],(\beta ^{m})|_{j}]=e,\forall i,j\in
Y,\forall z\in \mathbb{Z}  \label{oo1}
\end{equation}

Since $\sigma _{\beta ^{m}}=e,$ we have by Corollary \ref{inativo} 
\begin{equation}
\lbrack (\beta ^{m})|_{i},(\beta ^{m})|_{j}]=e,\forall i,j\in Y\text{.}
\label{com}
\end{equation}

Moreover, 
\begin{equation}
(\beta ^{m})|_{i}^{\displaystyle\tau }=\left( (\beta ^{m})|_{i}\right)
[(\beta ^{m})|_{i},\tau ]\text{.}  \label{oos}
\end{equation}

Since $[\beta ,\beta ^{\tau ^{x}}]=e,\forall x\in \mathbb{Z},$ it follows
that $[\beta ^{m},\beta ^{\tau ^{x}}]=e,\forall x\in \mathbb{Z}.$

Therefore, by (\ref{eq6}), 
\begin{equation*}
e=(\beta ^{m})|_{\left( i\right) {\sigma _{\left( \beta ^{\tau ^{x}}\right) }%
}}^{-1}(\beta ^{\tau ^{x}})|_{i}^{-1}(\beta ^{m})|_{i}(\beta ^{\tau
^{x}})|_{\left( i\right) {\sigma _{\beta ^{m}}}},\forall x\in \mathbb{Z}%
,\forall i\in Y.
\end{equation*}

Now, as $\sigma _{\beta }=\sigma _{\tau }^{s}$ and $\sigma _{\beta ^{m}}=e$,
we reach

\begin{equation}
(\beta ^{m})|_{i+s}=(\beta ^{m})|_{i}^{(\beta ^{\tau ^{x}})|_{i}},\forall
x\in \mathbb{Z},\forall i\in Y\text{.}  \label{facilita}
\end{equation}

By (\ref{eq4}), the following 
\begin{equation*}
(\beta ^{\tau ^{x}})|_{i}=(\tau ^{x})|_{\left( i\right) {\sigma _{\tau
^{x}}^{-1}}}^{-1}\left( \beta |_{\left( i\right) {\sigma _{\tau ^{x}}^{-1}}%
}\right) (\tau ^{x})|_{\left( i\right) \left( {\sigma _{\tau ^{x}}^{-1}}%
\right) {\sigma _{\beta }}}=(\tau ^{x})|_{i-x}^{-1}\left( \beta
|_{i-x}\right) (\tau ^{x})|_{i-x+s}
\end{equation*}%
holds for all $i\in Y$ and all $x\in \mathbb{Z}$.

From which we derive 
\begin{equation}
(\beta ^{\tau ^{x}})|_{i}=\tau ^{-\frac{x-\overline{x}}{n}-\delta
(i-x,x)}\left( \beta |_{\overline{i-x}}\right) \tau ^{\frac{x-\overline{x}}{n%
}+\delta (i-x+s,x)}  \label{facilita1}
\end{equation}%
for all $i\in Y$ and all $x\in \mathbb{Z}.$

Therefore, by (\ref{facilita}) and (\ref{facilita1}),

\begin{equation*}
(\beta ^{m})|_{\overline{i+s}}=(\beta ^{m})|_{i}^{\left( \tau ^{-\frac{x-%
\overline{x}}{n}-\delta (i-x,x)}\left( \beta |_{\overline{i-x}}\right) \tau
^{\frac{x-\overline{x}}{n}+\delta (i-x+s,x)}\right) },
\end{equation*}%
for all $i\in Y$ and all $x\in \mathbb{Z}..$

On writing $x=kn+\overline{x}=kn+r,r\in \mathbb{Z}$ in the above equation,
we obtain 
\begin{equation*}
(\beta ^{m})|_{i+s}=(\beta ^{m})|_{i}^{\displaystyle\tau ^{-k-\delta
(i-r,r)}\left( \beta |_{i-r}\right) \tau ^{k+\delta (i-r+s,r)}}
\end{equation*}%
\begin{equation*}
\Rightarrow (\beta ^{m})|_{i+s}^{\displaystyle\tau ^{-k-\delta
(i-r+s,r)}}=(\beta ^{m})|_{i}^{\displaystyle\left( \beta |_{i-r}\right) \tau
^{-k-\delta (i-r,r)}[\tau ^{-k-\delta (i-r,r)},\beta |_{i-r}]}
\end{equation*}%
\begin{equation*}
\Rightarrow (\beta ^{m})|_{i+s}^{\displaystyle\tau ^{-k-\delta
(i-r+s,r)}[\beta |_{i-r},\tau ^{-k-\delta (i-r,r)}]\tau ^{k+\delta
(i-r,r)}}=(\beta ^{m})|_{i}^{\left( \beta |_{i-r}\right) }
\end{equation*}%
for all $i,r\in Y$ and all $k\in \mathbb{Z}.$

By (\ref{oo1}), (\ref{oos}) and using the fact that $N$ is abelian and
normal in $H$, we find 
\begin{equation*}
(\beta ^{m})|_{i+s}^{\displaystyle\tau ^{\delta (i-r,r)-\delta
(i-r+s,r)}}=(\beta ^{m})|_{i}^{\left( \beta |_{i-r}\right) }
\end{equation*}%
\begin{equation*}
\Rightarrow (\beta ^{m})|_{i+s}^{\displaystyle\tau ^{\Delta_{-s}
(i-r,i-i)}}=(\beta ^{m})|_{i}^{\left( \beta |_{i-r}\right) }
\end{equation*}%
for all $i,r\in Y.$

On setting $j=\overline{i-r},$ we get 
\begin{equation}
(\beta ^{m})|_{i+s}^{\displaystyle\tau ^{\Delta_{-s} (j,i)}}=(\beta
^{m})|_{i}^{\left( \beta |_{j}\right) }  \label{oo2}
\end{equation}%
for all $i,j\in Y.$

Further, by using equations (\ref{oo1}),(\ref{com}) (\ref{oos}), (\ref{oo2})
and 
\begin{equation}
(\beta ^{m})|_{i}=\pi _{i},
\end{equation}%
we conclude that also $K$ is an abelian normal subgroup of $H$.

Now, $O=K\left\langle \tau \right\rangle $ is metabelian. Moreover it is
normal in $H$, because%
\begin{equation*}
\tau ^{\beta |_{i}}=\tau \tau ^{-1}\tau ^{\beta |_{i}}=\tau \lbrack \tau
,\beta |_{i}]\in O
\end{equation*}

for all $i\in Y$.

(ii) Consider the following Fibonacci type group 
\begin{equation*}
X=\left\langle 
\begin{array}{c}
b_{0},\ldots ,b_{n-1}\mid b_{i}b_{\overline{j+s}}=b_{j}b_{\overline{i+s}},
\\ 
b_{i}b_{\overline{i+s}}\cdots {}b_{\overline{i+(n-1)s}}=e,\forall i,j\in Y%
\end{array}%
\right\rangle
\end{equation*}%
where the bar notation indicates \ 
%TCIMACRO{\U{b4}}%
%BeginExpansion
\'{}%
%EndExpansion
modulo $n$ 
%TCIMACRO{\U{b4}}%
%BeginExpansion
\'{}%
%EndExpansion
.

The Equation (\ref{outra0}) shows that $\displaystyle\frac{H}{O}$ is a
homomorphic image of $X$. We will prove that $X$ is isomorphic to a subgroup
of the wreath product $C_{m}\wr C_{n}$.

As a matter of fact the group $C_{m}\wr C_{n}$ has the presentation

\begin{equation*}
\left\langle u,a\mid
u^{m}=e,a^{n}=e,u^{a^{i}}u^{a^{j}}=u^{a^{j}}u^{a^{i}}\right\rangle \text{.}
\end{equation*}

On defining $b=a^{s}u^{-1}$, we have

\begin{equation*}
\begin{array}{ll}
& u^{m}=e \\ 
\Rightarrow & (a^{-s}b)^{m}=e \\ 
\Rightarrow & (\underbrace{a^{-s}b\cdots {}a^{-s}b}_{m\;\mathrm{terms}%
})^{a^{-s+i}}=e \\ 
\Rightarrow & b^{a^{i}}b^{a^{i+s}}\cdots b^{a^{i+(m-1)s}}=e\text{.}%
\end{array}%
\end{equation*}%
Also, the commutation relation 
\begin{equation*}
u^{a^{i}}u^{a^{j}}=u^{a^{j}}u^{a^{i}}
\end{equation*}%
implies%
\begin{equation*}
\begin{array}{ll}
& 
(b^{-1}a^{s})^{a^{i}}(b^{-1}a^{s})^{a^{j}}=(b^{-1}a^{s})^{a^{j}}(b^{-1}a^{s})^{a^{i}}
\\ 
\Rightarrow & 
(a^{-s}b)^{a^{j}}(a^{-s}b)^{a^{i}}=(a^{-s}b)^{a^{i}}(a^{-s}b)^{a^{j}} \\ 
\Rightarrow & b^{a^{j}}a^{-s}b^{a^{i}}=b^{a^{i}}a^{-s}b^{a^{j}} \\ 
\Rightarrow & b^{a^{j}}b^{a^{i+s}}=b^{a^{i}}b^{a^{j+s}}\text{.}%
\end{array}%
\end{equation*}

By using Tietze transformations we conclude that $C_{m}\wr C_{n}$ has the
presentation

\begin{equation*}
\left\langle a,b\mid
a^{n}=e,b^{a^{j}}b^{a^{i+s}}=b^{a^{i}}b^{a^{j+s}},b^{a^{i}}b^{a^{i+s}}\cdots
b^{a^{i+(m-1)s}}=e,\forall i,j\in Y\right\rangle \text{.}
\end{equation*}

Then, on introducing $b_{i}=b^{a^{i}},i=0,\ldots ,n-1$, the above
presentation is expressed as

\begin{equation*}
\left\langle a,b_{0},\ldots ,b_{n-1}\mid
a^{n}=e,b_{i}=b_{0}^{a^{i}},\;b_{j}b_{\overline{i+s}}=b_{i}b_{\overline{j+s}%
},\;b_{i}b_{\overline{i+s}}\cdots {}b_{\overline{i+(m-1)s}}=e,\right.
\end{equation*}%
\begin{equation*}
\left. \forall i,j\in Y\right\rangle \text{.}
\end{equation*}
\end{proof}

The next results leads to a proof of Theorem C.

\begin{lemma}
\label{congexp} Let $\sigma =\left( 0,1,...,n-1\right) \in \Sigma _{n}$ and
let $L$ be the layer closure of $\left\langle \sigma \right\rangle $ in $%
\mathcal{A}_{n}$. Suppose $\beta =(\beta |_{0},\beta |_{1},\ldots ,\beta
|_{n-1})\sigma _{\beta }\in L$ satisfies $[\beta ,\beta ^{\tau ^{x}}]=e$ for
all $x\in \mathbb{Z}$. Write $\sigma _{\beta }=\sigma ^{s}$ and $\sigma
_{\beta |_{i}}=\sigma ^{m_{i}}$ for all $i\in Y$. Then for all $i,t\in Y$,
the following congruence holds

\begin{equation}
\Delta _{s}(i,t)+m_{\overline{i-s}}+m_{t}\equiv m_{\overline{t-s}%
}+m_{i}+\Delta _{s}(i+s,t+s)\bmod{n}.  \label{quasetudo}
\end{equation}
\end{lemma}

\begin{proof}
Since $\sigma _{\beta |_{i}}=\sigma ^{m_{i}}$, we conclude by (\ref{outra0}%
), 
\begin{equation*}
\sigma ^{\Delta _{s}(i,t)+m_{\overline{i-s}}+m_{t}}=\sigma ^{m_{\overline{t-s%
}}+m_{i}+\Delta _{s}(i+s,t+s)}
\end{equation*}

and therefore, $\displaystyle\Delta _{s}(i,t)+m_{\overline{i-s}}+m_{t}\equiv
m_{\overline{t-s}}+m_{i}+\Delta _{s}(i+s,t+s)\bmod{n}.$
\end{proof}

\begin{lemma}
\label{lemmasigma} Maintain the notation of the previous lemma and let $s=1$%
. Then,%
\begin{equation*}
\sigma _{\left( \beta ^{n}\right) |_{0}}=\sigma _{\left( \beta |_{0}\right)
\left( \beta |_{1}\right) \cdots \left( \beta |_{n-1}\right) }=\sigma \text{.%
}
\end{equation*}
\end{lemma}

\begin{proof}
The case $n=2$ is covered by Proposition 6 of \cite{Sid03}.

Now let $n$ be an odd prime. From 
\begin{equation*}
\Delta _{1}(i,t)+m_{\overline{i-1}}+m_{t}\equiv m_{\overline{t-1}%
}+m_{i}+\Delta _{1}(i+1,t+1)\bmod n
\end{equation*}%
we conclude 
\begin{equation*}
\begin{array}{ll}
& \displaystyle\sum_{i=0}^{n-2}\sum_{t=i+1}^{n-1}\left( \Delta _{1}(i,t)+m_{%
\overline{i-1}}+m_{t}\right)  \\ 
\equiv  & \displaystyle\sum_{i=0}^{n-2}\sum_{t=i+1}^{n-1}\left( m_{\overline{%
t-1}}+m_{i}+\Delta _{1}(i+1,t+1)\right) \bmod n\text{.}%
\end{array}%
\end{equation*}

Now,

\begin{equation*}
\sum_{i=0}^{n-2}\sum_{t=i+1}^{n-1}\Delta _{1}(i,t)\overset{\mathrm{Prop.\ref%
{PropDelta}(i)}}{=}\sum_{t=1}^{n-1}\Delta _{1}(0,t)\overset{\mathrm{Prop.\ref%
{PropDelta}(ii)}}{=}\sum_{t=0}^{n-1}\Delta _{1}(0,t)
\end{equation*}%
\begin{equation*}
\overset{\mathrm{Prop.\ref{PropDelta}(ii)}}{=}\sum_{t=0}^{n-1}-\Delta
_{1}(t,0)\overset{\mathrm{Prop.\ref{PropDelta}(vi)}}{=}-(n-1),
\end{equation*}%
\begin{equation*}
\sum_{i=0}^{n-2}\sum_{t=i+1}^{n-1}\Delta _{1}(i+1,t+1)\overset{\mathrm{Prop.%
\ref{PropDelta}(i)}}{=}\sum_{i=0}^{n-2}\Delta _{1}(i+1,0)\overset{\mathrm{%
Prop.\ref{PropDelta}(ii)}}{=}\sum_{i=0}^{n-1}\Delta _{1}(i,0)
\end{equation*}%
\begin{equation*}
\overset{\mathrm{Prop.\ref{PropDelta}(vi)}}{=}(n-1),
\end{equation*}%
\begin{equation*}
\sum_{i=0}^{n-2}\sum_{t=i+1}^{n-1}\left( m_{\overline{i-1}}+m_{t}\right)
=2(n-1)m_{n-1}+(n-2)\sum_{k=0}^{n-2}m_{k}
\end{equation*}%
and 
\begin{equation*}
\sum_{i=0}^{n-2}\sum_{t=i+1}^{n-1}\left( m_{\overline{t-1}}+m_{i}\right)
=n\sum_{k=0}^{n-2}m_{k}\text{.}
\end{equation*}%
Since $n$ is odd, we have 
\begin{equation*}
\sum_{k=0}^{n-1}m_{k}\equiv 1\bmod{n}
\end{equation*}
and therefore, $\sigma _{\left( \beta |_{0}\right) \cdots \left( \beta
|_{n-1}\right) }=\sigma ^{ m_{0}+\cdots + m_{n-1}}=\sigma $.
\end{proof}

Now we prove Theorem C.

\begin{theorem}
\label{=sigma} \textit{Let }$n$\textit{\ be an odd number,} $\sigma
=(0,\ldots ,n-1)\in \Sigma _{n}$\textit{\ and let }$L$\textit{\ be the layer
closure of }$\left\langle \sigma \right\rangle $\textit{\ in }$\mathcal{A}_{n}$\textit{%
.\ Let }$s$ be an integer relatively prime to $n$ and $\beta =(\beta
|_{0},\beta |_{1},\ldots ,\beta |_{n-1})\sigma ^{s}\in L$\textit{\ be such }%
that\textit{\ }$[\beta ,\beta ^{\tau ^{x}}]=e$\textit{\ for all }$x\in 
\mathbb{Z}.$\textit{\ Then }$\beta $\textit{\ is a conjugate of }$\tau $%
\textit{\ in } $\mathcal{A}_{n}$\textit{.}
\end{theorem}

\begin{proof}
We start with the case $s=1$. The element 
\begin{equation*}
\alpha (1)=(e,\beta |_{0}^{-1},(\left( \beta |_{0}\right) \left( \beta
|_{1}\right) )^{-1},\ldots ,(\left( \beta |_{0}\right) \cdots \left( \beta
|_{n-2}\right) )^{-1})\in \mathrm{Stab}(1)
\end{equation*}%
conjugates $\beta $ to 
\begin{equation*}
\beta ^{\alpha (1)}=(e,\ldots ,e,\left( \beta |_{0}\right) \cdots \left(
\beta |_{n-1}\right) )\sigma .
\end{equation*}

By Lemma \ref{lemmasigma} we find $\sigma _{\left( \beta |_{0}\right) \cdots
\left( \beta |_{n-1}\right) }=\sigma .$ Moreover by Proposition \ref%
{inativos}, 
\begin{equation*}
\lbrack (\beta ^{n})|_{0},(\left( \beta ^{n})|_{0}\right) ^{\tau
^{x}}]=[\left( \beta |_{0}\right) \cdots \left( \beta |_{n-1}\right)
,(\left( \beta |_{0}\right) \cdots \left( \beta |_{n-1}\right) )^{\tau
^{x}}]=e,
\end{equation*}
for all integers $x$. Therefore $\left( \beta |_{0}\right) \cdots \left(
\beta |_{n-1}\right) $ satisfies the hypothesis of the theorem. The process
can be repeated until we obtain a sequence $\displaystyle\left( \alpha
(k)\right) _{k\in \mathbb{N}}$ such that $\beta ^{\alpha (1)\alpha (2)\cdots
\alpha (k)\cdots }=\tau ,$ where $\alpha (k)\in \mathrm{Stab}(k)$
satisfies $\alpha (k)|_{u}=\alpha (k)|_{v}$ for all $u,v\in \mathcal{M}$
with $|u|=|v|=k-1.$

Now, suppose more generally $s$ is such that $\gcd \left( s,n\right) =1$ and let $
k$ be the minimum positive integer for which $sk\equiv 1$ mod$\left(
n\right) $. Then $\beta ^{k}$ satisfies the hypothesis of the first part and
so, there exists $\alpha \in L$ such that $(\beta ^{k})^{\alpha }=\tau $.
Since $k$ is invertible in $\mathbb{Z}_{n}$, there exists  $
\gamma  \in \mathcal{A}_{n} $ such that $\tau ^{\gamma }=\tau ^{k^{-1}}$. Thus, $
\beta ^{\alpha \gamma ^{-1}}=\tau $.
\end{proof}

\section{\bf Solvable groups for $n=p$, a prime number}

\vskip 0.4 true cm

We will prove in this section the case $n=p$ of Theorem A.

Let $B$ be an abelian subgroup of $Aut(T_{p})$ normalized by $\tau $ and let 
$\beta \in B$. By Proposition \ref{casop^k_1}, $\sigma _{\beta }\in
\left\langle \sigma _{\tau }\right\rangle $ and therefore we have in effect
two cases, namely, $\sigma _{\beta }=e,\sigma _{\tau }$.

\begin{proposition}
Suppose $\sigma _{\beta }\in \left<\sigma _{\tau }\right>$. Then, $\sigma _{\left( \beta
|_{i}\right) }\in \left\langle \sigma _{\tau }\right\rangle $ for all $i\in
Y $.
\end{proposition}

\begin{proof}
By Theorem \ref{principal2}, $K$ is an abelian normal subgroup of $H$ and $\frac{H}{O}
$ is homomorphic to a subgroup of $C_{p}\wr C_{p}$ for $O = K\langle \tau \rangle.$ 

By Proposition \ref{casop^k_1}, $K$ is a subgroup of $\left\langle \sigma
_{\tau }\right\rangle $ modulo $Stab(1).$ So the same is true for $O = K\langle \tau \rangle.$ 

Therefore, $H$ is a $p$-group modulo $Stab(1).$ Since  $H$ is a $p$-group modulo $Stab(1)$ and since $\tau \in H, $ it follows that $H$ coincides with $\langle \sigma_{\tau}\rangle$ modulo $Stab(1), $  by Proposition \ref{casop^k_1}. Hence, necessarily
we have $\sigma _{\left( \beta |_{i}\right) }\in \left\langle \sigma _{\tau
}\right\rangle $.
\end{proof}

\begin{theorem}
\label{csn=p_1} Let $p$ be a prime number and $\beta \in \mathrm{Aut}(T_{p})$
such that $\sigma _{\beta }=\sigma _{\tau }^{s}$ for some integer $s$
relatively prime to $p$. Suppose $[\beta ,\beta ^{\tau ^{x}}]=e$ for all $%
x\in \mathbb{Z}.$ Then $\beta $ is conjugate to $\tau $ in $\mathrm{Aut}%
(T_{p})$.
\end{theorem}

\begin{proof} As the second author showed the case p=2 in \cite{Sid03}, we will show the case p odd.

Suppose $s=1$. Recall that 
\begin{equation*}
\alpha (1)=(e,\beta |_{0}^{-1},(\beta |_{0}\beta |_{1})^{-1},\ldots ,(\beta
|_{0}\cdots \beta |_{p-2})^{-1})\in \mathrm{Stab}(1)
\end{equation*}%
conjugates $\beta $ to its normal form 
\begin{equation*}
\beta ^{\alpha (1)}=(e,\ldots ,e,\beta |_{0}\cdots \beta |_{p-1})\sigma .
\end{equation*}

By Lemma \ref{lemmasigma} we have $\sigma _{\beta |_{0}\beta |_{1}\cdots
\beta |_{p-1}}=\sigma _{\tau }.$ Moreover by Proposition \ref{inativos}, 
\begin{equation*}
\lbrack \beta ^{p}|_{0},(\beta ^{p}|_{0})^{\tau ^{x}}]=[\beta |_{0}\beta
|_{1}\cdots \beta |_{p-1},(\beta |_{0}\beta |_{1}\cdots \beta |_{p-1})^{\tau
^{x}}]=e,
\end{equation*}%
for all integers $x.$ Therefore $\beta |_{0}\beta |_{1}\cdots \beta |_{p-1}$
satisfies the condition of the theorem. This process can be repeated to
produce a sequence $\displaystyle\left( \alpha (k)\right) _{k\in \mathbb{N}}$
such that $\beta ^{\alpha (1)\alpha (2)\cdots \alpha (k)\cdots }=\tau ,$
where $\alpha (k)\in \mathrm{Stab}(k)$ satisfies $\alpha (k)|_{u}=\alpha
(k)|_{v}$ for all $u,v\in \mathcal{M}$ where $|u|=|v|=k-1$.

In the general case, $s$ is such that $\gcd (p,s)=1$. Let $k$ be the minimum
positive integer which is the inverse of $s$ modulo $p$. Then, $\sigma
|_{\beta ^{k}}=\sigma _{\tau }$ and $\beta ^{k}$ satisfies the hypotheses.
Thus there exists $\alpha \in \mathcal{A}_{p}$ such that $\left( \beta
^{k}\right) ^{\alpha }=\tau $. Let $k^{-1}$ be the inverse of $k$ in $%
U\left( \mathbb{Z}_{n}\right) $; then $\beta ^{\alpha }=\tau ^{k^{-1}}$.
There exists $\gamma \in N_{\mathcal{A}_{p}}\left(\overline{<\tau >}\right)$
which conjugates $\tau $ to $\tau ^{k^{-1}}$ and so, $\left( \beta ^{\alpha
}\right) ^{\gamma ^{-1}}=\tau $.
\end{proof}

\begin{lemma}
\label{csn=p_3} Let $p$ be a prime number and $\beta \in \mathrm{Aut}(T_{p})$
such that $[\beta ,\beta ^{\tau ^{x}}]=e$ for all $x\in \mathbb{Z}.$ Then,
there exists a tree level $m$ and a conjugate $\mu $ of $\tau $ such that $%
\beta \in \times _{p^{m}}\overline{\left\langle \mu \right\rangle }$ and
there exists an index $u$ of length $m$ such that $\beta |_{u}=\mu $.
\end{lemma}

\begin{proof}
Let $m$ be the minimum tree level such that $\sigma _{\beta |_{u}}\neq e$
for some $|u|=m$. Therefore, $\beta \in Stab(m)$ and $\sigma _{\beta |_{u}}=\sigma _{\tau }^{s}$ for
some integer $s$ such that $\gcd \left( p,s\right) =1.$ By Proposition \ref{inativos}, $[\beta|_{u}, \beta|_{v}^{\tau^{k}}] = e $ for all indices $v$ such that $|v| = m$ and 
for all $k \in \mathbb{Z}.$ So, by Theorem \ref{csn=p_1},  $\mu =\beta
|_{u}$ is conjugate to $\tau $ in $\mathrm{Aut}(T_{p})$ and  $\beta |_{v}\in \overline{%
\left\langle \mu \right\rangle }$ for all $v\ $such that $|v|=m,$ by Lemma \ref{centralizadortau}.
\end{proof}

\begin{theorem}
Let $p$ be a prime number, $\sigma =(0,1,\ldots ,p-1)\in \Sigma _{p}$, $%
F=N_{\Sigma _{p}}\left( \left\langle \sigma \right\rangle \right) $, $\Gamma
_{0}=N_{\mathcal{A}}(\overline{<\tau >})$ . Let $G$ be a  
solvable subgroup of $Aut\left( T_{p}\right) $ which contains the $p$-adic
adding machine $\tau $. Then, there exists an integer $t\geq 1$ such that $G$
is conjugate to a subgroup of 
\begin{equation*}
\times _{p}\left( \cdots \left( \times _{p}\left( \times _{p}\Gamma
_{0}\rtimes F\right) \rtimes \right) \cdots \right) \rtimes F,
\end{equation*}
where $\times _{p}$ appears $t$ times.
\end{theorem}

\begin{proof}
We may suppose $G$ has derived length $d\geq 2$. Let $B$ be the $(d-1)$-th
term of the derived series of $G.$ By Lemma \ref{csn=p_3}, there exists a
level $t$ such that $B$ is a subgroup of $V=\times _{p^{t}}\overline{%
\left\langle \mu \right\rangle }$ where $\mu =\tau ^{\alpha }$ for some $%
\alpha \in Aut\left( T_{p}\right) $.

We will show that $G$ is a subgroup of 
\begin{equation*}
\dot{J}=\times _{p}\left( \cdots \left( \times _{p}\left( \times _{p}\left(
\Gamma _{0}\right) ^{\alpha }\rtimes \Sigma _{p}\right) \rtimes \Sigma
_{p}\right) \cdots \right) \rtimes \Sigma _{p},
\end{equation*}%
where $\times _{p}$ appears $t$ times.

Let $\gamma \in G\backslash \dot{J}.$ Then there exists an index $w$ of
length $t$ such that $\gamma |_{w}\not\in \left( \Gamma _{0}\right) ^{\alpha
}$. Since $B$ is an abelian subgroup normalized by $\tau $ and $\tau $ is transitive on all levels of the tree,   by Lemma \ref%
{csn=p_3}, there exists $\beta \in B$ such that $\beta |_{w}=\mu ^{\eta }$
for some $\eta \in U(\mathbb{Z}_{p})$.

Write $v=w^{\gamma }$. Then, 
\begin{equation*}
\left( \beta ^{\gamma }\right) |_{v}\overset{(\ref{eq9})}{=}\left( \beta
|_{v^{\gamma ^{-1}}}\right) ^{\left( \gamma |_{v^{\gamma ^{-1}}}\right)
}=\left( \beta |_{w}\right) ^{\gamma |_{w}}\not\in \overline{\left\langle
\mu \right\rangle },
\end{equation*}%
and this implies $\beta ^{\gamma }\not\in B\leq \times_{p^{t}} \overline{%
\left\langle \mu \right\rangle }$ and $\gamma \not\in G$. Hence, $G$ is a
subgroup of $\dot{J} $.

Now, since $G$ is a solvable group containing $\tau $, there exist $G_{i}$ $%
\left( 0\leq i\leq t\right) $ solvable subgroups of $\Sigma _{p}$ containing 
$\sigma =(0,1,\ldots ,p-1)$ such that $G$ is a subgroup of 
\begin{equation*}
R_{t}\left( \alpha \right) =\times _{p}\left( \cdots \left( \times
_{p}\left( \times _{p}\left( \Gamma _{0}\right) ^{\alpha }\rtimes
G_{1}\right) \rtimes G_{2}\right) \cdots \right) \rtimes G_{t}\text{.}
\end{equation*}%
Since for all $i$, we have $G_{i}\leq $ $F$ we may substitute every the $%
G_{i}$ by $F$. Finally, $R_{t}\left( \alpha \right) $ is a conjugate of $%
R_{t}\left( 1\right) $ by the diagonal automorphism $\alpha ^{\left(
t\right) }$.
\end{proof}

\section{\bf Two cases for $n$ even}

\vskip 0.4 true cm

We prove in this section part \textit{(II)} of Theorem B.

\subsection{The case $\protect\sigma _{\protect\beta }=\left( \protect\sigma %
_{\protect\tau }\right) ^{\frac{n}{2}}$}

\label{casosigma2}

\begin{theorem}
\label{2n} Let $n$ be an even number, $\beta \in \mathcal{A}_{n}$ such that $%
\sigma _{\beta }=\sigma _{\tau }^{\frac{n}{2}}$ and $\displaystyle\lbrack
\beta ,\beta ^{\tau ^{x}}]=e$ for all $x\in \mathbb{Z}.$ Then $%
H=\left\langle \beta |_{i}\text{ }\left( 0\leq i\leq n-1\right) ,\tau
\right\rangle $ is a metabelian subgroup of $\mathcal{A}_{n}.$
\end{theorem}

\begin{proof}
Denote $\Delta _{\frac{n}{2}}(i,j)$ by $\Delta (i,j)$.

Define the subgroup 
\begin{equation*}
R=\left\langle [\beta |_{t},\tau ^{k}],\text{ }\beta |_{i}\beta |_{i+\frac{n%
}{2}},\text{ }\beta |_{j}^{2}\tau ^{-\Delta (j,j+\frac{n}{2})}\mid k\in 
\mathbb{Z}\text{ and }i,j,t\in Y\right\rangle \text{.}
\end{equation*}

We will prove that $R$ is an abelian normal subgroup of $H$.

\begin{itemize}
\item[(I)] $R$ is normal in $H:$

\begin{itemize}
\item $\left\langle [\beta |_{i},\tau ^{k}]\right\rangle ^{H}\leq R:$ 
\begin{equation*}
\lbrack \beta |_{i+\frac{n}{2}},\tau ^{k}]^{\beta |_{j}}\overset{(\ref{outra}%
)}{=}[\beta |_{i},\tau ^{k}]^{\tau ^{\Delta (j,i)}};
\end{equation*}

\item $\left\langle \beta |_{i}\beta|_{i+\frac{n}{2}}\right\rangle ^{H}\leq
R:$

\begin{eqnarray*}
(\beta |_{i}\beta |_{i+\frac{n}{2}})^{\tau ^{k}} &=&\left( \beta |_{i}\beta
|_{i+\frac{n}{2}}\right) .[\beta |_{i}\beta |_{i+\frac{n}{2}},\tau ^{k}] \\
&=&\left( \beta |_{i}\beta |_{i+\frac{n}{2}}\right) [\beta |_{i},\tau
^{k}]^{\beta |_{i+\frac{n}{2}}}[\beta |_{i+\frac{n}{2}},\tau ^{k}]
\end{eqnarray*}%

\begin{equation}
(\beta |_{i}\beta |_{i+\frac{n}{2}})^{\beta |_{j}}=\left( \beta
|_{j}^{-1}\beta |_{i}\beta |_{i+\frac{n}{2}}\beta |_{j}\right) \tau ^{\Delta
(j+\frac{n}{2},i+\frac{n}{2})}\tau ^{-\Delta (j+\frac{n}{2},i+\frac{n}{2})}
\label{xyz0}
\end{equation}%
\begin{equation*}
\overset{(\ref{outra0})}{=}\left( \beta |_{j}^{-1}\beta |_{i}\right) \tau
^{\Delta (j,i)}\left( \beta |_{j+\frac{n}{2}}\beta |_{i}\right) \tau
^{-\Delta (j+\frac{n}{2},i+\frac{n}{2})}
\end{equation*}%
\begin{equation*}
=\left( \beta |_{j}^{-1}\beta |_{i}\beta |_{j+\frac{n}{2}}\right) \tau
^{\Delta (j,i)}.[\tau ^{\Delta (j,i)},\beta |_{j+\frac{n}{2}}].\beta
|_{i}\tau ^{-\Delta (j+\frac{n}{2},i+\frac{n}{2})}
\end{equation*}%
\begin{eqnarray*}
&&\overset{(\ref{outra0})}{=}\left( \beta |_{j}^{-1}\right) \tau ^{\Delta (j+%
\frac{n}{2},i+\frac{n}{2})}\left( \beta |_{j}\beta |_{i+\frac{n}{2}}\right) .
\\
&&\lbrack \tau ^{\Delta (j,i)},\beta |_{j+\frac{n}{2}}].\beta |_{i}\tau
^{-\Delta (j+\frac{n}{2},i+\frac{n}{2})}
\end{eqnarray*}%
\begin{eqnarray*}
&=&\tau ^{\Delta (j+\frac{n}{2},i+\frac{n}{2})}.[\tau ^{\Delta (j+\frac{n}{2}%
,i+\frac{n}{2})},\beta |_{j}]. \\
&&\beta |_{i+\frac{n}{2}}[\tau ^{\Delta (j,i)},\beta |_{j+\frac{n}{2}}]\beta
|_{i}\tau ^{-\Delta (j+\frac{n}{2},i+\frac{n}{2})}
\end{eqnarray*}%
\begin{eqnarray*}
&&\overset{Prop.\ref{PropDelta}}{=}\tau ^{-\Delta (j,i)}[\tau ^{-\Delta
(j,i)},\beta |_{j}].\beta |_{i+\frac{n}{2}}. \\
&&\lbrack \tau ^{\Delta (j,i)},\beta |_{j+\frac{n}{2}}]\beta |_{i}\tau
^{\Delta (j,i)}
\end{eqnarray*}%
\begin{equation*}
\overset{(\ref{outra})}{=}\tau ^{-\Delta (j,i)}\beta |_{i+\frac{n}{2}}.[\tau
^{-\Delta (j,i)},\beta |_{j+\frac{n}{2}}]^{\tau ^{\Delta (j,i)}}.[\tau
^{\Delta (j,i)},\beta |_{j+\frac{n}{2}}].\beta |_{i}\tau ^{\Delta (j,i)}
\end{equation*}%
\begin{equation*}
\overset{(\ref{basico})}{=}\left( \beta |_{i+\frac{n}{2}}\beta |_{i}\right)
^{\tau ^{\Delta (j,i)}}\text{.}
\end{equation*}

\item $\left\langle \beta |_{j}^{2}\tau ^{-\Delta (j,j+\frac{n}{2}%
)}\right\rangle ^{H}\leq R:$

\begin{equation*}
\begin{array}{ll}
& (\beta |_{j}^{2}\tau ^{-\Delta (j,j+\frac{n}{2})})^{\tau ^{k}}=\beta
|_{j}^{2}\tau ^{-\Delta (j,j+\frac{n}{2})}.[\beta |_{j}^{2}\tau ^{-\Delta
(j,j+\frac{n}{2})},\tau ^{k}] \\ 
= & \beta |_{j}^{2}\tau ^{-\Delta (j,j+\frac{n}{2})}.[\beta |_{j}^{2},\tau
^{k}]^{\tau ^{-\Delta (j,j+\frac{n}{2})}} \\ 
= & \beta |_{j}^{2}\tau ^{-\Delta (j,j+\frac{n}{2})}\left( [\beta |_{j},\tau
^{k}]^{\beta |_{j}}.[\beta |_{j},\tau ^{k}]\right) ^{\tau ^{-\Delta (j,j+%
\frac{n}{2})}} \\ 
\end{array}%
\end{equation*}

By Proposition \ref{PropDelta} and \ref{ptc}, we can show

\begin{equation}
\left( \beta |_{j}^{2}\tau ^{-\Delta (j,j+\frac{n}{2})}\right) ^{\beta
|_{i}}=\left( \beta |_{j+\frac{n}{2}}^{2}\tau ^{-\Delta (j+\frac{n}{2}%
,j)}[\tau ^{-\Delta (j+\frac{n}{2},j)},\beta |_{j+\frac{n}{2}}]\right)
^{\tau ^{\Delta (i,j)}}\text{.}  \label{xyz1}
\end{equation}
\end{itemize}

\item[(II)] The subgroup $R$ is abelian$:$ 
\begin{equation}
\lbrack \beta |_{i},\tau ^{k}]^{\beta |_{j}\tau ^{t}}\overset{Prop.\ref{ptc}}%
{=}[\beta |_{i},\tau ^{k}]^{\tau ^{t}\beta |_{j}};  \label{xyz2}
\end{equation}

\begin{equation}  \label{xyz3}
[\beta|_{i}, \tau^{k}]^{\beta|_{j}\beta|_{j+\frac{n}{2}}} \overset{(\ref%
{outra})}{=} [\beta|_{i+\frac{n}{2}}, \tau^{k}]^{\tau^{\Delta(j, i+\frac{n}{2%
})}\beta|_{j+\frac{n}{2}}} \overset{(\ref{xyz2})}{=} [\beta|_{i+\frac{n}{2}%
}, \tau^{k}]^{\beta|_{j+\frac{n}{2}}\tau^{\Delta(j, i+\frac{n}{2})}}
\end{equation}
\begin{equation*}
\overset{(\ref{outra})}{=} [\beta|_{i}, \tau^{k}]^{\tau^{\Delta(j+\frac{n}{2}%
,i)+\Delta(j, i+\frac{n}{2})} } \overset{\mathrm{Prop.} \ref{PropDelta}}{=}
[\beta|_{i}, \tau^{k}]
\end{equation*}

\begin{equation}  \label{xyz4}
\begin{array}{ll}
& [\beta|_{i},\tau^{k}]^{\beta|_{j}^{2}\tau^{-\Delta(j, j+\frac{n}{2})}} 
\overset{(\ref{outra})}{=} [\beta|_{i+\frac{n}{2}},
\tau^{k}]^{\tau^{\Delta(j, i+\frac{n}{2})}\beta|_{j}\tau^{-\Delta(j,j+\frac{n%
}{2})}} \\ 
\overset{(\ref{xyz2})}{=} & [\beta|_{i+\frac{n}{2}},
\tau^{k}]^{\beta|_{j}\tau^{\Delta(j, i+\frac{n}{2})-\Delta(j,j+\frac{n}{2}%
)}} \overset{(\ref{outra})}{=} [\beta|_{i},
\tau^{k}]^{\tau^{\Delta(j,i)+\Delta(j, i+\frac{n}{2})-\Delta(j, j+\frac{n}{2}%
)} } \\ 
\overset{\mathrm{Prop.} \ref{PropDelta}}{=} & [\beta|_{i}, \tau^{k}]%
\end{array}%
\end{equation}

\begin{equation*}
\begin{array}{ll}
\left( \beta |_{i}\beta |_{i+\frac{n}{2}}\right) ^{\beta |_{j}\beta |_{j+%
\frac{n}{2}}} & \overset{(\ref{xyz0})}{=}\left( \beta |_{i+\frac{n}{2}}\beta
|_{i}\right) ^{\tau ^{\Delta (j,i)}\beta |_{j+\frac{n}{2}}} \\ 
& =\left( \beta |_{i+\frac{n}{2}}\beta |_{i}\right) ^{\left( \beta |_{j+%
\frac{n}{2}}\tau ^{\Delta (j,i)}[\tau ^{\Delta (j,i)},\beta |_{j+\frac{n}{2}%
}]\right) } \\ 
& \overset{(\ref{xyz0})}{=}\left( \beta |_{i}\beta |_{i+\frac{n}{2}}\right)
^{\left( \tau ^{\Delta (j+\frac{n}{2},i+\frac{n}{2})+\Delta (j,i)}.[\tau
^{\Delta (j,i)},\beta |_{j+\frac{n}{2}}]\right) } \\ 
& \overset{\mathrm{Prop.}\ref{PropDelta}}{=}\left( \beta |_{i}\beta |_{i+%
\frac{n}{2}}\right) ^{[\tau ^{\Delta (j,i)},\beta |_{j+\frac{n}{2}}]} \\ 
& \overset{(\ref{xyz3})}{=}\beta |_{i}\beta |_{i+\frac{n}{2}}%
\end{array}%
\end{equation*}

\begin{equation*}
\begin{array}{ll}
(\beta|_{i}\beta|_{i+\frac{n}{2}})^{\beta|_{j}^{2}\tau^{-\Delta(j,j+\frac{n}{%
2})}} & \overset{(\ref{xyz0})}{=} (\beta|_{i+\frac{n}{2}}\beta|_{i})^{\tau^{%
\Delta(j,i)}\beta|_{j}\tau^{-\Delta(j,j+\frac{n}{2})}} \\ 
& = (\beta|_{i+\frac{n}{2}}\beta|_{i})^{\beta|_{j}\tau^{\Delta(j,i)}[\tau^{%
\Delta(j,i)},\beta|_{j}]\tau^{-\Delta(j,j+\frac{n}{2})}} \\ 
& =(\beta|_{i}\beta|_{i+\frac{n}{2}})^{\tau^{\Delta(j,i+\frac{n}{2}%
)+\Delta(j,i)}[\tau^{\Delta(j,i)},\beta|_{j}]\tau^{-\Delta(j,j+\frac{n}{2})}}
\\ 
& \overset{\mathrm{Prop.} \ref{PropDelta}}{=} (\beta|_{i}\beta|_{i+\frac{n}{2%
}})^{[\tau^{\Delta(j,i)},\beta|_{j}]^{\tau^{\Delta(j+\frac{n}{2},j)}}} \\ 
& \overset{(\ref{basico}),(\ref{xyz3})}{=} \beta|_{i}\beta|_{i+\frac{n}{2}}%
\end{array}%
\end{equation*}

Let 
\begin{equation}
\alpha =\beta |_{j}^{2}\tau ^{-\Delta (j,j+\frac{n}{2})}[\tau ^{-\Delta (j,j+%
\frac{n}{2})},\beta |_{j}]\text{.}  \label{reduz}
\end{equation}%
Then, 
\begin{equation*}
\begin{array}{ll}
& \left( \beta |_{j}^{2}\tau ^{-\Delta (j,j+\frac{n}{2})}\right) ^{\beta
|_{i}^{2}\tau ^{-\Delta (i,i+\frac{n}{2})}} \\ 
\overset{(\ref{xyz1})}{=} & \left( \beta |_{j+\frac{n}{2}}^{2}\tau ^{-\Delta
(j+\frac{n}{2},j)}.[\tau ^{-\Delta (j+\frac{n}{2},j)},\beta |_{j+\frac{n}{2}%
}]\right) ^{\tau ^{\Delta (i,j)}\beta |_{i}\tau ^{-\Delta (i,i+\frac{n}{2})}}%
\end{array}%
\end{equation*}%
\begin{equation*}
=\left( \beta |_{j+\frac{n}{2}}^{2}\tau ^{-\Delta (j+\frac{n}{2},j)}.[\tau
^{-\Delta (j+\frac{n}{2},j)},\beta |_{j+\frac{n}{2}}]\right) ^{\left( \beta
|_{i}\tau ^{\Delta (i,j)}.[\tau ^{\Delta (i,j)},\beta |_{i}].\tau ^{-\Delta
(i,i+\frac{n}{2})}\right) }
\end{equation*}%
\begin{equation*}
=\left( \left( \beta |_{j+\frac{n}{2}}^{2}\tau ^{-\Delta (j+\frac{n}{2}%
,j)}\right) ^{\beta |_{i}}.[\tau ^{-\Delta (j+\frac{n}{2},j)},\beta |_{j+%
\frac{n}{2}}]^{\beta |_{i}}\right) ^{\left( \tau ^{\Delta (i,j)}.[\tau
^{\Delta (i,j)},\beta |_{i}].\tau ^{-\Delta (i,i+\frac{n}{2})}\right) }
\end{equation*}%
\begin{equation*}
\overset{(\ref{outra})}{=}\left( \left( \beta |_{j+\frac{n}{2}}^{2}\tau
^{-\Delta (j+\frac{n}{2},j)}\right) ^{\beta |_{i}}.[\tau ^{-\Delta (j+\frac{n%
}{2},j)},\beta |_{j}]^{\tau ^{\Delta (i,j)}}\right) ^{\left( \tau ^{\Delta
(i,j)}.[\tau ^{\Delta (i,j)},\beta |_{i}].\tau ^{-\Delta (i,i+\frac{n}{2}%
)}\right) }
\end{equation*}%
\begin{equation*}
\overset{(\ref{xyz1})}{=}\left( \alpha ^{\tau ^{\Delta (i,j+\frac{n}{2}%
)}}.[\tau ^{-\Delta (j+\frac{n}{2},j)},\beta |_{j}]^{\tau ^{\Delta
(i,j)}}\right) ^{\left( \tau ^{\Delta (i,j)}.[\tau ^{\Delta (i,j)},\beta
|_{i}].\tau ^{-\Delta (i,i+\frac{n}{2})}\right) }
\end{equation*}%
\begin{equation*}
=\left( \alpha .[\tau ^{-\Delta (j+\frac{n}{2},j)},\beta |_{j}]^{\tau
^{\Delta (i,j)-\Delta (i,j+\frac{n}{2})}}\right) ^{\left( \tau ^{\Delta (i,j+%
\frac{n}{2})+\Delta (i,j)}.[\tau ^{\Delta (i,j)},\beta |_{i}].\tau ^{-\Delta
(i,i+\frac{n}{2})}\right) }
\end{equation*}%
\begin{equation*}
\overset{\mathrm{Prop.}\ref{PropDelta}}{=}\left( \alpha .[\tau ^{-\Delta (j+%
\frac{n}{2},j)},\beta |_{j}]^{\tau ^{\Delta (j+\frac{n}{2},j)}}\right)
^{\left( \tau ^{\Delta (i,i+\frac{n}{2})}[\tau ^{\Delta (i,j)},\beta
|_{i}]\tau ^{-\Delta (i,i+\frac{n}{2})}\right) }
\end{equation*}%
\begin{equation*}
\overset{(\ref{reduz})}{=}\left( \beta |_{j}^{2}\tau ^{-\Delta (j,j+\frac{n}{%
2})}[\tau ^{-\Delta (j,j+\frac{n}{2})},\beta |_{j}][\tau ^{\Delta (j+\frac{n%
}{2},j)},\beta |_{j}]^{-1}\right) ^{[\tau ^{\Delta (i,j)},\beta |_{i}]^{\tau
^{-\Delta (i,i+\frac{n}{2})}}}
\end{equation*}%
\begin{equation*}
\overset{\mathrm{Prop.}\ref{PropDelta}}{=}\left( \beta |_{j}^{2}\tau
^{-\Delta (j,j+\frac{n}{2})}\right) ^{[\tau ^{\Delta (i,j)},\beta
|_{i}]^{\tau ^{-\Delta (i,i+\frac{n}{2})}}}
\end{equation*}%
\begin{equation*}
\overset{(\ref{basico}),(\ref{xyz4})}{=}\beta
|_{j}^{2}\tau ^{-\Delta (j,j+\frac{n}{2})}\text{.}
\end{equation*}
\end{itemize}

Moreover, since 
\begin{equation*}
\begin{array}{ll}
R\left( \beta |_{i}\right) R\left( \beta |_{j}\right) & =R\left( \beta
|_{i}\right) \left( \beta |_{j}\right) \overset{\mathrm{Prop.}\ref{ptc}}{=}%
R\tau ^{\Delta (j,i+\frac{n}{2})}\beta |_{j+\frac{n}{2}}\beta |_{i+\frac{n}{2%
}}\tau ^{\Delta (j,i+\frac{n}{2})} \\ 
&  \\ 
& =R\beta |_{j+\frac{n}{2}}\beta |_{i+\frac{n}{2}}\tau ^{2\Delta (j,i+\frac{n%
}{2})}=R\beta |_{j}^{-1}\beta |_{i}^{-1}\tau ^{2\Delta (j,i+\frac{n}{2})} \\ 
&  \\ 
& =R\beta |_{j}^{-1}\beta |_{j}^{2}\tau ^{-\Delta (j,j+\frac{n}{2})}\beta
|_{i}^{-1}\beta |_{i}^{2}\tau ^{-\Delta (i,i+\frac{n}{2})}\tau ^{2\Delta
(j,i+\frac{n}{2})} \\ 
&  \\ 
& =R\beta |_{j}\beta |_{i}\tau ^{-\Delta (j,j+\frac{n}{2})-\Delta (i,i+\frac{%
n}{2})+2\Delta (j,i+\frac{n}{2})} \\ 
&  \\ 
& \overset{\mathrm{Prop.}\ref{PropDelta}}{=}R\beta |_{j}\beta |_{i}=R\beta
|_{j}R\beta |_{i}%
\end{array}%
\end{equation*}%
and 
\begin{equation*}
R\beta |_{i}=R\beta |_{i+\frac{n}{2}}^{-1},\;\;R\beta |_{i}^{2}=R\tau
^{\Delta (i,i+\frac{n}{2})},\forall i\in Y,
\end{equation*}%
we conclude $\displaystyle\frac{H}{R}$ is a homomorphic image of 
\begin{equation*}
\mathbb{Z}\times \underbrace{C_{2}\times \cdots \times C_{2}}_{\frac{n}{2}%
\;\;\mathrm{terms}}\text{.}
\end{equation*}
\end{proof}

\subsection{The case $\protect\sigma _{\protect\beta }$ transposition}

\label{2ciclo}

\begin{theorem}
\label{me(T1)} Let $n$ be an even number and $B$ an abelian subgroup of $%
\mathcal{A}_{n}$ normalized by $\tau $. Suppose $\beta =(\beta |_{0},\beta
|_{1},\ldots ,\beta |_{n-1})\sigma _{\beta }\in B$ where $\sigma _{\beta }$
is a transposition. Then $H=\left\langle \beta |_{i}\text{ }\left( 0\leq
i\leq n-1\right) ,\tau \right\rangle $ is a metabelian group.
\end{theorem}

We prove progressively that 
\begin{eqnarray*}
N &=&\left\langle [\beta |_{i},\tau ^{k}]\mid k\in \mathbb{Z},i\in
Y\right\rangle , \\
U &=&\left\langle N,\text{ \ }\beta |_{j}\mid j\neq 0,\frac{n}{2}%
\right\rangle , \\
V &=&\left\langle U,\text{ \ }\beta |_{\frac{n}{2}}\beta |_{0},\text{ \ }%
\tau \left( \beta |_{0}\right) ^{2}\right\rangle
\end{eqnarray*}%
are normal abelian subgroups of $H$, from which it follows that $\frac{H}{V}$
is cyclic and therefore $H$ metabelian.

\begin{lemma}
If the degree of the tree $n$ is even and $\sigma _{\beta }$ is a transposition, then $\sigma_{\beta} $ is  $\left\langle
\sigma _{\tau }\right\rangle $-conjugate to the transposition $\left( 0,%
\frac{n}{2}\right) $.
\end{lemma}

\begin{proof}
On conjugating by an appropriate power of $\sigma _{\tau }$, we may assume $%
\sigma _{\beta }=\left( 0,j\right) $. The conjugate of $\sigma _{\beta }$ by 
$\sigma _{\tau }^{i}$ is the transposition $\left( i,j+i\right) $. In
particular, $\left( j,2j\right) $ is a conjugate which is supposed to
commute with $\left( 0,j\right) $. Therefore, $\left\{ 0,j\right\} =\left\{
j,2j\right\} $, $2j=0$ modulo$\left( n\right) $, $n=2n^{\prime }$ and $%
j=n^{\prime }$.
\end{proof}

We go back to part (I) of the Proposition \ref{abp1}, 
\begin{eqnarray*}
&&\left( \tau ^{v}|_{\left( i\right) {\sigma _{\tau }^{-v}}}\right)
^{-1}\left( \beta |_{\left( i\right) {\sigma _{\tau }^{-v}}}\right) \left(
\tau ^{v}|_{\left( i\right) {\sigma _{\tau }^{-v}\sigma _{\beta }}}\right)
\left( \beta |_{\left( i\right) {\sigma _{\tau }^{-v}\sigma _{\beta }\sigma
_{\tau }^{v}}}\right) \\
&=&\left( \beta |_{i}\right) \left( \tau ^{v}|_{\left( i\right) {\sigma
_{\beta }\sigma _{\tau }^{-v}}}\right) ^{-1}\left( \beta |_{\left( i\right) {%
\sigma _{\beta }\sigma _{\tau }^{-v}}}\right) \left( \tau ^{v}|_{\left(
i\right) {\sigma _{\beta }\sigma _{\tau }^{-v}\sigma _{\beta }}}\right)
\end{eqnarray*}%
and set in it $j=\left( i\right) {\sigma _{\tau }^{-v}}$, $v=kn+r$, $r=%
\overline{v}$ to obtain%
\begin{eqnarray}
&&(\tau ^{v})|_{j}^{-1}\beta |_{j}(\tau ^{v})|_{\left( j\right) {\sigma
_{\beta }}}\beta |_{\left( j\right) {\sigma _{\beta }\sigma _{\tau }^{v}}} \\
&=&\beta |_{\left( j\right) {\sigma _{\tau }^{v}}}(\tau ^{v})|_{\left(
j\right) {\sigma _{\tau }^{v}\sigma _{\beta }\sigma _{\tau }^{-v}}%
}^{-1}\beta |_{\left( j\right) {\sigma _{\tau }^{v}\sigma _{\beta }\sigma
_{\tau }^{-v}}}(\tau ^{v})|_{\left( j\right) {\sigma _{\tau }^{v}\sigma
_{\beta }\sigma _{\tau }^{-v}\sigma _{\beta }}}\text{.}
\end{eqnarray}

\begin{proposition}
The following cases hold for different pairs $\left( j,r\right) $.

\begin{itemize}
\item For $j=0$ there are $3$ subcases

\begin{itemize}
\item If $r =0, $ then 
\begin{equation}  \label{eqme2}
[\beta |_{0},\tau ^{k}]^{\beta |_{\frac{n}{2}}} = [\beta |_{\frac{n}{ 2}%
},\tau ^{k}],\;\forall k\in \mathbb{Z}\text{;}
\end{equation}

\item If $r=\frac{n}{2},$ then 
\begin{equation}
\beta |_{0}\tau \beta |_{0}=\beta |_{\frac{n}{2}}\tau ^{-1}\beta |_{\frac{n}{%
2}},  \label{eqme3}
\end{equation}%
and 
\begin{equation}
\lbrack \beta |_{0},\tau ^{k}]^{\tau \left( \beta |_{0}\right) }=[\beta |_{%
\frac{n}{2}},\tau ^{k}],\forall k\in \mathbb{Z}.  \label{eqme4}
\end{equation}

\item If $r \not=0$ and $r \not=\frac{n}{2},$ then

\begin{equation}
\tau ^{\delta (\frac{n}{2},r)}\beta |_{0}\beta |_{\left( \frac{n}{2}%
+r\right) }=\beta |_{r}\tau ^{\delta (\frac{n}{2},r)}\beta |_{0},\forall
r\in Y-\{0,\frac{n}{2}\}  \label{eqme5}
\end{equation}%
and 
\begin{equation}
\lbrack \beta |_{0},\tau ^{k}]^{\beta |_{r}}=[\beta |_{0},\tau ^{k}],\forall
k\in \mathbb{Z}.  \label{eqme6}
\end{equation}
\end{itemize}

\item For $j=\frac{n}{2}$ there are $3$ subcases

\begin{itemize}
\item If $r =0, $ then 
\begin{equation}  \label{eqme8}
[ \beta |_{\frac{n}{2}},\tau ^{k}]^{\beta |_{0}} =[\beta |_{0},\tau
^{k}],\;\forall k\in \mathbb{Z}\text{;}
\end{equation}

\item If $r=\frac{n}{2}$, then 
\begin{equation}
\tau ^{-1}\beta |_{\frac{n}{2}}^{2}=\beta |_{0}^{2}\tau ,  \label{eqme9}
\end{equation}%
and 
\begin{equation}
\lbrack \beta |_{\frac{n}{2}},\tau ^{k}]^{\left( \beta |_{\frac{n}{2}%
}\right) \tau ^{-1}}=[\beta |_{0},\tau ^{k}],\forall k\in \mathbb{Z}\text{;}
\label{eqme10}
\end{equation}

\item If $r \not=0$ and $r \not=\frac{n}{2},$ then 
\begin{equation}  \label{eqme11}
\tau^{-\delta(\frac{n}{2}, r)} \beta|_{\frac{n}{2}} \beta|_{r} = \beta|_{%
\frac{n}{2} + r} \tau^{-\delta(\frac{n}{2}, r)} \beta|_{\frac{n}{2}},
\forall r \in Y - \{0, \frac{n}{2}\}
\end{equation}
and 
\begin{equation}  \label{eqme12}
[ \beta |_{\frac{n}{2}},\tau ^{k}]^{\beta |_{r}} =[\beta |_{\frac{n}{2}%
},\tau ^{k}],\forall k\in \mathbb{Z}, \forall r \in Y - \{0, \frac{n}{2} \}.
\end{equation}
\end{itemize}

\item For $j\not=0$ and $j \not=\frac{n}{2}$, there are $5$ subcases:

\begin{itemize}
\item If $j \not=n - r $ and $j \not= \frac{n}{2} - r,$ then, by substitution $t = j + r, $  we have
\begin{equation}  \label{eqme14}
\beta|_{j} \beta_{t} = \beta|_{t}\beta|_{j}, \forall j, t \in Y - \{0, \frac{%
n}{2}\}
\end{equation}
and 
\begin{equation}  \label{eqme15}
[\beta|_{j}, \tau^{k}]^{\beta|_{t}} = [\beta|_{j}, \tau^{k}], \forall j, t
\in Y - \{0, \frac{n}{2} \}
\end{equation}

\item If $j = n - r$ and $0 < r < \frac{n}{2}, $ then, by substitution $ t = j - \frac{n}{2}, $ we have 
\begin{equation}  \label{eqme16}
\tau^{-1} \beta|_{t + \frac{n}{2}} \tau \beta|_{0} = \beta|_{0}\beta|_{t},
\forall t\in \{1, 2, \cdots, \frac{n}{2} - 1\}
\end{equation}
and 
\begin{equation}  \label{eqme17}
[\beta|_{t+\frac{n}{2}}, \tau^{k}]^{\tau\beta|_{0}} = [\beta|_{t},
\tau^{k}], \forall j \in \{ 1, 2, \cdots, \frac{n}{2} - 1\}
\end{equation}

\item If $j = n - r$ and $\frac{n}{2} < r \leq n - 1, $ then 
\begin{equation}  \label{eqme18}
\beta|_{j}\beta|_{0} = \beta|_{0}\beta|_{\frac{n}{2} + j}, \forall j \in
\{1, \cdots, \frac{n}{2} - 1\}
\end{equation}
and 
\begin{equation}  \label{eqme19}
[\beta|_{j}, \tau^{k}]^{\beta|_{0}} = [\beta|_{\frac{n}{2} + j}, \tau^{k}],
\forall k \in \mathbb{Z}, \forall j \in \{1, \cdots, \frac{n}{2} - 1\}
\end{equation}

\item If $j=\frac{n}{2}-r$ and $0<r<\frac{n}{2},$ then 
\begin{equation}
\beta |_{j}\beta |_{\frac{n}{2}}=\beta |_{\frac{n}{2}}\tau ^{-1}\beta |_{j+%
\frac{n}{2}}\tau ,\forall j\in \{1,\ldots ,\frac{n}{2}-1\}  \label{eqme20}
\end{equation}%
\begin{equation*}
\end{equation*}%
and 
\begin{equation}
\lbrack \beta |_{j},\tau ^{k}]^{\beta |_{\frac{n}{2}}\tau ^{-1}}=[\beta |_{%
\frac{n}{2}+j},\tau ^{k}],\forall k\in \mathbb{Z},\forall j\in \{1,\ldots ,%
\frac{n}{2}-1\}  \label{eqme21}
\end{equation}

\item If $j = \frac{n}{2} - r $ and $\frac{n}{2} < r \leq n - 1,$ then 
\begin{equation}  \label{eqme22}
\beta|_{\frac{n}{2}} \beta|_{j} = \beta|_{\frac{n}{2} + j} \beta|_{\frac{n}{2%
}}, \forall j \in \{1, \cdots, \frac{n}{2} - 1\}
\end{equation}
and 
\begin{equation}  \label{eqme23}
[\beta|_{j}, \tau^{k}] = [\beta|_{\frac{n}{2} + j}, \tau^{k}]^{\beta|_{\frac{%
n}{2}}}, \forall k \in \mathbb{Z}, \forall j \in \{1, \cdots, \frac{n}{2} -
1\}.
\end{equation}
\end{itemize}
\end{itemize}
\end{proposition}

\begin{proof}
We will prove (\ref{eqme14}) and (\ref{eqme15}). As $j\not\in \{0,\frac{n}{2},n-r,\frac{n}{2%
}-r\}$, we have%
\begin{eqnarray*}
\left( j\right) {\sigma _{\tau }^{v}} &=&\left( j\right) {\sigma _{\beta
}\sigma _{\tau }^{v}}=j+r, \\
\left( j\right) {\sigma _{\beta }} &=&\left( j\right) {\sigma _{\tau
}^{v}\sigma _{\beta }\sigma _{\tau }^{-v}}=\left( j\right) {\sigma _{\tau
}^{v}\sigma _{\beta }\sigma _{\tau }^{-v}\sigma _{\beta }}=j.
\end{eqnarray*}

Therefore, 
\begin{equation*}
\begin{array}{ll}
& \displaystyle\left( (\tau ^{v})|_{j}^{-1}\beta |_{j}(\tau ^{v})|_{j}\beta
|_{j+r}=\beta |_{j+r}(\tau ^{v})|_{j}^{-1}\beta |_{j}(\tau ^{v})_{j},\forall
v\in \mathbb{Z}\right) \\ 
&  \\ 
\Leftrightarrow & \displaystyle\left( \tau ^{-k-\delta (j,r)}\beta |_{j}\tau
^{k+\delta (j,r)}\beta |_{j+r}=\beta |_{j+r}\tau ^{-k-\delta (j,r)}\beta
|_{j}\tau ^{k+\delta (j,r)},\forall k\in \mathbb{Z}\right) \\ 
&  \\ 
\Leftrightarrow & \displaystyle\left( \beta |_{j}[\beta |_{j},\tau
^{k+\delta (j,r)}]\beta |_{j+r}=\beta |_{j+r}\beta |_{j}[\beta |_{j},\tau
^{k+\delta (j,r)}],\forall k\in \mathbb{Z}\right) \text{,}%
\end{array}%
\end{equation*}

\begin{equation}
\beta |_{j}\beta _{t}=\beta |_{t}\beta |_{j},\forall j,t\in Y-\{0,\frac{n}{2}%
\}
\end{equation}%
and 
\begin{equation}
\lbrack \beta |_{j},\tau ^{k}]^{\beta |_{t}}=[\beta |_{j},\tau ^{k}],\forall
j,t\in Y-\{0,\frac{n}{2}\}\text{.}
\end{equation}
\end{proof}

\begin{lemma}
\label{me(L1)} The group $N=\left\langle [\beta |_{i},\text{ \ }\tau
^{k}]\mid k\in \mathbb{Z},i\in Y\right\rangle $ is an abelian normal
subgroup of $H$.
\end{lemma}

\begin{proof}
Define 
\begin{equation*}
N_{i}=\left\langle [\beta |_{i},\tau ^{k}]\mid k\in \mathbb{Z}\right\rangle
\end{equation*}%
for each $i\in Y$. Then, $N=\left\langle N_{i}\mid i\in Y\right\rangle $,
each $N_{i}$ is an abelian subgroup normalized by $\tau $ and 
\begin{equation}
\lbrack \beta |_{i},\tau ^{k}]^{\beta |_{j}^{-1}}=[\beta |_{i},\tau
^{k}],\forall k\in \mathbb{Z},\forall i,j\in Y,j\neq 0,\frac{n}{2}
\label{eqme24}
\end{equation}%
We have $[N_{i},N_{j}]=1,\forall i,j\in Y,j\neq 0,\frac{n}{2}$, because 
\begin{eqnarray*}
\lbrack \beta |_{i},\tau ^{k}]^{[\beta |_{j},\tau ^{t}]} &=&[\beta
|_{i},\tau ^{k}]^{\beta |_{j}^{-1}\tau ^{-t}\beta |_{j}\tau ^{t}}\overset{(%
\ref{eqme24})}{=}[\beta |_{i},\tau ^{k}]^{\tau ^{-t}\beta |_{j}\tau ^{t}} \\
&&\overset{(\ref{basico})}{=}\left( [\beta |_{i},\tau ^{-t}]^{-1}[\beta
|_{i},\tau ^{k-t}]\right) ^{\beta |_{j}\tau ^{t}}
\end{eqnarray*}%
\begin{eqnarray*}
&&\overset{(\ref{eqme24})}{=}\left( [\beta |_{i},\tau ^{-t}]^{-1}[\beta
|_{i},\tau ^{k-t}]\right) ^{\tau ^{t}} \\
\overset{(\ref{basico})}{=}[\beta |_{i},\tau ^{k}]^{\tau ^{-t}\tau ^{t}}
&=&[\beta |_{i},\tau ^{k}],\forall k,t\in \mathbb{Z}\text{,}
\end{eqnarray*}%
$\forall i,j\in Y,j\neq 0,\frac{n}{2}$.

Furthermore, $[N_{0},N_{\frac{n}{2}}]=1$, because 
\begin{eqnarray*}
\lbrack \beta |_{\frac{n}{2}},\tau ^{k}]^{[\beta |_{0},\tau ^{t}]} &=&[\beta
|_{\frac{n}{2}},\tau ^{k}]^{\beta |_{0}^{-1}\tau ^{-t}\beta |_{0}\tau ^{t}}%
\overset{(\ref{eqme4})}{=}[\beta |_{0},\tau ^{k}]^{\tau \tau ^{-t}\beta
|_{0}\tau ^{t}} \\
&&\overset{(\ref{basico})}{=}\left( [\beta |_{0},\tau ^{-t}]^{-1}[\beta
|_{0},\tau ^{k-t}]\right) ^{\tau \beta |_{0}\tau ^{t}}
\end{eqnarray*}%
\begin{equation*}
\overset{(\ref{eqme4})}{=}\left( [\beta |_{\frac{n}{2}},\tau
^{-t}]^{-1}[\beta |_{\frac{n}{2}},\tau ^{k-t}]\right) ^{\tau ^{t}}
\end{equation*}%
\begin{equation*}
\overset{(\ref{basico})}{=}[\beta |_{\frac{n}{2}},\tau ^{k}]^{\tau ^{-t}\tau
^{t}}=[\beta |_{\frac{n}{2}},\tau ^{k}],\forall k,t\in \mathbb{Z}\text{.}
\end{equation*}%
Therefore $N$ is abelian.

Now, equation (\ref{eqme24}) implies%
\begin{equation}
N_{i}=N_{i}^{\beta |_{j}}=N_{i}^{\beta |_{j}^{-1}},\forall i,j\in Y,j\neq 0,%
\frac{n}{2}\text{;}  \label{eqme26}
\end{equation}%
equation (\ref{eqme4}) imply
\begin{equation}
\left\{ N_{\frac{n}{2}}=N_{0}^{\beta |_{\frac{n}{2}}},\text{ }N_{0}=N_{\frac{n}{2}%
}^{\beta |_{\frac{n}{2}}^{-1}}\right. \text{;}  \label{eqme26.5}
\end{equation}%
equations (\ref{basico}), (\ref{eqme4}) imply 
\begin{equation}
\left\{ N_{\frac{n}{2}}=N_{0}^{\beta |_{0}},\text{ }N_{0}=N_{\frac{n}{2}%
}^{\beta |_{0}^{-1}}\right. \text{;}  \label{eqme27}
\end{equation}%
equation (\ref{eqme8}) implies 
\begin{equation}
\left\{ N_{0}=N_{\frac{n}{2}}^{\beta |_{0}},\text{ }N_{\frac{n}{2}%
}=N_{0}^{\beta |_{0}^{-1}}\right. \text{;}  \label{eqme28}
\end{equation}%
equations (\ref{basico}), (\ref{eqme10}) imply 
\begin{equation}
\left\{ N_{0}=N_{\frac{n}{2}}^{\beta |_{\frac{n}{2}}},\text{ }N_{\frac{n}{2}%
}=N_{0}^{\beta |_{\frac{n}{2}}^{-1}}\right. \text{;}  \label{eqme29}
\end{equation}%
equations (\ref{basico}), (\ref{eqme17}) imply%
\begin{equation}
\left\{ N_{j}=N_{j+\frac{n}{2}}^{\beta |_{0}}\text{, }N_{j+\frac{n}{2}%
}=N_{j}^{\beta |_{0}^{-1}}\right. ,\forall j\in \{1,\ldots ,\frac{n}{2}-1\}%
\text{;}  \label{eqme30}
\end{equation}%
equation (\ref{eqme19}) imply 
\begin{equation}
\left\{ N_{j+\frac{n}{2}}=N_{j}^{\beta |_{0}},\text{ }N_{j}=N_{j+\frac{n}{2}%
}^{\beta |_{0}^{-1}}\right. ,\forall j\in \{1,\ldots ,\frac{n}{2}-1\}\text{;}
\label{eqme31}
\end{equation}%
equations (\ref{basico}) and (\ref{eqme21}) imply 
\begin{equation}
\left\{ N_{j+\frac{n}{2}}=N_{j}^{\beta |_{\frac{n}{2}}},\text{ }N_{j}=N_{j+%
\frac{n}{2}}^{\beta |_{\frac{n}{2}}^{-1}}\right. ,\forall j\in \{1,\ldots ,%
\frac{n}{2}-1\}\text{;}  \label{eqme32}
\end{equation}%
equation (\ref{eqme23}) imply 
\begin{equation}
\left\{ N_{j}=N_{j+\frac{n}{2}}^{\beta |_{\frac{n}{2}}},\text{ }N_{j+\frac{n%
}{2}}=N_{j}^{\beta |_{\frac{n}{2}}^{-1}}\right. ,\forall j\in \{1,\ldots ,%
\frac{n}{2}-1\}\text{.}  \label{eqme33}
\end{equation}%
Thus (\ref{eqme24})-(\ref{eqme33}) prove%
\begin{eqnarray*}
N &=&\left\langle N_{i}\mid i\in Y\right\rangle \\
&=&\left\langle [\beta |_{i},\tau ^{k}]\mid \forall i,k\in \mathbb{Z}%
\right\rangle
\end{eqnarray*}%
is an abelian normal subgroup of $H$.
\end{proof}

\begin{lemma}
\label{me(L2)} The group $U=\left\langle N,\text{ \ }\beta |_{j}\mid j\neq 0,%
\frac{n}{2}\right\rangle $ is a normal abelian subgroup of $H$.
\end{lemma}

\begin{proof}
Lemma \ref{me(L1)} and equations (\ref{eqme6}), (\ref{eqme12}), (\ref{eqme14}%
) and (\ref{eqme15}) show that $U$ is abelian.

The fact that $N$ is normal in $H$, together with the following assertions
prove that $U$ is normal in $H$.

Let $J=\left\langle \beta _{0},\beta _{\frac{n}{2}},\tau \right\rangle $.
Then, for $j\in Y-\{0,\frac{n}{2}\}$, we have

\begin{itemize}
\item[(I)] $\left\langle \beta |_{j}\right\rangle ^{J}\leq U:$ 
\begin{equation*}
\beta |_{j}^{\tau ^{t}}=\beta |_{j}[\beta |_{j},\tau ^{t}];
\end{equation*}%
\begin{equation*}
\beta |_{j}^{\beta |_{0}}\overset{(\ref{eqme18})}{=}\beta |_{j+\frac{n}{2}};
\end{equation*}%
\begin{equation*}
\beta |_{j}^{\beta |_{0}^{-1}}\overset{(\ref{eqme16})}{=}\tau ^{-1}\beta
|_{j+\frac{n}{2}}\tau =\beta |_{j+\frac{n}{2}}[\beta |_{j+\frac{n}{2}},\tau
];
\end{equation*}%
\begin{equation*}
\beta |_{j}^{\beta |_{\frac{n}{2}}}\overset{(\ref{eqme20})}{=}\tau
^{-1}\beta |_{j+\frac{n}{2}}\tau =\beta |_{j+\frac{n}{2}}[\beta |_{j+\frac{n%
}{2}},\tau ];
\end{equation*}%
\begin{equation*}
\beta |_{j}^{\beta |_{\frac{n}{2}}^{-1}}\overset{(\ref{eqme22})}{=}\beta
|_{j+\frac{n}{2}};
\end{equation*}

\item[(II)] $\left\langle \beta |_{j+\frac{n}{2}}\right\rangle ^{J}\leq U:$ 
\begin{equation*}
\beta |_{j+\frac{n}{2}}^{\tau ^{t}}=\beta |_{j+\frac{n}{2}}[\beta |_{j+\frac{%
n}{2}},\tau ^{t}];
\end{equation*}%
\begin{eqnarray*}
&&\beta |_{j+\frac{n}{2}}^{\beta |_{0}}\overset{(\ref{eqme16})}{=}\beta
|_{0}^{-1}\tau \beta |_{0}\beta |_{j}\beta |_{0}^{-1}\tau ^{-1}\beta |_{0} \\
&=&\left( [\beta |_{0},\tau ]^{-1}\right) ^{\tau ^{-1}}\beta |_{j}^{\tau
^{-1}}[\beta |_{0},\tau ]^{\tau ^{-1}}\in U;
\end{eqnarray*}%
\begin{equation*}
\beta |_{j+\frac{n}{2}}^{\beta |_{0}^{-1}}\overset{(\ref{eqme18})}{=}\beta
|_{j}\in U;
\end{equation*}%
\begin{equation*}
\beta |_{j+\frac{n}{2}}^{\beta |_{\frac{n}{2}}}\overset{(\ref{eqme22})}{=}%
\beta |_{j}\in U;
\end{equation*}%
\begin{eqnarray*}
&&\beta |_{j+\frac{n}{2}}^{\beta |_{\frac{n}{2}}^{-1}}\overset{(\ref{eqme20})%
}{=}\beta |_{\frac{n}{2}}\tau \beta |_{\frac{n}{2}}^{-1}\beta |_{j}\beta |_{%
\frac{n}{2}}\tau ^{-1}\beta |_{\frac{n}{2}}^{-1} \\
&=&[\beta |_{\frac{n}{2}},\tau ]^{\beta |_{\frac{n}{2}}^{-1}\tau ^{-1}}\beta
|_{j}^{\tau ^{-1}}\left( [\beta |_{\frac{n}{2}},\tau ]^{-1}\right) ^{\beta
|_{\frac{n}{2}}^{-1}\tau ^{-1}}\text{.}
\end{eqnarray*}
\end{itemize}

Hence, $U$ is a normal abelian subgroup of $H.$
\end{proof}

\begin{lemma}
\label{me(L3)} $V=\left\langle U,\text{ \ }\beta |_{\frac{n}{2}}\beta |_{0},%
\text{ \ }\tau \beta |_{0}^{2}\right\rangle $ is a normal abelian subgroup
of $H.$
\end{lemma}

\begin{proof}
Lemma \ref{me(L2)} together with the following assertions proves that $V$ is
a normal abelian subgroup of $H.$

Given $j\in Y-\{0,\frac{n}{2}\},k\in \mathbb{Z},$ and $J=\left\langle \beta
|_{0},\beta _{\frac{n}{2}},\tau \right\rangle $, we prove

\begin{itemize}
\item[(I)] $\beta |_{\frac{n}{2}}\beta |_{0}\in C_{H}(U):$

\begin{equation*}
(\beta |_{j})^{\beta |_{\frac{n}{2}}\beta |_{0}}\overset{(\ref{eqme20})}{=}%
(\beta |_{j+\frac{n}{2}})^{\tau \beta |_{0}}\overset{(\ref{eqme16})}{=}\beta
|_{j};
\end{equation*}%
\begin{equation*}
(\beta |_{j+\frac{n}{2}})^{\beta |_{\frac{n}{2}}\beta |_{0}}\overset{(\ref%
{eqme22})}{=}(\beta |_{j})^{\beta |_{0}}\overset{(\ref{eqme18})}{=}\beta
|_{j+\frac{n}{2}};
\end{equation*}%
\begin{equation*}
\lbrack \beta |_{j},\tau ^{k}]^{\beta |_{\frac{n}{2}}\beta |_{0}}=[\beta
|_{j},\tau ^{k}]^{\beta |_{\frac{n}{2}}\tau ^{-1}\tau \beta |_{0}}\overset{(%
\ref{eqme21})}{=}[\beta |_{j+\frac{n}{2}},\tau ^{k}]^{\tau \beta |_{0}}
\end{equation*}%
\begin{equation*}
\overset{(\ref{eqme17})}{=}[\beta |_{j},\tau ^{k}];
\end{equation*}%
\begin{equation*}
\lbrack \beta |_{j+\frac{n}{2}},\tau ^{k}]^{\beta |_{\frac{n}{2}}\beta |_{0}}%
\overset{(\ref{eqme23})}{=}[\beta |_{j},\tau ^{k}]^{\beta |_{0}}\overset{(%
\ref{eqme19})}{=}[\beta |_{j+\frac{n}{2}},\tau ^{k}];
\end{equation*}%
\begin{equation*}
\lbrack \beta |_{0},\tau ^{k}]^{\beta |_{\frac{n}{2}}\beta |_{0}}\overset{(%
\ref{eqme2})}{=}[\beta |_{\frac{n}{2}},\tau ^{k}]^{\beta |_{0}}\overset{(\ref%
{eqme8})}{=}[\beta |_{0},\tau ^{k}];
\end{equation*}%
\begin{eqnarray*}
\lbrack \beta |_{\frac{n}{2}},\tau ^{k}]^{\beta |_{\frac{n}{2}}\beta |_{0}}
&=&[\beta |_{\frac{n}{2}},\tau ^{k}]^{\beta |_{\frac{n}{2}}\tau ^{-1}\tau
\beta |_{0}} \\
&&\overset{(\ref{eqme10})}{=}[\beta |_{0},\tau ^{k}]^{\tau \beta |_{0}}%
\overset{(\ref{eqme4})}{=}[\beta |_{\frac{n}{2}},\tau ^{k}];
\end{eqnarray*}

\item[(II)] $\tau \beta |_{0}^{2}\in C_{H}(U):$

\begin{equation*}
\beta |_{j}^{\tau \beta |_{0}^{2}}=(\beta |_{j}[\beta |_{j},\tau ])^{\beta
|_{0}^{2}}=(\beta |_{j}^{\beta |_{0}}[\beta |_{j},\tau ]^{\beta
|_{0}})^{\beta |_{0}}
\end{equation*}%
\begin{equation*}
\overset{(\ref{eqme18}),(\ref{eqme19})}{=}(\beta |_{j+\frac{n}{2}}[\beta
|_{j+\frac{n}{2}},\tau ])^{\beta |_{0}}=\beta |_{j+\frac{n}{2}}^{\tau \beta
|_{0}}\overset{(\ref{eqme16})}{=}\beta |_{j};
\end{equation*}%
\begin{equation*}
(\beta |_{j+\frac{n}{2}})^{\tau \beta |_{0}^{2}}\overset{(\ref{eqme16})}{=}%
\beta |_{j}^{\beta |_{0}}\overset{(\ref{eqme18})}{=}\beta |_{j+\frac{n}{2}};
\end{equation*}%
\begin{equation*}
\lbrack \beta |_{0},\tau ^{k}]^{\tau \beta |_{0}^{2}}\overset{(\ref{eqme4})}{%
=}[\beta |_{\frac{n}{2}},\tau ^{k}]^{\beta |_{0}}\overset{(\ref{eqme8})}{=}%
[\beta |_{0},\tau ^{k}];
\end{equation*}%
\begin{eqnarray*}
&&\lbrack \beta |_{\frac{n}{2}},\tau ^{k}]^{\tau \beta |_{0}^{2}}\overset{(%
\ref{basico})}{=}([\beta |_{\frac{n}{2}},\tau ]^{-1}[\beta |_{\frac{n}{2}%
},\tau ^{k+1}])^{\beta |_{0}^{2}} \\
&&\overset{(\ref{eqme8})}{=}([\beta |_{0},\tau ]^{-1}[\beta |_{0},\tau
^{k+1}])^{\beta |_{0}}
\end{eqnarray*}%
\begin{equation*}
\overset{(\ref{basico})}{=}[\beta |_{0},\tau ^{k}]^{\tau \beta |_{0}}\overset%
{(\ref{eqme4})}{=}[\beta |_{\frac{n}{2}},\tau ^{k}];
\end{equation*}%
\begin{equation*}
\lbrack \beta |_{j},\tau ^{k}]^{\tau \beta |_{0}^{2}}\overset{(\ref{basico})}%
{=}([\beta |_{j},\tau ]^{-1}[\beta |_{j},\tau ^{k+1}])^{\beta |_{0}^{2}}
\end{equation*}%
\begin{equation*}
\overset{(\ref{eqme19})}{=}([\beta |_{j+\frac{n}{2}},\tau ]^{-1}[\beta |_{j+%
\frac{n}{2}},\tau ^{k+1}])^{\beta |_{0}}
\end{equation*}%
\begin{equation*}
\overset{(\ref{basico})}{=}[\beta |_{j+\frac{n}{2}},\tau ^{k}]^{\tau \beta
|_{0}}\overset{(\ref{eqme17})}{=}[\beta |_{j},\tau ^{k}];
\end{equation*}%
\begin{equation*}
\lbrack \beta |_{j+\frac{n}{2}},\tau ^{k}]^{\tau \beta |_{0}^{2}}\overset{(%
\ref{eqme17})}{=}[\beta |_{j},\tau ^{k}]^{\beta |_{0}}\overset{(\ref{eqme19})%
}{=}[\beta |_{j+\frac{n}{2}},\tau ^{k}];
\end{equation*}

\item[(III)] $\tau \beta |_{0}^{2}\in C_{H}(\beta |_{\frac{n}{2}}\beta
|_{0}):$

\begin{eqnarray*}
(\beta |_{\frac{n}{2}}\beta |_{0})^{\tau \beta |_{0}^{2}} &=&\beta
|_{0}^{-2}\tau ^{-1}\beta |_{\frac{n}{2}}\beta |_{0}\tau \beta |_{0}^{2} \\
&&\overset{(\ref{eqme3})}{=}\beta |_{0}^{-2}\tau ^{-1}\beta |_{\frac{n}{2}%
}\beta |_{\frac{n}{2}}\tau ^{-1}\beta |_{\frac{n}{2}}\beta |_{0}
\end{eqnarray*}%
\begin{eqnarray*}
&=&\beta |_{0}^{-2}\tau ^{-1}\beta |_{\frac{n}{2}}^{2}\tau ^{-1}\beta |_{%
\frac{n}{2}}\beta |_{0}=(\tau \beta |_{0}^{2})^{-1}\beta |_{\frac{n}{2}%
}^{2}\tau ^{-1}\beta |_{\frac{n}{2}}\beta |_{0} \\
&&\overset{(\ref{eqme9})}{=}\beta |_{\frac{n}{2}}\beta |_{0};
\end{eqnarray*}
.

\item[(IV)] $\left\langle \tau \beta |_{0}^{2}\right\rangle ^{J}\leq V:$

\begin{equation*}
(\tau \beta |_{0}^{2})^{\tau ^{k}}=\tau (\beta |_{0}^{2})^{\tau ^{k}}=\tau
\beta |_{0}^{2}[\beta |_{0}^{2},\tau ^{k}]=\tau \beta |_{0}^{2}[\beta
|_{0},\tau ^{k}]^{\beta |_{0}}[\beta |_{0},\tau ^{k}];
\end{equation*}%
\begin{equation*}
(\tau \beta |_{0}^{2})^{\beta |_{0}}=\beta |_{0}^{-1}\tau \beta
|_{0}^{2}\beta |_{0}=\tau \tau ^{-1}\beta |_{0}^{-1}\tau \beta |_{0}\beta
|_{0}^{2}=\tau \lbrack \tau ,\beta |_{0}]\beta |_{0}^{2}
\end{equation*}%
\begin{equation*}
=\tau \lbrack \tau ,\beta |_{0}]\tau ^{-1}\tau \beta |_{0}^{2}=([\beta
|_{0},\tau ]^{-1})^{\tau ^{-1}}\tau \beta |_{0}^{2};
\end{equation*}%
\begin{equation}  \label{*(p)}
(\tau \beta |_{0}^{2})^{\beta |_{0}^{-1}}=\beta |_{0}\tau \beta |_{0}=\tau
\beta |_{0}[\beta |_{0},\tau ]\beta |_{0}=\tau \beta |_{0}^{2}[\beta
|_{0},\tau ]^{\beta |_{0}};
\end{equation}%
\begin{equation*}
(\tau \beta |_{0}^{2})^{\beta |_{\frac{n}{2}}^{-1}}\overset{(\ref{*(p)})}{=}%
\left( (\tau \beta |_{0}^{2})^{\beta |_{0}^{-1}}([\beta |_{0},\tau
]^{-1})^{\beta |_{0}}\right) ^{\beta |_{\frac{n}{2}}^{-1}}
\end{equation*}%
\begin{equation*}
=(\tau \beta |_{0}^{2})^{\beta |_{0}^{-1}\beta |_{\frac{n}{2}}^{-1}}([\beta
|_{0},\tau ]^{-1})^{\beta |_{0}\beta |_{\frac{n}{2}}^{-1}}
\end{equation*}%
\begin{equation*}
=(\tau \beta |_{0}^{2})^{(\beta |_{\frac{n}{2}}\beta |_{0})^{-1}}([\beta
|_{0},\tau ]^{-1})^{\beta |_{0}\beta |_{\frac{n}{2}}^{-1}}\overset{(III)}{=}%
\tau \beta |_{0}^{2}([\beta |_{0},\tau ]^{-1})^{\beta |_{0}\beta |_{\frac{n}{%
2}}^{-1}};
\end{equation*}%
\begin{equation*}
(\tau \beta |_{0}^{2})^{\beta |_{\frac{n}{2}}} =
(\tau\beta|_{0}^{2})^{\beta|_{\frac{n}{2}}\beta|_{0}\beta|_{0}^{-1}}\overset{%
(III)}{=} (\tau\beta|_{0}^{2})^{\beta|_{0}^{-1}} \overset{(\ref{*(p)})}{=}
\tau\beta|_{0}^{2}[\beta|_{0}, \tau]^{\beta|_{0}}.
\end{equation*}

\item[(V)] $\left\langle \beta |_{\frac{n}{2}}\beta |_{0}\right\rangle
^{J}\leq V:$

\begin{equation*}
(\beta |_{\frac{n}{2}}\beta |_{0})^{\tau ^{k}}=\beta |_{\frac{n}{2}}\beta
|_{0}[\beta |_{\frac{n}{2}}\beta |_{0},\tau ^{k}]=\beta |_{\frac{n}{2}}\beta
|_{0}[\beta |_{\frac{n}{2}},\tau ^{k}]^{\beta |_{0}}[\beta |_{0},\tau ^{k}];
\end{equation*}%
\begin{equation}  \label{*(t)}
(\beta |_{\frac{n}{2}}\beta |_{0})^{\beta |_{0}}=\beta |_{0}^{-1}\beta |_{%
\frac{n}{2}}\beta |_{0}^{2}=\beta |_{0}^{-1}\beta |_{\frac{n}{2}}\tau
^{-1}\tau \beta |_{0}^{2}
\end{equation}%
\begin{equation*}
=\beta |_{0}^{-1}\beta |_{\frac{n}{2}}^{-1}\beta |_{\frac{n}{2}}^{2}\tau
^{-1}\tau \beta |_{0}^{2}\overset{(\ref{eqme9})}{=}(\beta |_{\frac{n}{2}%
}\beta |_{0})^{-1}(\tau \beta |_{0}^{2})^{2};
\end{equation*}%
\begin{equation}  \label{*(u)}
\beta |_{\frac{n}{2}}\beta |_{0}\overset{(\ref{*(t)})}{=}(\tau \beta
|_{0}^{2})^{2}((\beta |_{\frac{n}{2}}\beta |_{0})^{-1})^{\beta |_{0}};
\end{equation}%
\begin{equation*}
(\beta |_{\frac{n}{2}}\beta |_{0})^{\beta |_{0}^{-1}}\overset{(\ref{*(u)})}{=%
}((\tau \beta |_{0}^{2})^{2})^{\beta |_{0}^{-1}}(\beta |_{\frac{n}{2}}\beta
|_{0})^{-1};
\end{equation*}%
\begin{equation}  \label{*(x)}
(\beta |_{\frac{n}{2}}\beta |_{0})^{\beta |_{\frac{n}{2}}^{-1}}=\beta |_{%
\frac{n}{2}}^{2}\beta |_{0}\beta |_{\frac{n}{2}}^{-1}=\beta |_{\frac{n}{2}%
}^{2}\tau ^{-1}\tau \beta |_{0}\beta |_{0}\beta |_{0}^{-1}\beta |_{\frac{n}{2%
}}^{-1}
\end{equation}%
\begin{equation*}
\overset{(\ref{eqme9})}{=}(\tau \beta |_{0}^{2})^{2}\beta |_{0}^{-1}\beta |_{%
\frac{n}{2}}^{-1}=(\tau \beta |_{0}^{2})^{2}(\beta |_{\frac{n}{2}}\beta
|_{0})^{-1};
\end{equation*}%
\begin{equation*}
(\beta |_{\frac{n}{2}}\beta |_{0})^{\beta |_{\frac{n}{2}}}\overset{(\ref%
{*(x)})}{=}(\beta |_{\frac{n}{2}}\beta |_{0})^{-1}((\tau \beta
|_{0}^{2})^{2})^{\beta |_{\frac{n}{2}}}
\end{equation*}
\end{itemize}
\end{proof}

\section{\bf Solvable groups for $n=4$}

\vskip 0.4 true cm

Let $B$ be an abelian subgroup of $\mathcal{A}_{4}=Aut(T_{4})$ normalized by 
$\tau $ and let $\beta \in B$. Then, by Proposition \ref{casop^k_1} , $%
\sigma _{\beta }\in D=\left\langle (0,1,2,3),(0,2)\right\rangle $, the
unique Sylow $2$-subgroup of $\Sigma _{4}$ which contains $\sigma =\sigma
_{\tau }=(0,1,2,3)$.

The normalizer of $\overline{\left\langle \tau \right\rangle }$ here is $%
\Gamma _{0}=N_{\mathcal{A}_{4}}\left( \overline{\left\langle \tau
\right\rangle }\right) =\left\langle \Psi ,\iota \right\rangle $ where $%
\Psi$ is the monic normalizer and where $\iota =\iota ^{\left( 1\right)
}\left( 0,3\right) \left( 1,2\right) $ inverts $\tau $.

Given a group $W$, the subgroup generated by  squares of its elements is
denoted by $W^{2}$.

\begin{lemma}
\label{xcaso4_1}Let $L=L\left( D\right) $ be the layer closure of $D$ above.
If $\gamma \in L^{2}$ then $\gamma \tau $ is conjugate to $\tau $.
\end{lemma}

\begin{proof}
If $\alpha \in L$ then $\sigma _{\alpha ^{2}}\in \left\langle \sigma
^{2}\right\rangle $ and the product in some order of the states $\left(
\alpha ^{2}\right) |_{i}$ $\left( 0\leq i\leq 3\right) $ belongs to $S=L^{2}$%
.

Let $\gamma \in S$. Then $\gamma \tau $ is transitive on the $1st$ level of
the tree and $(\gamma \tau )^{4}$ is inactive with conjugate $1$st level
states, where the first state is%
\begin{equation*}
\left( \gamma |_{0}\right) \left( \gamma |_{1}\right) \left( \gamma
|_{2}\right) \left( \gamma |_{3}\right) \tau \text{ if }\sigma _{\gamma }=e,
\end{equation*}%
and 
\begin{equation*}
\left( \gamma |_{0}\right) \left( \gamma |_{3}\right) \left( \gamma
|_{2}\right) \left( \gamma |_{1}\right) \tau \text{ if }\sigma _{\gamma
}=\sigma ^{2}\text{;}
\end{equation*}%
in both cases the element is contained in $S\tau $. Therefore, $\gamma
\tau $ is transitive on the $2$nd level of the tree. Now use induction to
prove that $\gamma \tau $ is transitive on all levels of the tree. As $\gamma\tau $ is transitive on all levels of the tree, then  $\gamma\tau$ is conjugate to $\tau.$ 
\end{proof}

\subsection{ Cases $\protect\sigma _{\protect\beta }\in
\{(0,3)(1,2),(0,1)(2,3)\}$}

\label{caso4_1} We will show that these cases cannot occur. We note that $%
\sigma _{\tau }$ conjugates $(0,1)(2,3)$ to $(0,3)(1,2)$. Since the argument
for $\beta $ applies to $\beta ^{\tau }$, it is sufficient to consider the
first case.

Suppose $\sigma _{\beta }=(0,1)(2,3)$. Then, 
\begin{equation*}
\beta ^{\tau }=\left( \tau ^{-1}\left( \beta |_{3}\right) ,\beta |_{0},\beta
|_{1},\beta |_{2}\tau \right) \left( \sigma _{\beta }\right) ^{\sigma _{\tau
}}\text{.}
\end{equation*}

On substituting $\alpha =\beta ^{\tau }$ in $\theta =[\beta ,\alpha ]$ and
in ( \ref{eq6}) 
\begin{equation}
\theta |_{\left( i\right) \sigma _{\alpha \beta }}=\left( \beta |_{\left(
i\right) \sigma _{\alpha }}\right) ^{-1}\left( \alpha |_{i}\right)
^{-1}\left( \beta |_{i}\right) \left( \alpha |_{\left( i\right) \sigma
_{\beta }}\right) ,\forall i\in Y\text{.}
\end{equation}%
we get $\theta =e$ and 
\begin{equation}
e=\left( \beta |_{\left( i\right) \sigma _{\beta ^{\tau }}}\right)
^{-1}\left( \beta ^{\tau }|_{i}\right) ^{-1}\left( \beta |_{i}\right) \left(
\beta ^{\tau }|_{\left( i\right) \sigma _{\beta }}\right) ,\forall i\in Y
\end{equation}%
and so for the index $i=0$, we obtain%
\begin{eqnarray*}
e &=&\left( \beta |_{3}\right) ^{-1}\left( \tau ^{-1}\left( \beta
|_{3}\right) \right) ^{-1}\left( \beta |_{0}\right) \left( \beta
|_{0}\right) , \\
e &=&\left( \beta |_{3}\right) ^{-2}\tau \left( \beta |_{0}\right) ^{2}
\end{eqnarray*}%
which is impossible.

\subsection{Cases $\protect\sigma _{\protect\beta }\in \{(0,2),(1,3)\}$}

\label{caso4_2}

\begin{lemma}
\label{xcaso4_2} Let $\alpha ,\gamma \in Aut(T_{4})$ be such that%
\begin{eqnarray*}
\sigma _{\alpha },\sigma _{\gamma } &\in &\left\langle
(0,1,2,3),(0,2)\right\rangle , \\
\tau ^{-1}\alpha ^{2} &=&\gamma ^{2}\tau , \\
\lbrack \alpha ,\tau ^{k}]^{\gamma } &=&[\gamma ,\tau ^{k}]
\end{eqnarray*}%
for all $k\in \mathbb{Z}$. Then, 
\begin{equation*}
\sigma _{\alpha },\sigma _{\gamma }\in \left\langle \sigma \right\rangle 
\text{, \ \ \ \ }\sigma _{\alpha }\sigma _{\gamma }=\sigma ^{\pm 1}\text{.}
\end{equation*}
\end{lemma}

\begin{proof}
From the second and third equations above, we have $\sigma ^{-1}\sigma
_{\alpha }^{2}=\sigma _{\gamma }^{2}\sigma $ and $[\sigma _{\alpha },\sigma
^{k}]^{\sigma _{\gamma }}=[\sigma _{\gamma },\sigma ^{k}]$.

(i) Suppose $\sigma _{\gamma }^{2}=e$. Then $\sigma _{\alpha }^{2}=\sigma
^{2}\ $and therefore, $\sigma _{\alpha }=\sigma ^{\pm 1}$, $[\sigma _{\alpha
},\sigma ^{k}]^{\sigma _{\gamma }}=[\sigma _{\gamma },\sigma ^{k}]=e$ for
all $k$; thus, $\sigma _{\gamma }\in \left\langle \sigma \right\rangle $ and 
$\sigma _{\gamma }\in \left\langle \sigma ^{2}\right\rangle $, $\sigma
_{\alpha }\sigma _{\gamma }=\sigma ^{\pm 1}$ follows.

(ii) Suppose $o(\sigma _{\gamma })=4$. Then, $\sigma _{\gamma }=\sigma ^{\pm
1}\ $and $\sigma _{\alpha }^{2}=e$. Since $[\sigma _{\alpha },\sigma
^{k}]^{\sigma _{\gamma }}=e$ for all $k$, we obtain $\sigma _{\alpha }\in
\left\langle \sigma \right\rangle $, $\sigma _{\alpha }^{2}=e$ and $\sigma
_{\alpha }\in \left\langle \sigma ^{2}\right\rangle .$Therefore, $\sigma
_{\alpha }\sigma _{\gamma }=\sigma ^{\pm 1}$.
\end{proof}

(1) Suppose $\sigma _{\beta }=(0,2)$. Then by the analysis in Section \ref%
{2ciclo}, we conclude 
\begin{equation*}
V=\left\langle [\beta |_{i},\tau ^{k}],\beta |_{1},\beta |_{3},\beta
|_{2}\beta |_{0},\tau \beta |_{0}^{2}\mid i\in Y, k \in \mathbb{Z}\right\rangle
\end{equation*}%
is an abelian normal subgroup of $H.$

By Lemma \ref{xcaso4_1} , $\tau \beta |_{0}^{2}=\mu $ is a conjugate of $%
\tau $. As $V$ is abelian, there exist $\xi ,t_{1},t_{2}\in \mathbb{Z}_{4}$
such that 
\begin{equation*}
\mu =\tau \beta |_{0}^{2},\beta |_{2}\beta |_{0}=\mu ^{\xi },\beta |_{1}=\mu
^{t_{1}},\beta |_{3}=\mu ^{t_{2}}\text{.}
\end{equation*}%
Therefore, 
\begin{equation*}
\beta |_{2}=\mu ^{\xi }\beta |_{0}^{-1},\tau =\mu \beta |_{0}^{-2}\text{.}
\end{equation*}%

On substituting $\gamma =\beta |_{0}$ and $\alpha =\beta |_{2}$ in (\ref{eqme8}) and (\ref{eqme9}), by  Lemma \ref%
{xcaso4_2}, we obtain $\sigma _{\alpha \gamma }=\sigma _{\beta |_{2}\beta
|_{0}}=\sigma ^{\pm 1}$. Thus, from $\beta |_{2}\beta |_{0}=\mu ^{\xi }$, we
reach $\xi \in U(\mathbb{Z}_{4})$.{\Large \ }

By (\ref{eqme9}), we have 
\begin{equation*}
\beta |_{2}^{2}\tau ^{-1}=\tau \beta |_{0}^{2}\text{.}
\end{equation*}%
It follows then that%
\begin{eqnarray*}
\mu ^{\xi }\beta |_{0}^{-1}\mu ^{\xi }\beta |_{0}^{-1}\beta |_{0}^{2}\mu
^{-1} &=&\mu , \\
\left( \mu ^{\xi }\right) ^{\beta |_{0}} &=&\mu ^{2-\xi }\text{.}
\end{eqnarray*}%
Therefore, 
\begin{equation}
\mu ^{\beta |_{0}}=\mu ^{\frac{2-\xi }{\xi }}
\end{equation}%
where $\frac{2-\xi }{\xi }\in \mathbb{Z}_{4}^{1}.$

By Equation (\ref{eqme18}) we have

\begin{equation*}
\beta |_{1}^{\beta |_{0}}=\beta |_{3}\text{.}
\end{equation*}%
It follows that 
\begin{equation*}
\left( \mu ^{t_{1}}\right) ^{\beta |_{0}}=\mu ^{t_{2}},\text{ }\mu ^{t_{1}%
\frac{2-\xi }{\xi }}=\mu ^{t_{2}},\text{ }t_{2}=t_{1}\frac{2-\xi }{\xi }%
\text{.}
\end{equation*}

We have reached the form of $\beta $,

\begin{equation*}
\beta =(\beta |_{0},\mu ^{t_{1}},\mu ^{\xi }\beta |_{0}^{-1},\mu ^{t_{1}%
\frac{2-\xi }{\xi }})(0,2)
\end{equation*}%
where $\mu =\tau ^{\alpha }$ for some $\alpha \in Aut(T_{4}).$

Since $\mu^{\beta|_{0}} = \mu^{\frac{2- \xi}{\xi}}, $ we have 
$\displaystyle 
\beta |_{0}=\left( \lambda _{\frac{2-\xi }{\xi }}\tau ^{m}\right) ^{\alpha }
$ 
for some $m\in \mathbb{Z}_{4}$. 

Hence,  
\begin{equation*}
\mu ^{t_{1}}=(\tau ^{t_{1}})^{\alpha },
\end{equation*}%
\begin{equation*}
\begin{array}{ll}
\mu ^{\xi }\beta |_{0}^{-1} & =\left( \tau ^{\xi }\left( \lambda _{\frac{%
2-\xi }{\xi }}\tau ^{m}\right) ^{-1}\right) ^{\alpha } \\ 
& =\left( \lambda _{\frac{\xi }{2-\xi }}\tau ^{(\xi -m)\frac{\xi }{2-\xi }%
}\right) ^{\alpha }\text{.}%
\end{array}%
\end{equation*}%
Thus 
\begin{equation*}
\beta =(\lambda _{\frac{2-\xi }{\xi }}\tau ^{m},\tau ^{t_{1}},\lambda _{%
\frac{\xi }{2-\xi }}\tau ^{(\xi -m)\frac{\xi }{2-\xi }},\tau ^{t_{1}\frac{%
2-\xi }{\xi }})^{\alpha ^{(1)}}(0,2)
\end{equation*}%
and 
\begin{equation*}
\begin{array}{ll}
\tau & =\mu \beta |_{0}^{-2} \\ 
& =\left( \tau \left( \lambda _{\frac{2-\xi }{\xi }}\tau ^{m}\right)
^{-2}\right) ^{\alpha } \\ 
& =\left( \lambda _{(\frac{\xi }{2-\xi })^{2}}\tau ^{\left( 1-\frac{2m}{\xi }%
\right) \left( \frac{\xi }{2-\xi }\right) ^{2}}\right) ^{\alpha }%
\end{array}%
\end{equation*}

We note that in case $\xi =1$ and $\beta $ has the form

\begin{equation*}
\beta =(\tau ^{m},\tau ^{t_{1}},\tau ^{1-m},\tau ^{t_{1}})^{\alpha
^{(1)}}(0,2)
\end{equation*}%
where $\tau =\left( \tau ^{1-2m}\right) ^{\alpha }$; therefore, 
\begin{equation*}
\beta =(\tau ^{\frac{m}{1-2m}},\tau ^{\frac{t_{1}}{1-2m}},\tau ^{\frac{1-m}{%
1-2m}},\tau ^{\frac{t_{1}}{1-2m}})(0,2)\text{.}
\end{equation*}

(2) Suppose $\sigma _{\beta }=(1,3)$. Then, $\gamma =\beta ^{\tau }$
satisfies $[\gamma ,\gamma ^{\tau ^{k}}]=e$. Therefore, the previous case
applies and

\begin{equation*}
\gamma =(\lambda _{\frac{2-\xi }{\xi }}\tau ^{m},\tau ^{t_{1}},\lambda _{%
\frac{\xi }{2-\xi }}\tau ^{(\xi -m)\frac{\xi }{2-\xi }},\tau ^{t_{1}\frac{%
2-\xi }{\xi }})^{\alpha ^{(1)}}(0,2),
\end{equation*}%
where 
\begin{equation*}
\tau =\left( \lambda _{(\frac{\xi }{2-\xi })^{2}}\tau ^{\left( 1-\frac{2m}{%
\xi }\right) \left( \frac{\xi }{2-\xi }\right) ^{2}}\right) ^{\alpha
}=(e,e,e,\left( \lambda _{(\frac{\xi }{2-\xi })^{2}}\tau ^{\left( 1-\frac{2m%
}{\xi }\right) \left( \frac{\xi }{2-\xi }\right) ^{2}}\right) ^{\alpha
})\sigma _{\tau }\text{.}
\end{equation*}

Hence, $\beta $ has the form 
\begin{equation*}
\beta =\gamma ^{\tau ^{-1}}=(\tau ^{t_{1}},\lambda _{\frac{2-\xi }{\xi }%
}\tau ^{1+m-\xi },\tau ^{t_{1}\frac{2-\xi }{\xi }},\lambda _{\frac{\xi }{%
2-\xi }}\tau ^{(1-m)\frac{\xi }{2-\xi }})^{\alpha ^{(1)}}(1,3)\text{.}
\end{equation*}

\subsection{The case $\protect\sigma _{\protect\beta }=\left( \protect\sigma %
_{\protect\tau }\right) ^{2}=\left( 0,2\right) \left( 1,3\right) $}

\label{caso4_3}

We know that 
\begin{equation*}
V=\left\langle N,\beta |_{i}\beta |_{i+2},\beta |_{j}^{2}\tau ^{-\Delta
(j,j+2)}\mid i,j\in Y\text{ and }k\in \mathbb{Z}\right\rangle
\end{equation*}%
is an abelian normal subgroup of $H$ and

\begin{equation}
\tau ^{\Delta (i,j)}\beta |_{i+2}\beta |_{j}\tau ^{\Delta (i,j)}=\beta
|_{j+2}\beta |_{i}\text{,}  \label{base}
\end{equation}%
by analysis of the case \ref{casosigma2}.

From Lemmas \ref{xcaso4_1} and \ref{xcaso4_2}, we have 
\begin{equation*}
\tau \beta |_{0}^{2}=\mu \text{, }\beta |_{2}\beta |_{0}=\mu ^{\xi _{0}}%
\text{, }\beta |_{3}\beta |_{1}=\mu ^{\xi _{1}}\text{, }\tau \beta
|_{1}^{2}=\mu ^{\xi _{2}}
\end{equation*}%
where $\mu =\tau ^{\alpha }$ and $\xi _{0},\xi _{1},\xi _{2}\in U(\mathbb{Z}%
_{4})$. Therefore, 
\begin{equation}
\tau =\mu \beta |_{0}^{-2}  \label{sigma2_1}
\end{equation}%
\begin{equation}
\beta |_{2}=\mu ^{\xi _{0}}\beta |_{0}^{-1}  \label{sigma2_2}
\end{equation}%
\begin{equation}
\beta |_{3}=\mu ^{\xi _{1}}\beta |_{1}^{-1}  \label{sigma2_3}
\end{equation}%
\begin{equation}
\tau =\mu ^{\xi _{2}}\beta |_{1}^{-2}\text{.}  \label{sigma2_4}
\end{equation}

Now, we let $i,j$ take their values from $Y$ in (\ref{base}). Note that $%
\left( i,j\right) $ and $\left( j,i\right) $ produce equivalent equations
and the case where $i=j$ is a tautology. Thus we have to treat the cases $%
\left( i,j\right) =\left( 0,1\right) ,\left( 0,2\right) ,\left( 1,3\right)
,\left( 2,3\right) ,\left( 0,3\right) ,\left( 1,2\right) $. Indeed, the last
two cases turn out to be superfluous.

(i) Substitute $i=0,j=2$ in (\ref{base}), to obtain 
\begin{equation}
\beta |_{2}^{2}\tau ^{-1}=\tau \beta |_{0}^{2}  \label{base1}
\end{equation}

Use (\ref{sigma2_1}) and (\ref{sigma2_2}) in (\ref{base1}) to get

\begin{equation*}
\mu ^{\xi _{0}}\beta |_{0}^{-1}\mu ^{\xi _{0}}\beta |_{0}^{-1}\beta
|_{0}^{2}\mu ^{-1}=\mu
\end{equation*}%
and so,

\begin{equation*}
(\mu ^{\xi _{0}})^{\beta |_{0}}=\mu ^{2-\xi _{0}}\text{.}
\end{equation*}%
Therefore, 
\begin{equation}
\mu ^{\beta |_{0}}=\mu ^{\frac{2-\xi _{0}}{\xi _{0}}}  \label{sigma2_5}
\end{equation}%
Since $\frac{2-\xi _{0}}{\xi _{0}}\in \mathbb{Z}_{4}^{1}$, we find 
\begin{equation}
\beta |_{0}=\left( \lambda _{\frac{2-\xi _{0}}{\xi _{0}}}\tau
^{m_{0}}\right) ^{\alpha }\text{.}  \label{sigma2_6}
\end{equation}%
From (\ref{sigma2_2}),

\begin{equation}
\beta |_{2}=\mu ^{\xi _{0}}\beta |_{0}^{-1}=\left( \tau ^{\xi _{0}}\tau
^{-m_{0}}\lambda _{\frac{\xi _{0}}{2-\xi _{0}}}\right) ^{\alpha }=\left(
\lambda _{\frac{\xi _{0}}{2-\xi _{0}}}\tau ^{(\xi _{0}-m_{0})\frac{\xi _{0}}{%
2-\xi _{0}}}\right) ^{\alpha }\text{.}  \label{sigma2_7}
\end{equation}

(ii) Substitute $i=1,j=3$ in (\ref{base}) to get 
\begin{equation}
\beta |_{3}^{2}\tau ^{-1}=\tau \beta |_{1}^{2}\text{.}  \label{base2}
\end{equation}

On using (\ref{sigma2_3}) and (\ref{sigma2_4}) in (\ref{base2}), we obtain

\begin{equation*}
\mu ^{\xi _{1}}\beta |_{1}^{-1}\mu ^{\xi _{1}}\beta |_{1}^{-1}\beta
|_{1}^{2}\mu ^{-\xi _{2}}=\mu ^{\xi _{2}}
\end{equation*}%
and so,

\begin{equation*}
(\mu ^{\xi _{1}})^{\beta |_{1}}=\mu ^{2\xi _{2}-\xi _{1}}\text{.}
\end{equation*}%
Therefore, 
\begin{equation}
\mu ^{\beta |_{1}}=\mu ^{\frac{2\xi _{2}-\xi _{1}}{\xi _{1}}\text{ }}\text{.}
\label{sigma2_8}
\end{equation}%
Since $\frac{2\xi _{2}-\xi _{1}}{\xi _{1}}\in \mathbb{Z}_{4}^{1}$, we have

\begin{equation}
\beta |_{1}=\left( \lambda _{\frac{2\xi _{2}-\xi _{1}}{\xi _{1}}}\tau
^{m_{1}}\right) ^{\alpha }\text{.}  \label{sigma2_9}
\end{equation}%
By (\ref{sigma2_3}), we find

\begin{equation}
\beta |_{3}=\mu ^{\xi _{1}}\beta |_{1}^{-1}=\left( \tau ^{\xi _{1}}\tau
^{-m_{1}}\lambda _{\frac{\xi _{1}}{2\xi _{2}-\xi _{1}}}\right) ^{\alpha
}=\left( \lambda _{\frac{\xi _{1}}{2\xi _{2}-\xi _{1}}}\tau ^{(\xi
_{1}-m_{1})\frac{\xi _{1}}{2\xi _{2}-\xi _{1}}}\right) ^{\alpha }\text{.}
\label{sigma2_10}
\end{equation}

(iii) Substitute $i=0,j=1$ in (\ref{base}) to get 
\begin{equation}
\beta |_{2}\beta |_{1}=\beta |_{3}\beta |_{0}\text{.}  \label{base3}
\end{equation}

Use (\ref{sigma2_6}), (\ref{sigma2_7}), (\ref{sigma2_9}) and (\ref{sigma2_10}%
) in (\ref{base3}), to obtain 
\begin{equation*}
\lambda _{\frac{\xi _{0}}{2-\xi _{0}}}\tau ^{(\xi _{0}-m_{0})\frac{\xi _{0}}{%
2-\xi _{0}}}\lambda _{\frac{2\xi _{2}-\xi _{1}}{\xi _{1}}}\tau
^{m_{1}}=\lambda _{\frac{\xi _{1}}{2\xi _{2}-\xi _{1}}}\tau ^{(\xi
_{1}-m_{1})\frac{\xi _{1}}{2\xi _{2}-\xi _{1}}}\lambda _{\frac{2-\xi _{0}}{%
\xi _{0}}}\tau ^{m_{0}}
\end{equation*}%
and so, 
\begin{equation*}
\lambda _{\frac{\xi _{0}}{2-\xi _{0}}\frac{2\xi _{2}-\xi _{1}}{\xi _{1}}%
}\tau ^{(\xi _{0}-m_{0})\frac{\xi _{0}}{2-\xi _{0}}\frac{2\xi _{2}-\xi _{1}}{%
\xi _{1}}+m_{1}}=\lambda _{\frac{\xi _{1}}{2\xi _{2}-\xi _{1}}\frac{2-\xi
_{0}}{\xi _{0}}}\tau ^{(\xi _{1}-m_{1})\frac{\xi _{1}}{2\xi _{2}-\xi _{1}}%
\frac{2-\xi _{0}}{\xi _{0}}+m_{0}}\text{.}
\end{equation*}%
Therefore,

\begin{equation}
\left( \frac{\xi _{1}}{2\xi _{2}-\xi _{1}}\right) ^{2}=\left( \frac{\xi _{0}%
}{2-\xi _{0}}\right) ^{2}  \label{sigma2_11}
\end{equation}%
and 
\begin{equation}
(\xi _{0}-m_{0})\frac{\xi _{0}}{2-\xi _{0}}\frac{2\xi _{2}-\xi _{1}}{\xi _{1}%
}+m_{1}=(\xi _{1}-m_{1})\frac{\xi _{1}}{2\xi _{2}-\xi _{1}}\frac{2-\xi _{0}}{%
\xi _{0}}+m_{0}\text{.}  \label{sigma2_12}
\end{equation}

(iv) Substitute $i=2,j=3$ in (\ref{base}) to get 
\begin{equation}
\beta |_{0}\beta |_{3}=\beta |_{1}\beta |_{2}\text{.}  \label{base4}
\end{equation}

Use (\ref{sigma2_6}), (\ref{sigma2_7}), (\ref{sigma2_9}) and (\ref{sigma2_10}%
) in (\ref{base4}), to obtain 
\begin{equation*}
\lambda _{\frac{2-\xi _{0}}{\xi _{0}}}\tau ^{m_{0}}\lambda _{\frac{\xi _{1}}{%
2\xi _{2}-\xi _{1}}}\tau ^{(\xi _{1}-m_{1})\frac{\xi _{1}}{2\xi _{2}-\xi _{1}%
}}=\lambda _{\frac{2\xi _{2}-\xi _{1}}{\xi _{1}}}\tau ^{m_{1}}\lambda _{%
\frac{\xi _{0}}{2-\xi _{0}}}\tau ^{(\xi _{0}-m_{0})\frac{\xi _{0}}{2-\xi _{0}%
}}
\end{equation*}%
and so, 
\begin{equation*}
\lambda _{\frac{\xi _{0}}{2-\xi _{0}}\frac{\xi _{1}}{2\xi _{2}-\xi _{1}}%
}\tau ^{m_{0}\frac{\xi _{1}}{2\xi _{2}-\xi _{1}}+(\xi _{1}-m_{1})\frac{\xi
_{1}}{2\xi _{2}-\xi _{1}}}=\lambda _{\frac{2\xi _{2}-\xi _{1}}{\xi _{1}}%
\frac{\xi _{0}}{2-\xi _{0}}}\tau ^{m_{1}\frac{\xi _{0}}{2-\xi _{0}}+(\xi
_{0}-m_{0})\frac{\xi _{0}}{2-\xi _{0}}}\text{.}
\end{equation*}%
Therefore, 
\begin{equation*}
\left( \frac{\xi _{1}}{2\xi _{2}-\xi _{1}}\right) ^{2}=\left( \frac{\xi _{0}%
}{2-\xi _{0}}\right) ^{2}
\end{equation*}%
and 
\begin{equation}
m_{0}\frac{\xi _{1}}{2\xi _{2}-\xi _{1}}+(\xi _{1}-m_{1})\frac{\xi _{1}}{%
2\xi _{2}-\xi _{1}}=m_{1}\frac{\xi _{0}}{2-\xi _{0}}+(\xi _{0}-m_{0})\frac{%
\xi _{0}}{2-\xi _{0}}\text{.}  \label{sigma2_13}
\end{equation}

We have from (\ref{sigma2_11}) 
\begin{equation}
\frac{\xi _{0}}{2-\xi _{0}}=\pm \frac{\xi _{1}}{2\xi _{2}-\xi _{1}}\text{.}
\label{sigma2_14}
\end{equation}

(a) If 
\begin{equation*}
\frac{\xi _{0}}{2-\xi _{0}}=\frac{\xi _{1}}{2\xi _{2}-\xi _{1}},
\end{equation*}%
then 
\begin{equation*}
2\xi _{2}\xi _{0}-\xi _{1}\xi _{0}=2\xi _{1}-\xi _{1}\xi _{0},
\end{equation*}%
and so, 
\begin{equation}
\xi _{2}=\frac{\xi _{1}}{\xi _{0}}\text{.}  \label{sigma2_15}
\end{equation}%
From (\ref{sigma2_12}), we get 
\begin{equation}
m_{1}=\frac{\xi _{1}-\xi _{0}}{2}+m_{0}\text{.}  \label{sigma2_16}
\end{equation}

(b) If 
\begin{equation*}
\frac{\xi _{0}}{2-\xi _{0}}=-\frac{\xi _{1}}{2\xi _{2}-\xi _{1}}
\end{equation*}%
then by (\ref{sigma2_12}) and (\ref{sigma2_13}), 
\begin{equation*}
m_{0}-\xi _{0}+m_{1}=m_{1}-\xi _{1}+m_{0}
\end{equation*}%
\begin{equation*}
m_{0}+\xi _{1}-m_{1}=-m_{1}-\xi _{0}+m_{0},
\end{equation*}%
which implies $\xi _{1}=\xi _{0}=0,$which is impossible.

Now by (\ref{sigma2_15}) and (\ref{sigma2_16}), we have 
\begin{equation}
\beta |_{1}=\left( \lambda _{\frac{2-\xi _{0}}{\xi _{0}}}\tau ^{\frac{\xi
_{1}-\xi _{0}}{2}+m_{0}}\right) ^{\alpha }  \label{sigma2_17}
\end{equation}

and 
\begin{equation}
\beta |_{3}=\left( \lambda _{\frac{\xi _{0}}{2-\xi _{0}}}\tau ^{\left( \frac{%
\xi _{1}+\xi _{0}}{2}-m_{0}\right) \frac{\xi _{0}}{2-\xi _{0}}}\right)
^{\alpha }\text{.}  \label{sigma2_18}
\end{equation}

Therefore,

\begin{equation*}
\beta =(\beta |_{0},\beta |_{1},\beta |_{2},\beta |_{3})(0,2)(1,3)
\end{equation*}%
where $\beta |_{0},\beta |_{1},\beta |_{2}$ and $\beta |_{3}$ are described
in (\ref{sigma2_6}),(\ref{sigma2_17}), (\ref{sigma2_7}) and (\ref{sigma2_18}%
), respectively, and

\begin{equation*}
\begin{array}{ll}
\tau & =\mu \beta |_{0}^{-2} \\ 
& =\left( \tau \left( \lambda _{\frac{2-\xi _{0}}{\xi _{0}}}\tau
^{m_{0}}\right) ^{-2}\right) ^{\alpha } \\ 
& =\left( \lambda _{(\frac{\xi _{0}}{2-\xi _{0}})^{2}}\tau ^{\left( 1-\frac{%
2m_{0}}{\xi _{0}}\right) \left( \frac{\xi _{0}}{2-\xi _{0}}\right)
^{2}}\right) ^{\alpha }\text{.}%
\end{array}%
\end{equation*}

(v) The cases $\left( i,j\right) =\left( 1,2\right) ,\left( 0,3\right) $ in (%
\ref{base}) do not add any more information about $\beta .$

Summarizing, we have found

\begin{equation}
\beta |_{0}=\left( \lambda _{\frac{2-\xi _{0}}{\xi _{0}}}\tau
^{m_{0}}\right) ^{\alpha }\text{, }\beta |_{1}=\left( \lambda _{\frac{2-\xi
_{0}}{\xi _{0}}}\tau ^{\frac{\xi _{1}-\xi _{0}}{2}+m_{0}}\right) ^{\alpha },
\end{equation}%
\begin{equation}
\beta |_{2}=\left( \lambda _{\frac{\xi _{0}}{2-\xi _{0}}}\tau ^{(\xi
_{0}-m_{0})\frac{\xi _{0}}{2-\xi _{0}}}\right) ^{\alpha },\beta |_{3}=\left(
\lambda _{\frac{\xi _{0}}{2-\xi _{0}}}\tau ^{\left( \frac{\xi _{1}+\xi _{0}}{%
2}-m_{0}\right) \frac{\xi _{0}}{2-\xi _{0}}}\right) ^{\alpha },
\end{equation}%
\begin{equation}
\tau =\left( \lambda _{(\frac{\xi _{0}}{2-\xi _{0}})^{2}}\tau ^{\left( 1-%
\frac{2m_{0}}{\xi _{0}}\right) \left( \frac{\xi _{0}}{2-\xi _{0}}\right)
^{2}}\right) ^{\alpha }\text{.}
\end{equation}

In the particular case where $\xi _{0}=1$, $\beta $ has the form%
\begin{equation*}
\beta =(\tau ^{\frac{m_{0}}{1-2m_{0}}},\tau ^{\frac{\frac{\xi _{1}-1}{2}%
+m_{0}}{1-2m_{0}}},\tau ^{\frac{1-m_{0}}{1-2m_{0}}},\tau ^{\frac{\frac{\xi
_{1}+1}{2}-m_{0}}{1-2m_{0}}})(0,2)(1,3)
\end{equation*}

where $\tau =\left( \tau ^{1-2m_{0}}\right) ^{\alpha }$.

\subsection{Cases $\protect\sigma _{\protect\beta }\in \{e,\protect\sigma _{%
\protect\tau },\protect\sigma _{\protect\tau }^{-1}\}$}

(1) \label{caso4_5} Suppose $\sigma _{\beta }=e$ and let $\beta $ stabilize
the $k$th level of the tree. Then by Proposition \ref{inativos}, we have 
\begin{equation*}
\lbrack \beta |_{u},\beta |_{v}^{\tau ^{\xi }}]=e,\text{ for all $u,v\in 
\mathcal{M}$ with $|u|=|v|=k$.}
\end{equation*}

Therefore, $\dot{N}=\left\langle \beta |_{w}\mid |w|=k,w\in \mathcal{M}%
\right\rangle $ is abelian and so is its normal closure $\dot{M}$ under $%
\left\langle \dot{N},\tau \right\rangle $. Also, active elements in $\dot{M}$
are characterized in \ref{caso4_1}, \ref{caso4_2}, \ref{caso4_3} and \ref%
{caso4_4}. In particular, there exists $\kappa \in \dot{M}$ such that $%
\sigma _{\kappa }=(0,2)(1,3)$ and $\beta \in \times _{p^{k}}C(\kappa )$.

(2) \label{caso4_4} Suppose $\sigma _{\beta }=\sigma _{\tau }=(0,1,2,3)$.
Then, clearly the element%
\begin{equation*}
\beta ^{2}=(\beta |_{0}\beta |_{1},\text{ }\beta |_{1}\beta |_{2},\text{ }%
\beta |_{2}\beta |_{3},\text{ }\beta |_{3}\beta |_{0})(0,2)(1,3)
\end{equation*}%
satisfies $[\beta ^{2},\left( \beta ^{2}\right) ^{\tau ^{k}}]=e$ for all $%
k\in \mathbb{Z}_{4}.$ Therefore, by the previous analysis, we have

\begin{equation}
\beta |_{0}\beta |_{1}=\left( \lambda _{\frac{2-\xi _{0}}{\xi _{0}}}\tau
^{m_{0}}\right) ^{\alpha },  \label{sigma1_1}
\end{equation}

\begin{equation}
\beta |_{1}\beta |_{2}=\left( \lambda _{\frac{2-\xi _{0}}{\xi _{0}}}\tau ^{%
\frac{\xi _{1}-\xi _{0}}{2}+m_{0}}\right) ^{\alpha },  \label{sigma1_2}
\end{equation}

\begin{equation}
\beta |_{2}\beta |_{3}=\left( \lambda _{\frac{\xi _{0}}{2-\xi _{0}}}\tau
^{(\xi _{0}-m_{0})\frac{\xi _{0}}{2-\xi _{0}}}\right) ^{\alpha },
\label{sigma1_3}
\end{equation}

\begin{equation}
\beta |_{3}\beta |_{0}=\left( \lambda _{\frac{\xi _{0}}{2-\xi _{0}}}\tau
^{\left( \frac{\xi _{1}+\xi _{0}}{2}-m_{0}\right) \frac{\xi _{0}}{2-\xi _{0}}%
}\right) ^{\alpha },  \label{sigma1_4}
\end{equation}

\begin{equation}
\tau =\left( \lambda _{(\frac{\xi _{0}}{2-\xi _{0}})^{2}}\tau ^{\left( 1-%
\frac{2m_{0}}{\xi _{0}}\right) \left( \frac{\xi _{0}}{2-\xi _{0}}\right)
^{2}}\right) ^{\alpha }\text{.}  \label{sigma1_5}
\end{equation}

Hence, multiplying (\ref{sigma1_1}) by (\ref{sigma1_3}), we obtain

\begin{equation} \label{produtosbeta}
\beta |_{0}\beta |_{1}\beta |_{2}\beta |_{3}=\left( \lambda _{\frac{2-\xi
_{0}}{\xi _{0}}}\tau ^{m_{0}}\lambda _{\frac{\xi _{0}}{2-\xi _{0}}}\tau
^{(\xi _{0}-m_{0})\frac{\xi _{0}}{2-\xi _{0}}}\right) ^{\alpha }=\left( \tau
^{\frac{\xi _{0}^{2}}{2-\xi _{0}}}\right) ^{\alpha }.
\end{equation}

We define 
\begin{equation}
\psi _{\eta }=\left\{ 
\begin{array}{ll}
\lambda _{\eta }, & \text{ if }\eta \in \mathbb{Z}_{4}^{1} \\ 
\theta \lambda _{-\eta }, & \text{ if }-\eta \in \mathbb{Z}_{4}^{1}%
\end{array}%
\right. ,  \label{defpsi}
\end{equation}%
\begin{equation*}
\theta =\theta ^{(1)}(e,\tau ^{-1},\tau ^{-1},\tau ^{-1})(1,3)
\end{equation*}%
(an invertor of $\tau $) and  
 $\gamma = (e, (\beta|_{0})^{-1}, (\beta|_{0}\beta|_{1})^{-1}, (\beta|_{0}\beta|_{1}\beta|_{2})^{-1})\left(\alpha^{-1}\psi_{\frac{2 - \xi_{0}}{\xi _{0}^{2}}}\right)^{(1)}.$ 
 
 We verify, by    (\ref{produtosbeta}), that $\gamma$ conjugates
$\beta$ to \[ (e, e, e, \beta|_{0}\beta|_{1}\beta|_{2}\beta|_{3})^{\left(\alpha^{-1}\psi_{\frac{2 - \xi_{0}}{\xi _{0}^{2}}}\right)^{(1)}}\sigma \]  
which is equal to  $\tau.$

(3) Suppose $\sigma _{\beta }=\sigma _{\tau }^{-1}=(0,3,2,1)$. Then, $\beta
^{-1}$ satisfies the previous case and $\beta^{-1} = \tau^{\gamma} $ for some $\gamma \in \mathcal{A}_{4}$. Therefore, as $ \theta$ inverts $\tau $,
we have 
\begin{equation}
\beta =\left( \beta ^{-1}\right) ^{-1}=\left( \tau ^{\gamma }\right)
^{-1}=\left( \tau \right) ^{\theta \gamma }  \label{sigma-1_1}
\end{equation}%

\subsection{Final Step}

We finish the proof of the second part of Theorem A. For the case where the activity of $\beta $ is a $4$-cycle, we use the
fact that $\beta ^{2}\in B,$ which we have already described. Next, from the
description of the centralizer of $\beta ^{2}$, we are able to pin down the
form of $\beta $.

\begin{proposition}
\label{conjtau2} Let $\beta =(\beta |_{0},$ $\beta |_{1},$ $\beta |_{2},$ $%
\beta |_{3})(0,2)(1,3)$ be such that $\left( \beta |_{0}\right) \left( \beta
|_{2}\right) =\tau ^{\theta _{1}}$ and $\left( \beta |_{1}\right) \left(
\beta |_{3}\right) =\tau ^{\theta _{2}},$ for some $\theta _{1},\theta
_{2}\in Aut(T_{4})$. Then, $\beta $ is conjugate to $\tau ^{2}.$
\end{proposition}

\begin{proof}
Let $\alpha =(e,e,\beta |_{0}^{-1},\beta |_{1}^{-1})$. Then, 
\begin{equation}
\beta ^{\alpha }=(e,e,\text{ }\beta |_{0}\beta |_{2}\text{ },\beta
|_{1}\beta |_{3})(0,2)(1,3).  \label{soluvel4_1}
\end{equation}

Therefore, substituting $\beta |_{0}\beta |_{2}=\tau ^{\theta _{1}}$ and $%
\beta |_{1}\beta |_{3}=\tau ^{\theta _{2}}$ in the above equation, we have 
\begin{equation*}
\beta ^{\alpha }=(e,e,\tau ^{\theta _{1}},\tau ^{\theta _{2}})(0,2)(1,3)%
\text{.}
\end{equation*}

Conjugating $\beta ^{\alpha }$ by $\gamma =(\theta _{1}^{-1},\theta
_{2}^{-1},\theta _{1}^{-1},\theta _{2}^{-1})$ we produce 
\begin{equation*}
\beta ^{\alpha \gamma }=\tau ^{2}\text{.}
\end{equation*}
\end{proof}

We show below that active elements of $B$ produce within $B$ elements
conjugate to $\tau ^{2}$.

\begin{proposition}
Let $\beta \in B$ with nontrivial $\sigma _{\beta }$. Then

\begin{itemize}
\item[(i)] If $\sigma _{\beta }=\sigma _{\tau }^{2},$ then $\beta $ is a
conjugate of $\tau ^{2}.$

\item[(ii)] If $\sigma _{\beta }\in \{(0,2),(1,3)\},$ then $\beta \beta
^{\tau }$ is a conjugate $\tau ^{2}.$

\item[(iii)] If $\sigma _{\beta }\in \{\sigma _{\tau },\sigma _{\tau
}^{-1}\},$ then $\beta ^{2}$ is a conjugate of $\tau ^{2}.$
\end{itemize}
\end{proposition}

\begin{proof}
It is enough to prove (i), since (ii), (iii) are just special cases.

If $\sigma _{\beta }=\sigma _{\tau }^{2},$ then

\begin{equation}
\beta |_{0}=\left( \lambda _{\frac{2-\xi _{0}}{\xi _{0}}}\tau
^{m_{0}}\right) ^{\alpha }\text{, }\beta |_{1}=\left( \lambda _{\frac{2-\xi
_{0}}{\xi _{0}}}\tau ^{\frac{\xi _{1}-\xi _{0}}{2}+m_{0}}\right) ^{\alpha },
\end{equation}%
\begin{equation}
\beta |_{2}=\left( \lambda _{\frac{\xi _{0}}{2-\xi _{0}}}\tau ^{(\xi
_{0}-m_{0})\frac{\xi _{0}}{2-\xi _{0}}}\right) ^{\alpha },\beta |_{3}=\left(
\lambda _{\frac{\xi _{0}}{2-\xi _{0}}}\tau ^{\left( \frac{\xi _{1}+\xi _{0}}{%
2}-m_{0}\right) \frac{\xi _{0}}{2-\xi _{0}}}\right) ^{\alpha },
\end{equation}%
\begin{equation}
\tau =\left( \lambda _{(\frac{\xi _{0}}{2-\xi _{0}})^{2}}\tau ^{\left( 1-%
\frac{2m_{0}}{\xi _{0}}\right) \left( \frac{\xi _{0}}{2-\xi _{0}}\right)
^{2}}\right) ^{\alpha },
\end{equation}%
where $\xi _{0},\xi _{1}\in U(\mathbb{Z}_{4}),\;m_{0}\in \mathbb{Z}_{4}.$

Therefore, 
\begin{equation*}
\beta |_{0}\beta |_{2}=\left( \lambda _{\frac{2-\xi _{0}}{\xi _{0}}}\tau
^{m_{0}}\lambda _{\frac{\xi _{0}}{2-\xi _{0}}}\tau ^{(\xi _{0}-m_{0})\frac{%
\xi _{0}}{2-\xi _{0}}}\right) ^{\alpha }=\left( \tau^{\frac{\xi _{0}^{2}}{%
2-\xi _{0}}}\right)^{\alpha }=\tau^{\left(\psi_{\frac{\xi
_{0}^{2}}{2-\xi _{0}}}\right) \alpha }
\end{equation*}

\begin{equation*}
\beta |_{1}\beta |_{3}=\left( \lambda _{\frac{2-\xi _{0}}{\xi _{0}}}\tau ^{%
\frac{\xi _{1}-\xi _{0}}{2}+m_{0}}\lambda _{\frac{\xi _{0}}{2-\xi _{0}}}\tau
^{\left( \frac{\xi _{1}+\xi _{0}}{2}-m_{0}\right) \frac{\xi _{0}}{2-\xi _{0}}%
}\right) ^{\alpha }=\left( \tau ^{\frac{\xi _{1}\xi _{0}}{2-\xi _{0}}%
}\right) ^{\alpha }=\tau ^{\left(\psi _{\frac{\xi _{1}\xi _{0}}{2-\xi _{0}}}\right)\alpha
}
\end{equation*}

It follows from Proposition \ref{conjtau2}, that $\beta $ is a conjugate of $%
\tau ^{2}.$
\end{proof}

\begin{corollary}
Suppose $\beta \in $ $B$ is an active element. Then, $B$ is conjugate to a
subgroup of the centralizer $C(\tau ^{2})$.
\end{corollary}

\begin{proposition}
\label{centralizat2} Let $\gamma \in C(\tau ^{2})$. Then, 
\begin{equation}
\gamma =(\tau ^{m_{0}},\tau ^{m_{1}},\tau ^{m_{0}+\delta ((0)\sigma _{\gamma
},\text{ }2)},\tau ^{m_{1}+\delta ((1)\sigma _{\gamma },\text{ }2)})\sigma
_{\gamma },
\end{equation}%
where $m_{0},m_{1}\in \mathbb{Z}_{4},\sigma _{\gamma }\in C_{\Sigma
_{4}}(\sigma ^{2})$.
\end{proposition}

\begin{proof}
Write $\gamma =(\gamma |_{0},\gamma |_{1},\gamma |_{2},\gamma |_{3})\sigma
_{\gamma }$. Then $\tau ^{2}\gamma =\gamma \tau ^{2}$ translates to 
\begin{equation*}
\begin{array}{ll}
& (e,e,\tau ,\tau )(0,2)(1,3)(\gamma |_{0},\gamma |_{1},\gamma |_{2},\gamma
|_{3})\sigma _{\gamma } \\ 
= & (\gamma |_{0},\gamma |_{1},\gamma |_{2},\gamma |_{3})\sigma _{\gamma
}(e,e,\tau ,\tau )(0,2)(1,3),%
\end{array}%
\end{equation*}%
and this in turn translates to 

\[\begin{array}{ll} & (\gamma |_{2},\gamma |_{3},\tau \gamma |_{0},\tau \gamma
|_{1})(0,2)(1,3)\sigma _{\gamma }  \\ \\
=  & (\gamma |_{0},\gamma |_{1},\gamma |_{2},\gamma |_{3})\sigma _{\gamma }(\tau ^{\delta (0,2)},\tau ^{\delta (1,2)},\tau ^{\delta
(2,2)},\tau ^{\delta (3,2)})(0,2)(1,3)\\ \\
= &  (\gamma |_{0},\gamma |_{1},\gamma |_{2},\gamma |_{3}) 
(\tau ^{\delta ((0)\sigma _{\gamma },2)},\tau ^{\delta ((1)\sigma _{\gamma
},2)},\tau ^{\delta ((2)\sigma _{\gamma },2)},\tau ^{\delta ((3)\sigma
_{\gamma },2)})\sigma _{\gamma }(0,2)(1,3)%
\\ \\
= &   (\gamma |_{0}\tau ^{\delta ((0)\sigma _{\gamma },2)},\gamma |_{1}\tau
^{\delta ((1)\sigma _{\gamma },2)},\gamma |_{2}\tau ^{\delta ((2)\sigma
_{\gamma },2)},\gamma |_{3}\tau ^{\delta ((3)\sigma _{\gamma },2)})\sigma
_{\gamma }(0,2)(1,3)%
\end{array}
\]

%%%%%%%%%%%%%%%%%%%%%%%%
Thus, we have 
\begin{equation*}
\left\{ 
\begin{array}{l}
\gamma |_{2}=\gamma |_{0}\tau ^{\delta ((0)\sigma _{\gamma },2)}, \\ 
\gamma |_{3}=\gamma |_{1}\tau ^{\delta ((1)\sigma _{\gamma },2)}, \\ 
\tau \gamma |_{0}=\gamma |_{2}\tau ^{\delta ((2)\sigma _{\gamma },2)}, \\ 
\tau \gamma |_{1}=\gamma |_{3}\tau ^{\delta ((3)\sigma _{\gamma },2)}\text{.}%
\end{array}%
\right.
\end{equation*}%
Hence,

\begin{equation*}
\left\{ 
\begin{array}{l}
\gamma |_{2}=\gamma |_{0}\tau ^{\delta ((0)\sigma _{\gamma },2)}\text{, }%
\gamma |_{3}=\gamma |_{1}\tau ^{\delta ((1)\sigma _{\gamma },2)}\text{,} \\ 
\tau ^{\gamma |_{0}}=\tau ^{\delta ((0)\sigma _{\gamma },2)+\delta
((2)\sigma _{\gamma },2)}=\tau \text{, }\tau ^{\gamma |_{1}}=\tau ^{\delta
((1)\sigma _{\gamma },2)+\delta ((3)\sigma _{\gamma },2)}=\tau%
\end{array}%
\right. \text{.}
\end{equation*}

Therefore, there exist $m_{0},m_{1}\in \mathbb{Z}_{4}$ such that 
\begin{equation*}
\left\{ 
\begin{array}{l}
\gamma |_{0}=\tau ^{m_{0}},\text{ }\gamma |_{1}=\tau ^{m_{1}}, \\ 
\gamma |_{2}=\tau ^{m_{0}+\delta ((0)\sigma _{\gamma },2)},\text{ }\gamma
|_{3}=\tau ^{m_{1}+\delta ((1)\sigma _{\gamma },2)}%
\end{array}%
\right. \text{.}
\end{equation*}

Hence, $\gamma $ has the form 
\begin{equation}
\gamma =(\tau ^{m_{0}},\tau ^{m_{1}},\tau ^{m_{0}+\delta ((0)\sigma _{\gamma
},2)},\tau ^{m_{1}+\delta ((1)\sigma _{\gamma },2)})\sigma _{\gamma },
\label{centralizatau2}
\end{equation}%
where $\sigma _{\gamma }\in C_{\Sigma _{4}}(\sigma ^{2})$.
\end{proof}

\begin{corollary}
The centralizer of $\tau ^{2}$ in $\mathcal{A}_{4}$ is 
\begin{equation*}
C(\tau ^{2})=\left\langle (e,e,\tau ,e)(0,2),\text{ }\tau ,\text{ }(\tau
^{m_{0}},\tau ^{m_{1}},\tau ^{m_{0}},\tau ^{m_{1}})\mid m_{0},m_{1}\in 
\mathbb{Z}_{4}\right\rangle .
\end{equation*}
\end{corollary}

\begin{corollary}
\label{corolcentralizat2} Let $\gamma \in C(\tau ^{2})$ be such that $\sigma
_{\gamma }\in \left\langle (0,2)(1,3)\right\rangle $. Then 
\begin{equation*}
\gamma \in \left\langle (\tau ^{m_{0}},\tau ^{m_{1}},\tau ^{m_{0}},\tau
^{m_{1}}),\text{ }\tau ^{2}\mid m_{0},m_{1}\in \mathbb{Z}_{4}\right\rangle 
\text{.}
\end{equation*}
\end{corollary}

\begin{proposition}
Let $\dot{H}=\left\langle (\tau ^{m_{0}},\tau ^{m_{1}},\tau ^{m_{0}},\tau
^{m_{1}}),\tau ^{2}\mid m_{0},m_{1}\in \mathbb{Z}_{4}\right\rangle $. Then
the normalizer $N_{\mathcal{A}_{4}}(\dot{H})$ is the group 
\begin{equation*}
\left\langle C(\tau ^{2}),(\psi _{2m_{0}+1},\psi _{2m_{1}+1},\psi
_{2m_{0}+1}\tau ^{m_{0}},\psi _{2m_{1}+1}\tau ^{m_{1}})\mid m_{0},m_{1}\in 
\mathbb{Z}_{4}\right\rangle ,
\end{equation*}%
where, for each $\eta \in U(\mathbb{Z}_{4}),\;\psi _{\eta }$ is defined by (%
\ref{defpsi}) and 
\begin{equation*}
\tau ^{\psi _{\eta }}=\tau ^{\eta }.
\end{equation*}
\end{proposition}

\begin{proof} As 
\begin{equation} \label{eqforma}
\alpha =(\psi
_{2m_{0}+1}\psi _{2m_{1}+1},\psi _{2m_{0}+1}\tau ^{m_{0}},\psi
_{2m_{1}+1}\tau ^{ m_{1}}),
\end{equation}%
conjugates $\tau^{2} $ to 
\begin{equation*} 
(\tau ^{m_{0}},\tau ^{m_{1}},\tau ^{m_{0}+1},\tau
^{m_{1}+1})(0,2)(1,3),
\end{equation*}%
where $m_{0},m_{1}\in \mathbb{Z}_{4},$ 
and any other element in $N_{\mathcal{A}_{4}}(\dot{H})$
is equal to an element in $C(\tau^{2}) $ times an element of the form
(\ref{eqforma}),  then
$N_{\mathcal{A}_{4}}(\dot{H})$ is the desired subgroup.
\end{proof}

\begin{theorem}
\label{principaln4} Let $G$ be a   solvable subgroup of $%
Aut(T_{4})$ which contains $\tau $. Then, $G$ is a subgroup of 
\begin{equation}
\times _{4}\left( \cdots \left( \times _{4}\left( \times _{4}T^{\alpha }\rtimes S_{4}\right) \rtimes S_{4}\right) \cdots \right)
\rtimes S_{4}  \label{eqprincipaln4}
\end{equation}%
for some $\alpha \in \mathcal{A}_{4},$ where $T$ is the normalizer  in $\mathcal{A}_{4} $ of $ C(\tau^{2})$ .
\end{theorem}

\begin{proof}
As in the case $n=p$, we assume $G$ has derived length $d\geq 2$ and let $B$
be the $(d-1)$th term of the derived series of $G.$ Then, $B$ is an abelian
group normalized by $\tau $. On analyzing the case \ref{caso4_5} and the
final step, there exists a level $t$ such that $B$ is a subgroup of $\dot{V}%
=\times _{4^{k}}C(\mu ^{2}),$where $\mu =\tau ^{\alpha }$ for some $\alpha
\in \mathcal{A}_{4}$ and where $\sigma _{\mu ^{2}}=(0,2)(1,3)$. There also
exists $\beta \in B$ such that $\beta |_{u}=\mu ^{2}$ for some index $u\in 
\mathcal{M}.$

Moreover, if $T$ is the normalizer of $C(\tau ^{2}),$ then clearly, $%
T^{\alpha }$ is the normalizer of $C(\mu ^{2})$.

We will show now that $G$ is a subgroup of 
\begin{equation*}
\dot{J}=\times _{4}\left( \cdots \left( \times _{4}\left( \times _{4}T^{\alpha }\rtimes S_{4}\right) \rtimes S_{4}\right)
\cdots \right) \rtimes S_{4}
\end{equation*}%
where the cartesian product $\times _{4}$appears $t$ times..

Let $\gamma \not\in \dot{J}.$ Since $\gamma \not\in \dot{J},$ there exists $%
w\in \mathcal{M}$ having $|w|=t$ and $\gamma |_{w}\not\in T^{\alpha }$.
Since $\tau $ is transitive on all levels of the tree, by Corollary \ref%
{corolcentralizat2} we can conjugate $\beta $ by an appropriate power of $%
\tau $ to get $\theta \in B$ such that 
\begin{equation*}
\theta |_{w}=\mu ^{2}\text{ or }\theta |_{w}=\left( \mu ^{2}\right) ^{\tau
}=\left( (\tau ^{m_{0}},\tau ^{m_{1}},\tau ^{m_{0}+1},\tau
^{m_{1}+1})(0,2)(1,3)\right) ^{\alpha },
\end{equation*}%
where $m_{0},m_{1}\in \mathbb{Z}_{4}.$ Thus, for $v=w^{\gamma }$ we have%
\begin{equation*}
\left( \theta ^{\gamma }\right) |_{v}\overset{(\ref{eq9})}{=}\theta
|_{v^{\gamma ^{-1}}}^{\gamma _{v^{\gamma ^{-1}}}}=\theta |_{w}^{\gamma
|_{w}}\not\in C(\mu ^{2})
\end{equation*}%
which implies $\theta ^{\gamma }\not\in B\leq \dot{V}$ and $\gamma \not\in G$%
. Hence, $G$ is a subgroup of $\dot{J}$.
\end{proof}

%------------------------------------------------------------------------------------%

\vskip 0.4 true cm

\newpage

\begin{center}{\textbf{Acknowledgments}}
\end{center}
The authors are grateful to Gustavo A. Fern\'andez-Alcober for numerous observations on an earlier version of this paper. 
We also thank the referee for further helpful remarks. \\ \\
\vskip 0.4 true cm

%------------------------------------------------------------------------------------%

%-----------------------------------------------------------------------------
%-----------------------------------------------------------------------------

\bigskip
\bigskip

{\footnotesize \pn{\bf Josimar da Silva Rocha}\; \\ {Coordination of
Mathematics (COMAT)}, {Universidade Tecnol\'ogica Federal do Paran\'a, 86300-000,} {Corn\'elio Proc\'opio-PR, Brazil}\\
{\tt Email: jsrocha74@gmail.com}\\

{\footnotesize \pn{\bf Said Najati Sidki}\; \\ {Department of
Mathematics}, {Universidade de Bras\'{\i}lia,
70910-900,} {Bras\'{\i}lia-DF, Brazil}\\
{\tt Email: ssidki@gmail.com}\\
\end{document}